\RequirePackage{amsthm}

\documentclass[sn-mathphys,Numbered]{sn-jnl}

\usepackage{graphicx}%
\usepackage{multirow}%
\usepackage{amsmath,amssymb,amsfonts}%
\usepackage{mathrsfs}%
\usepackage[title]{appendix}%
\usepackage{xcolor}%
\usepackage{textcomp}%
\usepackage{manyfoot}%
\usepackage{booktabs}%
\usepackage{algorithm}%
\usepackage{algorithmicx}%
\usepackage{algpseudocode}%
\usepackage{listings}%
\usepackage[capitalise,noabbrev]{cleveref}

\usepackage{xr}
\usepackage{bookmark}
\usepackage{paralist,enumitem}
\usepackage{verbatim}
\usepackage[T1]{fontenc}
\usepackage[utf8]{inputenc}
\usepackage{amsfonts,amssymb,amsmath,graphicx,hyperref}
\usepackage{tikz,pgfplots}
\pgfplotsset{compat=1.16} 
\usepackage{standalone}
\usepackage{mathtools}
\usepackage[caption=false]{subfig}
\usepackage{tabularx}
\usepackage{enumitem}
\usepackage{bm}
\usepackage{booktabs}
\usepackage{tcolorbox}
\usepackage{enumerate}
\usepackage{algpseudocode}
\usepackage{tikz-3dplot}

\newcommand{\Rb}{\mathbb{R}}
\newcommand{\Cb}{\mathbb{C}}

\newcommand{\dsum}{\displaystyle \sum}

\newcommand{\dcup}{\displaystyle \cup}

\newcommand{\dbcup}{\displaystyle \bigcup}

\newcommand{\abs}[1]{\left\lvert #1 \right\rvert}
\newcommand{\magn}[1]{\left\lVert #1 \right\rVert}
\newcommand{\ceil}[1]{\left\lceil #1 \right\rceil}

\usepackage{amsopn}

\newcommand{\iso}{0.4}

\newcommand{\cube}[2]{
	\draw [ultra thick] (#1,#2) rectangle (1+#1,1+#2);
	\draw [ultra thick, loosely dotted] (\iso+#1,\iso+#2) rectangle (1+\iso+#1,1+\iso+#2);
	\draw [ultra thick] (1+#1,1+#2) -- (1+\iso+#1,1+\iso+#2);
	\draw [ultra thick, loosely dotted] (1+#1,#2) -- (1+\iso+#1,\iso+#2);
	\draw [ultra thick, loosely dotted] (#1,1+#2) -- (\iso+#1,1+\iso+#2);
	\draw [ultra thick, loosely dotted] (#1,#2) -- (\iso+#1,\iso+#2);
}

\newcommand{\textib}[1]{\textbf{\textit{\text{#1}}}}
\newcommand{\Bb}{\textbf{\textit{\text{B}}}}
\newcommand{\Cl}{{\mathcal{C}^{(l)}_i}}

\newcommand{\bkt}[1]{
	\left(#1\right)
}

\DeclareRobustCommand{\cubex}[6]{
	\draw[#5] (#1,#2,#3) -- (#1,#2+#4,#3) -- (#1+#4,#2+#4,#3) -- (#1+#4,#2,#3) -- cycle;
	\draw[#5] (#1,#2,#3+#4) -- (#1,#2+#4,#3+#4) -- (#1+#4,#2+#4,#3+#4) -- (#1+#4,#2,#3+#4) -- cycle;
	\draw[#5] (#1,#2,#3) -- (#1,#2,#3+#4);
	\draw[#5] (#1,#2+#4,#3) -- (#1,#2+#4,#3+#4);
	\draw[#5] (#1+#4,#2,#3) -- (#1+#4,#2,#3+#4);
	\draw[#5] (#1+#4,#2+#4,#3) -- (#1+#4,#2+#4,#3+#4);
    \node at (#1 + 0.5*#4,#2 + 0.5*#4,#3 + 0.5*#4) {#6};
}
\newcommand{\cmt}[1]{\iffalse {#1} \fi}

\theoremstyle{thmstyleone}%
\newtheorem{remark}{Remark}%
\newtheorem{lemma}{Lemma}%
\newtheorem{corollary}{Corollary}%
\newtheorem{definition}{Definition}%
\newtheorem{theorem}{Theorem}

\begin{document}

\title[HODLR3D]{HODLR3D: Hierarchical matrices for $N$-body problems in three dimensions}

\author*[1]{\fnm{Kandappan} \sur{V A}}\email{kandappanva@gmail.com}
\equalcont{These authors contributed equally to this work.}

\author[1]{\fnm{Vaishnavi} \sur{Gujjula}}\email{vaishnavihp@gmail.com}
\equalcont{These authors contributed equally to this work.}

\author[1,2]{\fnm{Sivaram} \sur{Ambikasaran}}\email{sivaambi@iitm.ac.in}

\affil*[1]{\orgdiv{Department of Mathematics}, \orgname{Indian Institute of Technolgy Madras}, \orgaddress{\city{Chennai}, \postcode{600036}, \state{Tamil Nadu}, \country{India}}}

\affil[2]{\orgdiv{Robert Bosch Centre for Data Science and Artificial Intelligence}, \orgname{Indian Institute of Technolgy Madras}, \orgaddress{\city{Chennai}, \postcode{600036}, \state{Tamil Nadu}, \country{India}}}

\abstract{This article introduces HODLR3D, a class of hierarchical matrices arising out of $N$-body problems in three dimensions. HODLR3D relies on the fact that certain off-diagonal matrix sub-blocks arising out of the $N$-body problems in three dimensions are numerically low-rank. For the Laplace kernel in $3$D, which is widely encountered, we prove that all the off-diagonal matrix sub-blocks are rank deficient in finite precision. We also obtain the growth of the rank as a function of the size of these matrix sub-blocks. For other kernels in three dimensions, we numerically illustrate a similar scaling in rank for the different off-diagonal sub-blocks. We leverage this hierarchical low-rank structure to construct HODLR3D representation, with which we accelerate matrix-vector products. The storage and computational complexity of the HODLR3D matrix-vector product scales almost linearly with system size. We demonstrate the computational performance of HODLR3D representation through various numerical experiments. Further, we explore the performance of the HODLR3D representation on distributed memory systems. \emph{HODLR3D, described in this article, is based on a weak admissibility condition. Among the hierarchical matrices with different weak admissibility conditions in $3$D, only in HODLR3D did the rank of the admissible off-diagonal blocks not scale with any power of the system size. Thus, the storage and the computational complexity of the HODLR3D matrix-vector product remain tractable for $N$-body problems with large system sizes.}}



\maketitle

\section{Introduction}\label{sec:introduction}
Hierarchical matrix representations are used to accelerate $N$-body problems arising from applications such as particle simulations, high-order statistics, machine learning~\cite{gray2000n,litvinenko2019likelihood}, solving PDEs, radial basis function interpolation~\cite{coulier2016efficient,gumerov2007fast}, etc. In this article, we introduce a class of hierarchical matrix representation~\cite{hackbusch1999sparse,grasedyck2003construction,H2DKandappan} for matrices arising out of $N$ body problems in three dimensions to construct an algorithm that performs matrix-vector products that has almost linear computational cost~\footnote{We say a matrix algorithm has almost linear computational complexity if given $A \in \Cb^{N \times N}$, the computational cost of the algorithm scales as $\mathcal{O}\bkt{N^{1+\epsilon}}$ for all $\epsilon>0$.} both in time and space. One of the first almost linear complexity algorithms for $N$-body problems was introduced by Barnes and Hut~\cite{barnes1986hierarchical} (frequently addressed as the "Treecode"), which reduced the computational complexity of matrix-vector product from $\mathcal{O}\bkt{N^2}$ to $\mathcal{O}\bkt{N \log N}$. The Fast Multipole Method by Greengard and Rokhlin~\cite{greengard1987fast,greengard1988rapid,greengard1997new} further reduced the computational complexity of these $N$-body problems to $\mathcal{O}\bkt{N}$.

Algebraically, the Treecode and the FMM can be interpreted as the matrix corresponding to $N$-body problems having a hierarchical low-rank structure. This notion was pioneered by Hackbush~\cite{hackbusch1999sparse,grasedyck2003construction} in the late 1990's and early 2000's and these matrices were termed as hierarchical matrices ($\mathcal{H}$-matrix). These hierarchical matrices use a hierarchical tree to subdivide the matrix and identify the low-rank matrix sub-blocks at different levels in the hierarchical tree. A detailed description of these matrices is in~\Cref{sec:hmat}. These hierarchical matrices can be stored efficiently and matrix algorithms for these hierarchical matrices can be devised in almost linear computational cost.

There has been extensive research on hierarchical matrices in recent years. Some of the widely used hierarchical representations are Hierarchically Off-Diagonal Low-Rank (HODLR)~\cite{sa_thesis,SA_FDS_2013}, Hierarchically Semi-Separable (HSS)~\cite{chandrasekaran2005some,vandebril2005bibliography,vandebril2005note}, $\mathcal{H}^{2}$\cite{borm2003hierarchical,borm2003introduction,hackbusch2015hierarchical}, etc.   For a more detailed literature review on hierarchical matrices and their applications, we refer the readers to the articles~\cite{borm2003introduction,yokota2015fast,H2DKandappan} and the references therein. In addition to the hierarchical matrices, there are flat low-rank structures such as Block Low Rank (BLR) format, where the matrix is subdivided into $n_b\times n_b$ blocks with each block of size at most $b\times b$. Like, hierarchical representations, certain off-diagonal blocks in the BLR format are approximated and represented as low-rank matrices~\cite{amestoy2015improving,amestoy2017complexity}.

The particular hierarchical low-rank matrix which is of interest to this article is HODLR matrix representation of a dense matrix from $N$ body problems. In the HODLR matrix representation, all the off-diagonal blocks that result from the hierarchical subdivision of the underlying domain are approximated by low-rank matrices. Leveraging the HODLR matrix representation, we can construct an algorithm that performs matrix-vector in $\mathcal{O}(pN\log(N))$ where $p$ is the maximum rank of the off-diagonal sub-blocks~\cite{SA_FDS_2013}. The maximum rank of the off-diagonal blocks plays a critical role in the computational complexity of the HODLR matrix. For example, if we represent the dense matrix from a wide range of $N$ body problems in $1$D using HODLR representation, the off-diagonal blocks correspond to the vertex-sharing interactions. We know from~\cite{khan2022numerical,sa_thesis} that the ranks of the off-diagonal blocks due to the vertex sharing interaction scale as $\mathcal{O}\bkt{\log\bkt{N}}$. Using that fact we have the computational complexity of the matrix-vector product in HODLR representation scale as $\mathcal{O}\bkt{N\log^2\bkt{N}}$, i.e., it scales almost linearly. Whereas in the case of $N$ body problems from $2$D, the maximum rank of the off-diagonal blocks is due to the edge-sharing interactions whose rank scales roughly as $\mathcal{O}\bkt{\sqrt{N}}$~\cite{khan2022numerical,H2DKandappan}. So, the computational complexity of the matrix-vector product in HODLR representation scale roughly as $\mathcal{O}\bkt{N^{3/2}}$. In the case of the matrices from three-dimensional problems, the maximum rank of the off-diagonal blocks is due to the face-sharing interactions whose rank scales roughly as $\mathcal{O}\bkt{\sqrt[3]{N^2}}$ (refer \cref{tm:ranks3D} and~\cite{khan2022numerical} for more details), which results in the computational complexity for matrix-vector product using HODLR representation as roughly $\mathcal{O}\bkt{N^{5/3}}$.  Accordingly, the matrix-vector product through HODLR matrix representation does not scale almost linearly.

In our recent work~\cite{H2DKandappan}, we had proposed a new hierarchical matrix to represent matrices arising out of $N$-body problems in two dimensions which performs matrix-vector product in almost linear complexity. In this article, we extend this notion to $N$ body problems in three dimensions and introduce a new class of hierarchical matrices, HODLR3D. HODLR2D~\cite{H2DKandappan} and HODLR3D representation are the extensions of HODLR matrix representation in two and three dimensions by a judicious choice of admissibility condition. The complexity of matrix-vector product using HODLR3D representation scales as $\mathcal{O}(pN\log(N))$, where $p$ is the maximum rank of the compressed off-diagonal blocks. For HODLR3D, the maximum rank of the compressed off-diagonal blocks is due to vertex-sharing interactions. Through our numerical illustrations and~\cref{tm:ranks3D} we show that the vertex-sharing interaction scale as $\mathcal{O}\bkt{\log^{3}\bkt{N}}$. Hence, like HODLR2D, HODLR3D also scales almost linear in complexity for the matrix-vector product.

The main highlights of the article are:
\begin{itemize}
    \item 
    We illustrate the scaling of ranks of different off-diagonal blocks for handpicked kernel functions numerically. We state~\Cref{tm:ranks3D}, which provides bounds on the ranks of different off-diagonal blocks for the Laplacian kernel in $3$D. The proof can be found in~\Cref{sec:Proof}. We refer the readers to~\cite{khan2022numerical} for the bounds on the ranks of different off-diagonal blocks for a bigger class of kernel functions.
    \item
    Based on the outcomes of~\Cref{tm:ranks3D} and the numerical illustrations, we build the HODLR3D representation, wherein we only choose to compress the off-diagonal blocks corresponding to interactions between vertex sharing neighbours and well-separated boxes of the hierarchical tree to achieve an almost linear complexity algorithm.
    \item 
    We perform various numerical experiments to illustrate the performance of the HODLR3D matrix-vector product. 
\end{itemize}

\subsection{Main Result}
We state the main theorem that gives the scaling of ranks of the different off-diagonal blocks of the matrix arising out of Green's function of the 3D Laplace equation. The proof of~\Cref{tm:ranks3D} can be found in~\Cref{sec:Proof}.~\cref{tm:ranks3D} can be extended to a wide range of kernel functions, the proof of which can be found in~\cite{khan2022numerical}. On comparing with~\cite{khan2022numerical}, we get a slightly tighter bound in~\cref{tm:ranks3D} for the kernel $1/r$. This should not be surprising since~\cite{khan2022numerical} is applicable for a wide range of kernel functions as opposed to~\Cref{tm:ranks3D}, which is applicable only for $1/r$ kernel function.
\begin{tcolorbox}[colback=white!5!white,colframe=black!75!black,title=\textbf{Rank of different interactions for the Laplace kernel in $3$D}]
        \begin{theorem} \label{tm:ranks3D}
            Consider a box $\Bb\subseteq\mathbb{R}^3$. The boxes $\Bb_1\subseteq\Bb$ and $\Bb_2\subseteq\Bb$ are from the hierarchical subdivision of the box $\Bb$ using a balanced octree. Let $\{q_j\}_{j=1}^N$ be $N$ uniformly located charges at $\{z_j\}_{j=1}^N$ inside a box $\Bb_1$ of side length $l$. Let $Q = \dsum_{i=1}^N \abs{q_i}$. The potential due to these $N$ charges at $M$ locations $\{w_i\}_{i=1}^M$, with $M \geq N$, inside another box $\Bb_2$ is given by $\phi_i = \dsum_{j=1}^N \dfrac{q_j}{\abs{w_j-z_j}}$. In matrix-vector parlance, we have
            $$\vec{\phi}=A\vec{q}$$
            where $\vec{q} \in \Rb^{N \times 1}$, $\vec{\phi} \in \Rb^{M \times 1}$ and $A \in \Rb^{M \times N}$ with $A_{ij} = \dfrac{1}{\abs{w_j-z_j}}$. Then for a given $\epsilon > 0$, there exists a matrix $\tilde{A} \in \Rb^{M \times N}$ with rank at the most $\tau$ that approximates $\phi_i$ such that $\abs{\phi_i-\bkt{\tilde{A}q}_i} < \epsilon$ and $\tau \in \mathcal{O} \bkt{R(N)\log^2\bkt{\dfrac{S(N)Q}{l\epsilon}}}$, where
            \begin{itemize}
                \item $R(N) = 1,S(N)=1$ if the boxes $\Bb_1$ and $\Bb_2$ are at least one box away, i.e., $\text{dist}\bkt{\Bb_1,\Bb_2} \geq l$, where $dist(x,y)$ is the Euclidean distance between the sets $x$ and $y$.
                \item $R(N) = \log(N),S(N)=\sqrt[3]{N}$ if the boxes $\Bb_1$ and $\Bb_2$ are vertex sharing neighbors, i.e., $\Bb_1\cap\Bb_2$ is at the most a vertex.
                \item $R(N) = \sqrt[3]{N}, S(N)=\sqrt[3]{N^2}$ if the boxes $\Bb_1$ and $\Bb_2$ are edge sharing neighbors, i.e., $\Bb_1\cap\Bb_2$ is at the most a line.
                \item $R(N) = \sqrt[3]{N^2}, S(N)=N$ if the boxes $\Bb_1$ and $\Bb_2$ are face sharing neighbors, i.e., $\Bb_1\cap\Bb_2$ is at most a $2$D plane.
            \end{itemize}
        \end{theorem}
    \end{tcolorbox}

\subsection{Outline of the article}
The rest of the article is organized as follows: In~\Cref{sec:OctTree}, we describe the octree on which we build the hierarchical matrices and also give a brief on HODLR and a $\mathcal{H}$ matrix structures. Further in~\Cref{sec:rankGrowth}, we numerically illustrate the growth of ranks for various kinds of interactions and prove~\cref{tm:ranks3D}. And in~\Cref{sec:HODLR3D}, we describe the algorithm for HODLR3D matrix-vector product. We present the numerical results in~\Cref{sec:numericalResults} and compare the performance of HODLR3D with HODLR and a $\mathcal{H}$ matrix representation. 
We describe the parallel implementation of HODLR3D and present numerical results in~\Cref{sec:parH3D}.

\section{Preliminaries}\label{sec:OctTree}

In this section, we describe the construction of the tree that is used to build the HODLR3D matrix. We briefly describe the HODLR and $\mathcal{H}$ matrix representations, with which we compare the HODLR3D later in the article. 
\subsection{Construction of tree}
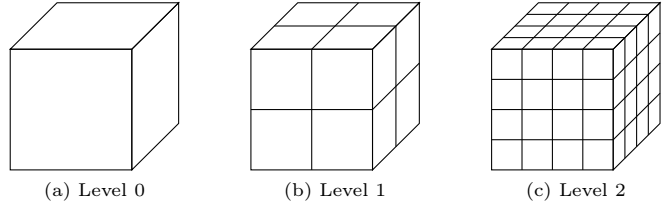
\begin{figure}[htbp]
    \centering
    \subfloat[Level 0]{
        \begin{tikzpicture}[scale=0.4]
            \draw (0,0,4)--(4,0,4)--(4,4,4)--(0,4,4)--cycle;
            \draw (4,0,0) -- (4,4,0)--(4,4,4)--(4,0,4)--cycle;
            \draw (4,4,0) -- (4,4,4)--(0,4,4)--(0,4,0)--cycle;
            \end{tikzpicture}
        }\qquad%
    \subfloat[Level 1]{
        \begin{tikzpicture}[scale=0.4]
            \draw (0,0,4)--(4,0,4)--(4,4,4)--(0,4,4)--cycle;
            \draw (4,0,0) -- (4,4,0)--(4,4,4)--(4,0,4)--cycle;
            \draw (4,4,0) -- (4,4,4)--(0,4,4)--(0,4,0)--cycle;
            \foreach \x in {2}
                {
            \draw(4,0,\x)--(4,4,\x); 
            \draw(4,\x,0)--(4,\x,4); 
            \draw(0,\x,4)--(4,\x,4); 
            \draw(\x,0,4)--(\x,4,4); 
            \draw(\x,4,0)--(\x,4,4); 
            \draw(0,4,\x)--(4,4,\x); 
            }
            \end{tikzpicture}
        }\qquad%
    \subfloat[Level 2]{
        \begin{tikzpicture}[scale=0.4]
            \draw (0,0,4)--(4,0,4)--(4,4,4)--(0,4,4)--cycle;
            \draw (4,0,0) -- (4,4,0)--(4,4,4)--(4,0,4)--cycle;
            \draw (4,4,0) -- (4,4,4)--(0,4,4)--(0,4,0)--cycle;
            \foreach \x in {1,2,3}
                {
            \draw(4,0,\x)--(4,4,\x); 
            \draw(4,\x,0)--(4,\x,4); 
            \draw(0,\x,4)--(4,\x,4); 
            \draw(\x,0,4)--(\x,4,4); 
            \draw(\x,4,0)--(\x,4,4); 
            \draw(0,4,\x)--(4,4,\x); 
            }
            \end{tikzpicture}
        }
\caption{Hierarchical sub-division of $\textbf{\textit{B}}$ using an octree} 
\label{fig:octTree}
\end{figure}

We assume the computational domain to be a cube $\textbf{\textit{B}} \subset \mathbb{R}^3$, that has the support of the particles. Level $0$ of the tree is the domain $\textbf{\textit{B}}$ itself. A cube at level $k$ is subdivided into $8$ cubes that belong to level $k+1$ of the tree. The former is said to be the parent of the latter, and the latter are said to be the children of the former. The sub-division is carried on until the depth of the tree is $L$, where a cube at level $L$, called a leaf, contains at most $N_{max}$ particles. $N_{max}$ is a user-specified threshold that defines the maximum number of particles in a leaf. We denote the octree by $\mathcal{T}^{L}$. We illustrate the construction of octree till level $2$ in Figure~\ref{fig:octTree}. In this article, we work with the balanced tree for simplicity, though an adaptive oct-tree or a K-d tree too can be considered.

\subsection{\texorpdfstring{$\mathcal{H}$}{Hc}-matrix}
\label{sec:hmat}
\begin{figure}[htbp]
    \subfloat[Level 1]{
        \includegraphics[scale=0.35]{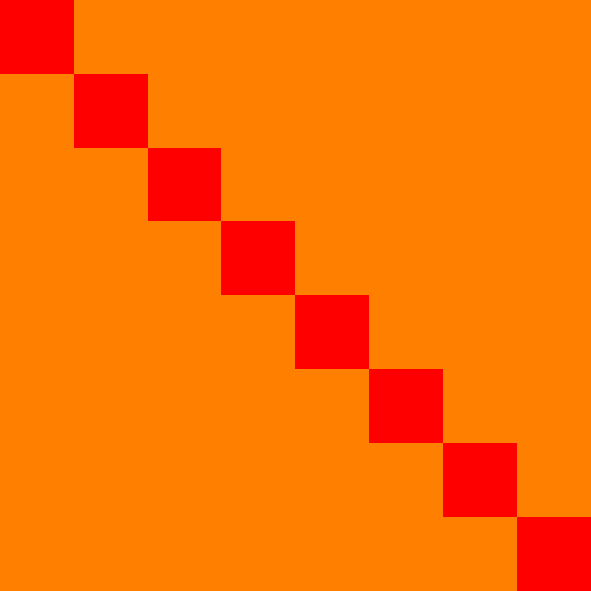}
        }\quad%
    \subfloat[Level 2]{
        \includegraphics[scale=0.35]{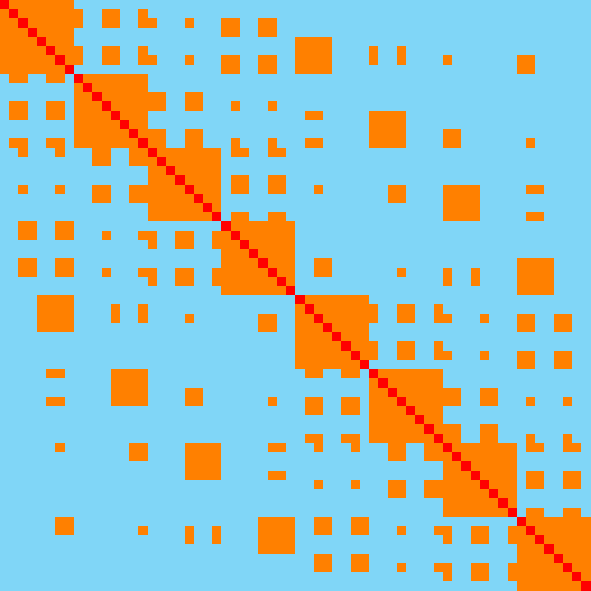}
        }\quad%
    \subfloat{
        \begin{tikzpicture}
            [
            box/.style={rectangle,draw=black, minimum size=0.25cm},scale=0.2
            ]
            \node[box,fill=red,label=right:{\scriptsize Full-rank (self-interaction)}, anchor=west] at (0,12){};
            \node[box,fill=orange,label=right:{\scriptsize Full-rank (neighbour-interaction)}, anchor=west] at (0,8){};
            \node[box,fill=cyan!50,label=right:{\scriptsize Low-rank}, anchor=west] at (0,4){};
        \end{tikzpicture}
    }
    \caption{A $\mathcal{H}$ matrix low-rank structure at levels 1 and 2}
    \label{fig:H_matrix}
\end{figure}
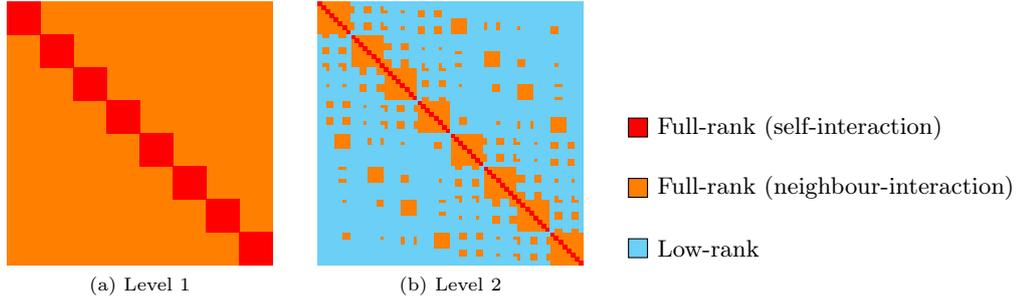

$\mathcal{H}$-matrix is the hierarchical low-rank representation of a particular class of dense matrices based on an admissible condition. The admissibility condition helps in identifying the off-diagonal blocks due to the interaction between blocks in a particular level to be represented as low-rank matrices with considerable accuracy. In literature, there exist two classes of hierarchical matrix representations based on the admissibility condition, viz., $\mathcal{H}$ matrix based on strong or standard admissibility condition and $\mathcal{H}$ based on weak admissibility condition~\cite{hackbusch2004hierarchical}. In the case of $\mathcal{H}$-matrix with weak admissibility condition, the admissible blocks are closer, whereas, in the case of the strong or standard admissible condition, the admissible blocks are well-separated. In this article, we consider the $\mathcal{H}$-matrix to be built on the hierarchical tree, $\mathcal{T}^{L}$ with the strong or standard admissibility condition given by~\Cref{eq:adm}. Consider cells $X$ and $Y$ that are at the same level of $\mathcal{T}^{L}$. Let the index sets of the particles lying in cells $X$ and $Y$ be $I_{X}$ and $I_{Y}$. The matrix sub-block $K(I_{X}, I_{Y})$, represented using MATLAB notation, is approximated by a low-rank matrix if the following admissibility condition is satisfied.
\begin{gather} 
\max(diam(X), diam(Y)) \leq \eta \quad dist(X,Y),\quad\text{where,}\label{eq:adm}\\
diam(X)= \sup\{\|x-y\|_{2}:x,y\in X\},\label{eq:diam}\\
dist(X,Y) = \inf\{\|x-y\|_{2}:x\in X, y\in Y\},\label{eq:dist}
\end{gather}
    and $\eta$ is a parameter that controls the admissibility condition. For $\eta=\sqrt{3}$, the interaction between particles in cells that are well-separated, i.e., those cells that do not share a boundary, is approximated by a low-rank matrix. In~\cref{fig:H_matrix}, we illustrate the low-rank structure of a $\mathcal{H}$ matrix, with $\eta=\sqrt{3}$, at levels 1 and 2.
    
\subsection{HODLR matrix}
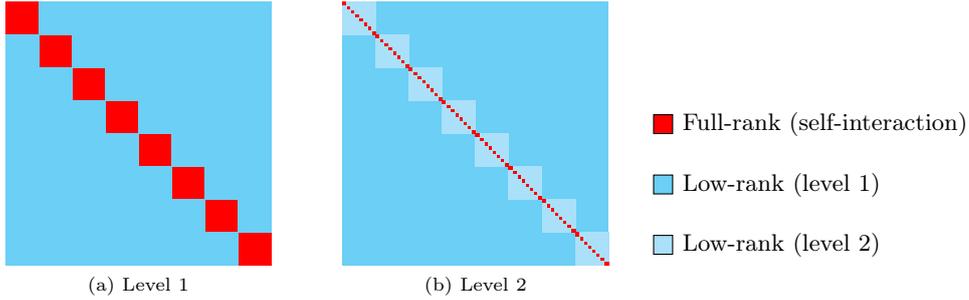
\begin{figure}[!htbp]
    \centering
    \subfloat[Level 1]{
        \includegraphics[scale=0.35]{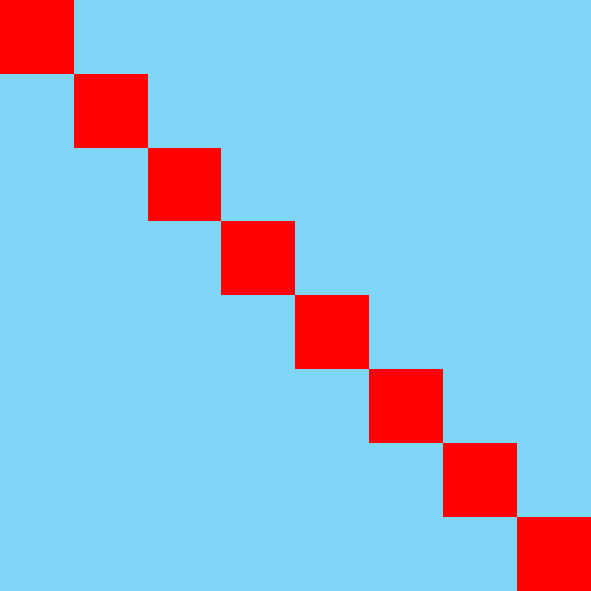}
        }\qquad%
    \subfloat[Level 2]{
        \includegraphics[scale=0.35]{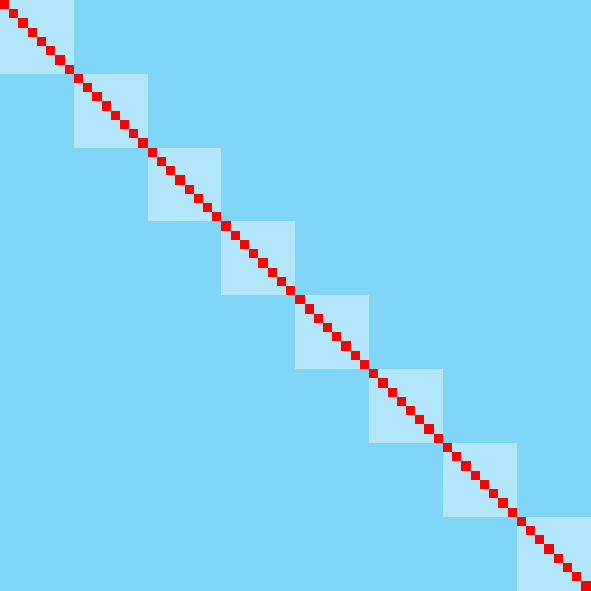}
        }\quad%
    \subfloat{
        \begin{tikzpicture}
            [
            box/.style={rectangle,draw=black, minimum size=0.25cm},scale=0.2
            ]
            \node[box,fill=red,label=right:{\scriptsize Full-rank (self-interaction)}, anchor=west] at (0,12){};
            \node[box,fill=cyan!50,label=right:{\scriptsize Low-rank (level 1)}, anchor=west] at (0,8){};
            \node[box,fill=cyan!30,label=right:{\scriptsize Low-rank (level 2)}, anchor=west] at (0,4){};
        \end{tikzpicture}
    }
    \caption{HODLR matrix low-rank structure at levels 1 and 2}
    \label{fig:HODLR_matrix}
\end{figure}

HODLR matrix is an example of a $\mathcal{H}$ matrix with weak admissibility condition. It is constructed by compressing all the off-diagonal sub-blocks~\cite{SA_FDS_2013}. Equivalently, other than the self-interaction, all other interactions are approximated by a low-rank matrix.
Usually, HODLR is built on a K-d tree, but in this article, since we develop HODLR3D on an octree, to be on the same terms, we consider HODLR to be built on an octree. In Figure~\ref{fig:HODLR_matrix}, we illustrate the low-rank structure of the HODLR matrix at levels 1 and 2.
    
\section{Rank growth of different interactions in 3D}\label{sec:rankGrowth}
In this section, we analyze numerically the ranks of sub-matrices corresponding to different types of interaction in $3$D for handpicked widely used kernel functions. We then prove the growth in ranks for various off-diagonal blocks of the matrix for Green's function of the Laplace equation in $3$D, which is $\dfrac1{r}$ kernel. 

\subsection{Numerical rank for different kernels}
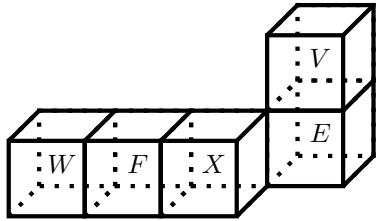
\begin{figure}[!htbp]
\begin{center}
 \begin{tikzpicture}
    \cube{0}{0}
    \cube{-1}{0}
    \cube{-2}{0}
    \cube{1+\iso}{\iso}
    \cube{1+\iso}{1+\iso}
    \draw [ultra thick] (-2,1) -- (-2+\iso,1+\iso) -- (1+\iso,1+\iso);
    \draw [ultra thick] (1,0) -- (1+\iso,\iso);
    \draw [ultra thick] (2+\iso,\iso) -- (2+2*\iso,2*\iso) -- (2+2*\iso,2+2*\iso) -- (1+2*\iso,2+2*\iso) -- (1+\iso,2+\iso);

    \node at (0.7,0.7) {$X$};
    \node at (-0.3,0.7) {$F$};
    \node at (-1.3,0.7) {$W$};
    \node at (1.7+\iso,0.7+\iso) {$E$};
    \node at (1.7+\iso,1.7+\iso) {$V$};
 \end{tikzpicture}
 \end{center}
 \caption{An illustration of cubes $V$, $E$, $F$, that share a vertex, edge, face with cube $X$ respectively and cube $W$ that is well-separated to $X$.} \label{fig:interactions}
\end{figure}

\begin{figure}[!htbp]
      \begin{center}
    \includegraphics[width=1\linewidth]{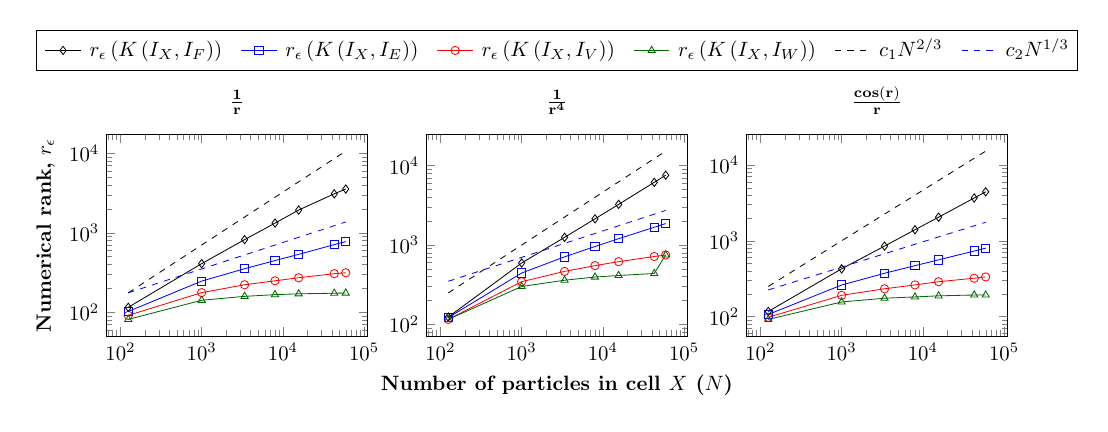}
        \end{center}
        \caption{Numerical ranks with $\epsilon=10^{-14}$ of different kinds of interactions for kernels a) $\frac{1}{r}$ b) $\frac{1}{r^{4}}$ c) $\frac{\cos(r)}{r}$ versus $N$.} 
        \label{NumericalRanks}
\end{figure}

We consider a cube $X$ belonging to the octree. And we further consider cubes $V$, $E$, and $F$ that share a vertex, edge, and face with cube $X$ and a cube $W$ that is well-separated to $X$ as shown in~\cref{fig:interactions}.
We uniformly distribute $N$ particles in each cube $W$, $F$, $X$, $E$, and $V$ at $\{r_{i}\}_{i=1}^{5N}$. We consider the interaction between two particles at $r_{i}$ and $r_{j}$ to be
\begin{equation}
    K\bkt{i,j} = \begin{cases}
    0, & \text{if} \quad i=j\\
   {f\bkt{\magn{r_{i}-r_{j}}_{2}}}, & \text{otherwise}
\end{cases}\label{eq:Laplace3D}
\end{equation}
where $f(r):\mathbb{R}\rightarrow\mathbb{R}$. Let the index sets of the particles lying in each of these cubes $W$, $F$, $X$, $E$, and $V$ be $I_{W}$, $I_{F}$, $I_{X}$, $I_{E}$, and $I_{V}$ respectively. The matrices $K(I_{X}, I_{W})$, $K(I_{X}, I_{F})$, $K(I_{X}, I_{E})$, and $K(I_{X}, I_{V})$ correspond to different off-diagonal sub-matrices. We illustrate in Figure~\ref{NumericalRanks}, the growth of numerical ranks\footnote{For an $\epsilon>0$, the numerical rank of matrix $K$, $r_{\epsilon}(K)$ is defined as $r_{\epsilon}(K) = \max\{k\in\{1, 2,\dots N\}:\frac{\sigma_{k}}{\sigma_{1}}>\epsilon\}$ where $\sigma_{1}\geq \sigma_{2}\geq \dots \sigma_{N}$ are the singular values of $K$.} of the off-diagonal sub-matrices from four kinds of interactions for the three kernel functions $\bm{K}(x,y) = f(r)$ where, $r = \magn{x-y}_{2}$, defined as i) $\frac{1}{r}$ ii) $\frac{1}{r^4}$ iii) $\frac{\cos(r)}{r}$. It is to be observed that the ranks of the face and edge-sharing interactions $K(I_{X}, I_{F}) \text{ and } K(I_{X}, I_{E})$, for the three kernels, scale roughly as $\mathcal{O}(N^{2/3})$ and $\mathcal{O}(N^{1/3})$ respectively, and that of the vertex sharing interaction $K(I_{X}, I_{V})$ scales roughly as $\mathcal{O}(\log^{3}(N))$ and the well-separated interaction $K(I_{X}, I_{W})$ does not scale with $N$.
We tabulate this observation in~\cref{tab:rankgrowth}.
\begin{table}[!htbp]
\centering
\caption{Numerical observation of the scaling of ranks of different off-diagonal blocks of size $N$}
\label{tab:rankgrowth}
\begin{tabular}{@{}ll@{}}
\toprule
Off-diagonal block & Rough scaling of numerical rank      \\ \midrule
$K(I_X,I_W)$        & $\mathcal{O}\bkt{1}$ \\
$K(I_X,I_V)$        & $\mathcal{O}\bkt{\log^{3}\bkt{N}}$ \\
$K(I_X,I_E)$        & $\mathcal{O}\bkt{\sqrt[3]{N}}$ \\
$K(I_X,I_F)$        & $\mathcal{O}\bkt{\sqrt[3]{N^2}}$ \\ \bottomrule
\end{tabular}
\end{table}
We have numerically shown the growth of ranks of various off-diagonal sub-blocks for widely used kernel functions in $3$D, or equivalently we have numerically established that the matrix possesses a low-rank structure which can be leveraged to construct fast matrix algorithms.

\subsection{Rank of off-diagonal blocks for Laplacian Kernel in \texorpdfstring{$3$}{3}D}
\label{sec:Proof}
In the previous subsection, we showed the growth of ranks of various off-diagonal blocks numerically for the Laplacian Kernel in $3$D. We will prove in~\Cref{tm:ranks3D}, the scaling of ranks for various off-diagonal blocks due to different interactions resulting from the hierarchical subdivision. 
\textbf{Proof of \cref{tm:ranks3D}.}
The proof of~\cref{tm:ranks3D} uses multipole expansions of the Laplacian kernel in $3$D. So we state here the multipole expansion lemma. The multipole expansions lemma approximates the potential at a point P due to $N$ charges using~\eqref{eq:multipole}. The proof of the lemma and its detailed analysis can be found in~\cite{beatson1997short,greengard1988rapid}. 
The potential at $M$ locations due to the $N$ charges can be represented as a matrix-vector product. By multipole expansions lemma, we can approximate the potential at $M$ locations through a low-rank representation as stated in~\cref{tm:prank}.
\begin{lemma}[Multipole Expansion]
    \label{tm:multipole}
        Suppose that $N$ charges of strengths $\{q_i\}_{i=1}^N$ are located at $\{\vec{z}_i\}_{i=1}^N$ whose spherical coordinates are $\{\bkt{\rho_i,\alpha_i,\beta_i}\}_{i=1}^N$ inside a sphere of radius $a$ and we have $\rho_i<a$. Then for any point $\vec{P} \in \Rb^3$ with the spherical coordinate $(r,\theta,\gamma)$ outside the sphere of radius $a$, the potential at $P$ due to the $1/\tilde{r}$ kernel (i.e., $\tilde{r}_i =\abs{\vec{P}-\vec{z}_i} $) is given by
        \begin{gather}
        \Phi\bkt{P} = \dsum_{n=0}^{\infty} \dsum_{m=-n}^{n} \dfrac{M_n^m}{r^{n+1}} Y_n^{m}\bkt{\theta,\gamma}, \quad \text{where,}\\
        M_n^m = \dsum_{i=1}^N q_i \rho_i^n Y_n^{-m}\bkt{\alpha_i,\beta_i},
        \end{gather}
        and $Y_n^{m}$ are the spherical harmonics. Furthermore, for any $p \geq 1$, we have
        \begin{equation}\label{eq:multipole}
            \abs{\Phi(P)-\dsum_{n=0}^p \dsum_{m=-n}^n \dfrac{M_n^{m}}{r^{n+1}}Y_n^{m}\bkt{\theta,\gamma}} \leq \dfrac{Q}{r-a}\bkt{\dfrac{a}r}^{p+1} = \dfrac{Q}a \bkt{\dfrac1{c-1}}\left(\dfrac1c\right)^{p+1}
        \end{equation}
        where $Q = \dsum_{i=1}^N \abs{q_i}$, $c=\dfrac{r}{a}$.
\end{lemma}
\begin{corollary}[Potential through rank-$\tau$ approximation]
 \label{tm:prank}
 From~\cref{tm:multipole}, the potential at $P\bkt{r_i,\theta_i,\gamma_i}$ using multipole expansions can be written in the matrix-vector format as follows:   
 \begin{equation}
     \phi_i = \dsum_{j=1}^N \dfrac{q_j}{\tilde{r}_{ij}} \approx \tilde{\phi}_i =
     \begin{bmatrix} 
         \dfrac{1}{r_i} & \dfrac{1}{r_i^2} & \dfrac{1}{r_i^3} & \cdots &  \dfrac{1}{r_i^{p+1}}
     \end{bmatrix}
     \begin{bmatrix}  
         \bkt{M_0^0}Y_0^0(\theta_i,\gamma_i) \\ 
         \dsum_{m=-1}^1 \bkt{M_1^m}Y_1^m(\theta_i,\gamma_i) \\ 
         \dsum_{m=-2}^2 \bkt{M_2^m}Y_2^m(\theta_i,\gamma_i)\\ 
         \vdots \\ 
         \dsum_{m=-p}^p \bkt{M_p^m}Y_p^m(\theta_i,\gamma_i)
     \end{bmatrix},
 \end{equation}     
 which can be rewritten as, 
  \begin{gather}
     \tilde{\phi}_i =  u_i V \vec{q},\quad\text{where,}\\
\begin{aligned}
u_i &=
\left[\begin{matrix}
\dfrac{Y_0^0(\theta_i,\gamma_i)}{r} & \dfrac{Y_1^{-1}(\theta_i,\gamma_i)}{r^2} & \dfrac{Y_1^0(\theta_i,\gamma_i)}{r^2} & \dfrac{Y_1^1(\theta_i,\gamma_i)}{r^2}&  \cdots
\end{matrix}\right.\\
&\qquad\qquad
\left.\begin{matrix}
\dfrac{Y_p^{-p}(\theta_i,\gamma_i)}{r^{p+1}}  & \cdots & \dfrac{Y_p^0(\theta_i,\gamma_i)}{r^{p+1}}   &  \cdots & \dfrac{Y_p^p(\theta_i,\gamma_i)}{r^{p+1}}
\end{matrix}\right]
\end{aligned}\label{theta-alt}\\    
     V = \begin{bmatrix}  
     Y_0^0(\alpha_1,\beta_1) & Y_0^0(\alpha_2,\beta_2) & \cdots & Y_0^0(\alpha_N,\beta_N) \\
     \rho_1 Y_1^{1}(\alpha_1,\beta_1) & \rho_2 Y_1^{1}(\alpha_2,\beta_2) & \cdots & \rho_N Y_1^{1} (\alpha_N,\beta_N)\\ 
     \rho_1 Y_1^{0}(\alpha_1,\beta_1) & \rho_2 Y_1^{0}(\alpha_2,\beta_2) & \cdots & \rho_N Y_1^{0} (\alpha_N,\beta_N)\\
     \rho_1 Y_1^{-1}(\alpha_1,\beta_1) & \rho_2 Y_1^{-1}(\alpha_2,\beta_2) & \cdots & \rho_N Y_1^{-1} (\alpha_N,\beta_N)\\
         \cdots & & &\\ 
         \rho_1^p Y_p^{p}(\alpha_1,\beta_1) & \rho_2^p Y_p^{p}(\alpha_2,\beta_2) & \cdots & \rho_N^p Y_p^{p} (\alpha_N,\beta_N)\\ 
         \vdots & & &\\ 
         \rho_1^p Y_p^{0}(\alpha_1,\beta_1) & \rho_2^p Y_p^{0}(\alpha_2,\beta_2) & \cdots & \rho_N^p Y_p^{0} (\alpha_N,\beta_N)\\
         \vdots & & &\\ 
     \rho_1^p Y_p^{-p}(\alpha_1,\beta_1) & \rho_2^p Y_p^{-p}(\alpha_2,\beta_2) & \cdots & \rho_N^p Y_p^{-p} (\alpha_N,\beta_N)\\
     \end{bmatrix},\label{eq:vapp}
 \end{gather} and 
 $\vec{q} = \begin{bmatrix} 
     q_1 & q_2 & \cdots & q_N
     \end{bmatrix}^T$.
 Then the potential $\{\phi_i\}_{i=1}^M$ at $M$ locations which are located outside the sphere of radius `$a$' due to $N$ particles (where $N\leq M$) located inside a sphere of radius `$a$' can be approximated as $\{\tilde{\phi}_i\}_{i=1}^M$ using a matrix $\tilde{A}$ with rank-$\tau$ through multipole expansions. 
 \begin{equation}
     \tilde{A} = UV,
 \end{equation}
 where $U = \begin{bmatrix} u_1 & u_2 & \cdots & u_M
 \end{bmatrix}^T\quad \text{and } U\in\mathbb{R}^{M\times(p+1)^2},$ $V\in\mathbb{R}^{(p+1)^2\times N}$ as in~\cref{eq:vapp}.
 The matrix $\tilde{A}$ can have a maximum rank of $\tau =\bkt{p+1}^2$ and for any $p$, we have $\abs{\phi_i - \bkt{\tilde{A}\vec{q}}_i} < \epsilon$, where $\epsilon > 0$ and $1\leq i \leq M$.
\end{corollary}
\begin{proof}[Proof of~\cref{tm:ranks3D}]
Consider the computational domain, a hypercube $\bm{B}$ in $3$D, which is hierarchically subdivided as in~\cref{fig:octTree} and represented using an octree, $\mathcal{T}^L$. Then each cube in a level has at the most $4$ types of interaction, viz., well-separated interaction, vertex-sharing interaction, edge-sharing interaction and face-sharing interaction. For a better interactive view of the figures used in the proof, please check \url{https://kandapva.github.io/hodlr3d/}.

\textbf{Case (i)}:(\underline{Boxes $\Bb_1$ and $\Bb_2$ are at least one box away}).
As shown in~\cref{fig:3d_far}, consider two boxes $\Bb_1 = [0,l]^3$ and $\Bb_2 = [2l,3l]\times[0,l]^2$, that is one box away. 
	\begin{figure}[H]
		\centering
		\tdplotsetmaincoords{40}{30}
		\begin{tikzpicture}[tdplot_main_coords]
		\cubex{0}{0}{0}{2}{line width=0.5mm,red}{$\Bb_1$};
		\cubex{4}{0}{0}{2}{line width=0.5mm,black}{$\Bb_2$};
		\draw[thick,<->,dashed, color=blue] (2,-0.2,0) -- (4,-0.2,0) node[midway, below]{$l$};
		\draw[thick,<->,dashed, color=blue] (2.2,0,0) -- (2.2,2,0) node[midway, below]{$l$};
		\draw[thick,<->,dashed, color=blue] (0,-0.2,0) -- (2,-0.2,0) node[midway, below]{$l$};
		\draw[thick,<->,dashed, color=blue] (4,-0.2,0) -- (6,-0.2,0) node[midway, below]{$l$};
		\end{tikzpicture}
		\caption{Far field boxes which in this case are "one box away"}
		\label{fig:3d_far}
	\end{figure}
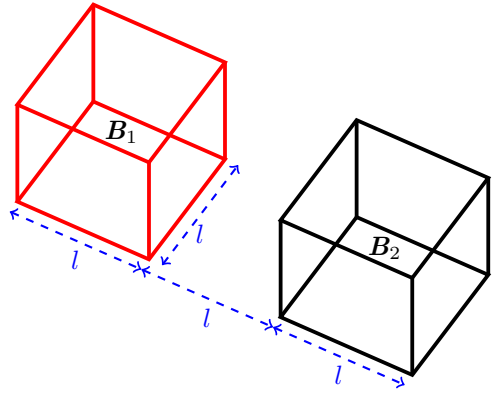
	Then we have a sphere $S:(r_s,\vec{c}_s)$ that circumscribes the box $\Bb_1$, where radius $r_s = \dfrac{l\sqrt{3}}{2}$  and center $\vec{c}_s = \bkt{\dfrac{l}{2},\dfrac{l}{2},\dfrac{l}{2}}$. The sphere $S$ encloses all the charges in the Box $\Bb_1$ such that $|\rho_i| < \dfrac{l\sqrt{3}}{2}$ and the minimum distance between the Box $B_2$ and the center of the sphere is $\dfrac{3l}{2}$. So in the Multipole Expansion lemma, we have $a  =  \dfrac{l\sqrt{3}}{2}$ and $c = \dfrac{\dfrac{3l}{2}}{\dfrac{l\sqrt{3}}{2}} = \sqrt{3}$. Then from~\cref{tm:prank}, we can approximate the potential using a rank-$\tau$ approximation such that,
 \begin{displaymath}
     \abs{\phi_i-\tilde{\phi}_i} \leq \dfrac{2Q}{3l} \bkt{\dfrac{1}{\sqrt3-1}} \bkt{\dfrac{1}{\sqrt{3}}}^{p+1},
 \end{displaymath}
and choosing $p+1 = \ceil{\dfrac{\log{\dfrac{2Q}{3l\epsilon (\sqrt3 - 1)}}}{\log{\sqrt3}}}$ ensures,
$\abs{\phi_i-\tilde{\phi}_i} \leq \epsilon.$ The rank $\tau$ is bounded above by $\ceil{\dfrac{\log{\dfrac{2Q}{3l\epsilon (\sqrt3 - 1)}}}{\log{\sqrt3}}}^2$.
So when the boxes are separated at least by $l$, we have a rank-$\tau$ matrix approximation $\tilde{A}$, where 
\begin{displaymath}
    \boxed{\tau \in  \mathcal{O} \bkt{\log^2\bkt{Q/l\epsilon}}}\quad \text{such that,}\quad \abs{\phi_i-\bkt{\tilde{A}q}_i}<\epsilon.
\end{displaymath}

\textbf{Case (ii)}:(\underline{Boxes $\Bb_1$ and $\Bb_2$ are Vertex sharing neighbors}).
Consider two boxes $\Bb_1 = [0,l]^3$ and $\Bb_2 = [l,2l]^3$ as shown in~\cref{fig:3d_ver}, that shares at most a vertex.
	\begin{figure}[H]
		\centering
		\subfloat[Vertex sharing boxes]{
			\label{fig:3d_ver}
			\tdplotsetmaincoords{40}{30}
			\begin{tikzpicture}[tdplot_main_coords]
				\cubex{0}{0}{0}{2}{line width=0.5mm,red}{$\Bb_1$};
				\cubex{2}{2}{2}{2}{line width=0.5mm,black}{$\Bb_2$};
				\draw[thick,<->,dashed, color=blue] (0,-0.2,0) -- (2,-0.2,0) node[midway, below]{$l$};
				\draw[thick,<->,dashed, color=blue] (2.2,0,0) -- (2.2,2,0) node[midway, below]{$l$};
			\end{tikzpicture}
		}\quad%
		\subfloat[Subdivision]{
			\label{fig:3d_ver_sub}
			\tdplotsetmaincoords{40}{30}
			\begin{tikzpicture}[tdplot_main_coords]
				\cubex{0}{0}{0}{2}{line width=0.5mm,red}{};
				\cubex{0}{0}{0}{1}{draw=none,line width=0.5mm,red}{$\Bb^{(1)}_1$};
				\cubex{1}{1}{1}{1}{line width=0.3mm,blue}{$\bar{\Bb}^{(1)}_1$};
				\cubex{2}{2}{2}{2}{line width=0.5mm,black}{$\Bb_2$};
				\draw[thick,<->,dashed, color=blue] (0,-0.2,0) -- (2,-0.2,0) node[midway, below]{$l$};
				\draw[thick,<->,dashed, color=blue] (2.2,0,0) -- (2.2,2,0) node[midway, below]{$l$};
			\end{tikzpicture}
		}
	\end{figure}
	To prove our theorem, we subdivide the cube $\Bb_1$ into 8 smaller cubes of length $l/2$. We form $\Bb^{(1)}_1$, which is the union of the 7 smaller cubes that do not share a boundary with $\Bb_2$. The remaining box, which shares a vertex with $\Bb_2$, is named as $\bar{\Bb}^{(1)}_1$.
    This is shown in~\cref{fig:3d_ver_sub}. For the charges in Box $\Bb^{(1)}_1$, we have sphere $S_1:(r_s^{(1)},\vec{c}_s^{(1)})$ that encloses all the charges which is centered at $\vec{c}_s^{(1)} = \bkt{\dfrac{l}{4},\dfrac{l}{4},\dfrac{l}{4}}$ and radius $r_s^{(1)} = \dfrac{l\sqrt{17}}{4}$. The minimum distance between the box $\Bb_2$ and the center of the sphere is $\dfrac{3l\sqrt{3}}{4}$ and hence in ~\cref{tm:multipole} we have $c = \sqrt{\dfrac{27}{19}}$. By using~\cref{tm:prank}, the potential, $\phi^{(1)}$ due to charges inside $\Bb^{(1)}_1$ can be approximated with $\tilde{\phi}^{(1)}$ through a rank-$\tau$ approximation. The corresponding matrix is $\tilde{A}^{(1)}$, and the potential $\tilde{\phi}^{(1)} = \tilde{A}^{(1)}\Vec{q}^{(1)}$ such that, 
 \begin{displaymath}
     \abs{\phi^{(1)}-\tilde{\phi}^{(1)}}\leq \dfrac{4Q}{l\bkt{\sqrt{27}-\sqrt{19}}}\bkt{\sqrt{\dfrac{19}{27}}}^{p+1}.
 \end{displaymath}
	If we choose $p+1 \in \ceil{\dfrac{\log\bkt{\dfrac{4Q}{l\epsilon_1\bkt{\sqrt{27}-\sqrt{19}}}}}{\log{\bkt{\sqrt{\dfrac{27}{19}}}}}}$ ensures $\abs{\phi^{(1)}_i-\bkt{\tilde{A}^{(1)}\vec{q}^{(1)}}_i} < \epsilon_1.$
  The rank $\tau$ of $\tilde{A}^{(1)}$ is bounded above by $\ceil{\dfrac{\log\bkt{\dfrac{4Q}{l\epsilon_1\bkt{\sqrt{27}-\sqrt{19}}}}}{\log{\bkt{\sqrt{\dfrac{27}{19}}}}}}^2$.
	We perform this hierarchical subdivision $\kappa$ times then, box $\Bb_1 = \bar{\Bb}^{(\kappa)}_1\dcup_{j=0}^{\kappa}\Bb^{(j)}_1$. Since we have assumed the $N$ charges be distributed uniformly in the box $\Bb_1$, the hierarchical subdivision with $\kappa>\log_8{\bkt{N}}$ results in the box $\bar{\Bb}^{(\kappa)}_1$ without any charges. Now for each box $\Bb^{(j)}_1|_{j=1}^{\kappa}$, we can circumscribe a sphere $S_j:(r_s^{(j)},\vec{c}_s^{(j)})$, with a radius of $r_s^{(j)}=\dfrac{l\sqrt{17}}{2^{j+1}}$. Then from~\cref{tm:multipole,tm:prank} the potential due to the charges inside the box $\Bb^{(j)}_1|_{j=1}^{\kappa}$ be $\phi^{(j)}_1|_{j=1}^{\kappa}$, for which we have approximate potential $\tilde{\phi}^{(j)}_1|_{j=1}^{\kappa}$ through a rank-$\tau$ matrix approximations. At each level of hierarchical subdivision we have $\tilde{A}^{(2)} \cdots \tilde{A}^{(\kappa)}$ with rank-$\tau$ approximation and at level $1\leq j\leq \kappa$ we have, $\abs{\phi^{(j)}-\tilde{\phi}^{(j)}}<\epsilon_j,$
 and $\epsilon_j = \dfrac{2^{j+1}Q}{l\bkt{\sqrt{27}-\sqrt{19}}}\bkt{\sqrt{\dfrac{19}{27}}}^{p_j+1}$,$\epsilon_j\quad \forall j\in\{1,\cdots,\kappa\},$ where $\epsilon_j = 2^{j-1}\epsilon_1$. $p_j+1$ is given by $\ceil{\dfrac{\log\bkt{\dfrac{2^{j+1}Q}{l\epsilon_1\bkt{\sqrt{27}-\sqrt{19}}}}}{\log{\bkt{\sqrt{\dfrac{27}{19}}}}}}$ and the rank $\tau_j$ at each level is bounded above by $\ceil{\dfrac{\log\bkt{\dfrac{2^{j+1}Q}{l\epsilon_1\bkt{\sqrt{27}-\sqrt{19}}}}}{\log{\bkt{\sqrt{\dfrac{27}{19}}}}}}^2$. Now for hierarchical subdivision of box $\kappa$ times , the potential $\phi = \dsum_{j=1}^{\kappa}\phi^{(j)}$ at $\Bb_2$ due to charges at $\Bb_1$ are approximated through rank-$\tau$ approximation $\dsum_{j=1}^{\kappa}\tilde{\phi}^{(j)}$ then,
	\begin{equation} \label{eq:ver_ineq}
		\begin{split}
		\abs{\phi_i - \tilde{\phi}_i} &= \abs{\dsum_{j=1}^{\kappa}\phi^{(j)}_i - \dsum_{j=1}^{\kappa}\tilde{\phi}^{(j)}_i} = \abs{\dsum_{j=1}^{\kappa}\bkt{\phi^{(j)}_i - \tilde{\phi}^{(j)}_i}}\\
		 & \leq \dsum_{j=1}^{\kappa}\abs{\phi^{(j)}_i - \tilde{\phi}^{(j)}_i} < \dsum_{j=1}^{\kappa}\epsilon_j
	\end{split}
	\end{equation}
	Using the relation between $\epsilon_j$, the sum $\dsum_{j=1}^{\kappa}\epsilon_j$ becomes $\epsilon_1\dsum_{j=1}^{\kappa}2^{j-1} = (2^{\kappa}-1)\epsilon_1$. We know that for uniform distribution of particles, the recursive subdivision of box $\Bb_1$ for $\kappa>\log_8\bkt{N}$ results in smaller boxes with no particles. So keeping $\kappa \approx \log_8\bkt{N}$ gets an error bound in $\epsilon$, i.e., we choose $\epsilon_1 = \dfrac{\epsilon}{\sqrt[3]{N}}$ to get, $\abs{\phi_i - \tilde{\phi}_i} < \epsilon$. 
	The rank-$\tau$ matrix approximation $\tilde{A}^{(j)}$, at a level $j$ is bounded above by
 \begin{displaymath}
     \ceil{\log^2 \bkt{\dfrac{2Q\sqrt[3]{N}}{l\epsilon}}}.
 \end{displaymath}
The total rank $\tau$ of the matrix involved in the approximation is bounded above by
\begin{displaymath}
    \kappa \ceil{\log^2 \bkt{\dfrac{2Q\sqrt[3]{N}}{l\epsilon}}}
\end{displaymath}
 and by setting the number of levels, $\kappa \approx \log_8\bkt{N}$, the rank of the matrix due to the interaction between two vertex sharing boxes  is bounded above by  
 \begin{displaymath}
     \boxed{\tau \in \mathcal{O}\bkt{ \log_8\bkt{N}\log^2 \bkt{\dfrac{2Q\sqrt[3]{N}}{l\epsilon}}},}\quad \text{where,} \quad \abs{\phi_i - \bkt{\tilde{A}q}_i} < \epsilon.
 \end{displaymath}

\textbf{Case (iii)}:(\underline{Boxes $\Bb_1$ and $\Bb_2$ are Edge sharing neighbors}).
Consider two boxes $\Bb_1 = [0,l]^3$ and $\Bb_2 = [l,2l]^2\times[0,l]$, that shares at most an edge as shown in~\cref{fig:3d_edge}.
	\begin{figure}[H]
		\centering
		\subfloat[Edge sharing boxes]{
			\label{fig:3d_edge}
			\tdplotsetmaincoords{40}{30}
			\begin{tikzpicture}[tdplot_main_coords]
				\cubex{0}{0}{0}{2}{line width=0.5mm,red}{$\Bb_1$};
				\cubex{2}{2}{0}{2}{line width=0.5mm,black}{$\Bb_2$};
				\draw[thick,<->,dashed, color=blue] (0,-0.2,0) -- (2,-0.2,0) node[midway, below]{$l$};
				\draw[thick,<->,dashed, color=blue] (2,1.8,0) -- (4,1.8,0) node[midway, below]{$l$};
			\end{tikzpicture}
		}\quad%
		\subfloat[Subdivision]{
			\label{fig:3d_edge_sub}
			\tdplotsetmaincoords{40}{30}
			\begin{tikzpicture}[tdplot_main_coords]
				\cubex{0}{0}{0}{2}{line width=0.5mm,red}{$\Bb_1^{(1)}$};
				\cubex{2}{2}{0}{2}{line width=0.5mm,black}{$\Bb_2$};
				\cubex{0}{0}{0}{1}{line width=0.3mm,blue}{};
				\cubex{0}{1}{0}{1}{line width=0.3mm,blue}{};
				\cubex{0}{0}{1}{1}{line width=0.3mm,blue}{};
				\cubex{1}{0}{0}{1}{line width=0.3mm,blue}{};
				\cubex{1}{0}{1}{1}{line width=0.3mm,blue}{};
				\cubex{0}{1}{1}{1}{line width=0.3mm,blue}{};
				\draw[thick,<->,dashed, color=blue] (0,-0.2,0) -- (2,-0.2,0) node[midway, below]{$l$};
				\draw[thick,<->,dashed, color=blue] (2,1.8,0) -- (4,1.8,0) node[midway, below]{$l$};
			\end{tikzpicture}
		}
	\end{figure}
	In this case, we will subdivide the box $\Bb_1$ into $8$ smaller boxes of side length $l/2$. We group the boxes that do not share a boundary with the box $\Bb_2$ as $\Bb_1^{(1)}=\bigcup_{j=0}^{6}\Bb_1^{(1,j)}$, whereas the remaining two boxes are $\bar{\Bb}_1^{(1,1)}$, $\bar{\Bb}_1^{(1,2)}$ (i.e., $\Bb_1-\Bb^{(1)}_1 = \dbcup_{i=1}^{2} \bar{\Bb}_1^{(1,i)}$).  So, for the box $\Bb^{(1)}$, that donot share a boundary with $\Bb_2$, we can circumscribe a sphere $S:(r_s,\vec{c}_s)$ with center $\vec{c}_s = \bkt{\dfrac{l}{2},\dfrac{l}{4},\dfrac{l}{4}}$ and radius $ r_s = l\sqrt{\dfrac{7}{8}}$ that encloses all the charges in the box $\Bb_1^{(1,j)}$. For the Multipole expansions lemma we have $c = \dfrac{3}{\sqrt{7}}$, then the potential $\phi^{(1)}$ at box $\Bb_2$ due to the charges at box $\Bb^{(1)}$ be given as $\phi^{(1)}$. Then from~\cref{tm:prank}, we have potential $\tilde{\phi}^{(1)}$ through a rank-$\tau$ approximation $A^{(1)}$ such that,
 \begin{displaymath}
     \abs{\phi_i^{(1)}-\tilde{\phi}_i^{(1)}}\leq \dfrac{2Q\sqrt{2}}{l\bkt{3-\sqrt{7}}}\bkt{\dfrac{\sqrt{7}}{3}}^{p+1}.
 \end{displaymath}
	Choosing $p+1 = \ceil{\dfrac{\log\bkt{\dfrac{2\sqrt{2}Q}{l\bkt{3-\sqrt{7}}\epsilon_1}}}{\log\bkt{3/\sqrt{7}}}}$, ensures 
	$\abs{\phi_i - \bkt{\tilde{A}^{(1)}q^{(1,j)}}_i} \leq \epsilon_1.$
 And, the rank of the matrix $\tilde{A}^{(1)}$ is bounded above by 
 \begin{displaymath}
\tau^{(1)}=\ceil{\dfrac{\log\bkt{\dfrac{2\sqrt{2}Q}{l\bkt{3-\sqrt{7}}\epsilon_1}}}{\log\bkt{3/\sqrt{7}}}}^2.
 \end{displaymath}	
	Again, each of the remaining boxes $\bar{\Bb}^{(1,1)}$, $\bar{\Bb}^{(1,2)}$ that share the edge with $\Bb_2$ are again subdivided into $8$ smaller boxes of side length $l/4$. We form the union of boxes that do not share a boundary with $\Bb_2$ as $\Bb^{(2,1)}$ and $\Bb^{(2,2)}$. The boxes $\bar{\Bb}^{(2,1)}$ and $\bar{\Bb}^{(2,2)}$ that share a boundary with $\Bb_2$. So, continuing this hierarchical subdivision, for a  level $k$, we will have boxes $\Bb^{(k,1)},\cdots,\Bb^{(k,2^{k-1})}$ that do not share a boundary and $\bar{\Bb}^{(k,1)},\cdots,\bar{\Bb}^{(k,2^{k-1})}$ that share a boundary. Hence for a box at level $k$ that does not share a boundary with $\Bb_2$, we can circumscribe a sphere such that $c = \dfrac{3}{\sqrt{7}}$. The potential $\phi^{(k,j)}$ at box $\Bb_2$ due to the charges at box $\Bb^{(k,j)}$ be given as $\phi^{(k,j)}$. Then from~\cref{tm:prank}, we have potential $\tilde{\phi}^{(k,j)}$ through a rank-$\tau$ approximation $\tilde{A}^{(k,j)}$ such that,
 \begin{displaymath}
     \abs{\phi_i^{(k,j)}-\tilde{\phi}_i^{(k,j)}}\leq \dfrac{2^jQ\sqrt{2}}{l\bkt{3-\sqrt{7}}}\bkt{\dfrac{\sqrt{7}}{3}}^{p+1}.
 \end{displaymath}
If we choose $p+1 = \ceil{\dfrac{\log\bkt{\dfrac{2^j\sqrt{2}Q}{l\bkt{3-\sqrt{7}}\epsilon_k}}}{\log\bkt{3/\sqrt{7}}}}$ ensures $\abs{\phi^{(k,j)}_i - \bkt{\tilde{A}^{(k,j)}q^{(k,j)}}_i} \leq \epsilon_k$.
 The rank of the matrix $\tilde{A}^{(k,j)}$ is bounded above by 
$\tau^{(k,j)}=\ceil{\dfrac{\log\bkt{\dfrac{2^j\sqrt{2}Q}{l\bkt{3-\sqrt{7}}\epsilon_k}}}{\log\bkt{3/\sqrt{7}}}}^2.$
 Using the approximate matrices at level $k$, we can form a matrix $\tilde{A}^{(k)}$ through the matrices $\tilde{A}^{(k,1)},\cdots,\tilde{A}^{(k,2^{k-1})}$ whose rank is bounded above by
 \begin{displaymath}
    \tau^{(k)}=2^{k-1}\ceil{\dfrac{\log\bkt{\dfrac{2^j\sqrt{2}Q}{l\bkt{3-\sqrt{7}}\epsilon_k}}}{\log\bkt{3/\sqrt{7}}}}^2, 
 \end{displaymath}
 such that $\abs{\phi^{(k)}_i - \bkt{\tilde{A}^{(k)}q^{(k)}}_i} < 2^{k-1}\epsilon_k.$ The hierarchical subdivision of box $\Bb_1$ for $\kappa$ times results in series of approximants $\tilde{A}^{(1)},\cdots,\tilde{A}^{(\kappa)}$ with which we can construct $\tilde{A}$ such that 
\begin{equation} \label{eq:edge_ineq}
    \begin{split}
    \abs{\phi_i - \bkt{\tilde{A}q}_i} & < \dsum_{k=1}^{\kappa}\abs{\bkt{\bkt{A_jq}_i - \bkt{\tilde{A}_jq}_i}}\\
     & < \dsum_{k=1}^{\kappa}2^{k-1}\epsilon_k.
\end{split}
\end{equation}
As per our previous construct, $\kappa\approx\log_8\bkt{N}$, then using the relation between $\epsilon_j$ and $\epsilon_1$, the sum $\dsum_{k=1}^{\kappa}\epsilon_j$ becomes $\epsilon_1\dsum_{k=1}^{\kappa}2^{k-1}2^k = \dfrac{2}{3}(4^{\kappa}-1)\epsilon_1 \approx \sqrt[3]{N^2}\epsilon_1$ and to get error bound in $\epsilon$, we choose $\epsilon_1 = \dfrac{\epsilon}{\sqrt[3]{N^2}}$ then,
\begin{displaymath}
    \abs{\phi_i - \bkt{\tilde{A}q}_i} < \epsilon.
\end{displaymath}
 The rank of $\tilde{A}$ can be at most $\dsum_{j=1}^{\kappa}2^j\ceil{\log^2\bkt{\dfrac{2QN}{l\epsilon}}}$.
Hence with $\kappa\approx\log_8\bkt{N}$, the rank of $\tilde{A}$ is given by 
\begin{displaymath}
    \boxed{\tau \in \mathcal{O}\bkt{N^{1/3}\log^2\bkt{\dfrac{2Q\sqrt[3]{N^2}}{l\epsilon}}},}\quad  \text{with} \quad \abs{\phi_i - \bkt{\tilde{A}q}_i} < \epsilon.
\end{displaymath}

\textbf{Case (iv)}:(\underline{Boxes $\Bb_1$ and $\Bb_2$ are Face sharing neighbors}).
 Consider two boxes sharing a face i.e., $\Bb_1= [0,l]^3$ and $\Bb_2 = [l,2l]\times[0,l]^2$, as shown in~\cref{fig:3d_face}.
	\begin{figure}[!htbp]
		\centering
		\subfloat[Face sharing boxes]{
			\label{fig:3d_face}
			\tdplotsetmaincoords{40}{30}
			\begin{tikzpicture}[tdplot_main_coords]
				\cubex{0}{0}{0}{2}{line width=0.5mm,red}{$\Bb_1$};
				\cubex{2}{0}{0}{2}{line width=0.5mm,black}{$\Bb_2$};
				\draw[thick,<->,dashed, color=blue] (0,-0.2,0) -- (2,-0.2,0) node[midway, below]{$l$};
				\draw[thick,<->,dashed, color=blue] (2,-0.2,0) -- (4,-0.2,0) node[midway, below]{$l$};
			\end{tikzpicture}
		}\quad%
		\subfloat[Subdivision]{
			\label{fig:3d_face_sub}
			\tdplotsetmaincoords{40}{30}
			\begin{tikzpicture}[tdplot_main_coords]
				\cubex{0}{0}{0}{2}{line width=0.1mm,red}{};
				\cubex{0}{0}{0}{1}{line width=0.3mm,blue}{$\Bb^{(1,1)}_1$};
				\cubex{0}{1}{0}{1}{line width=0.3mm,blue}{$\Bb^{(1,2)}_1$};
				\cubex{0}{0}{1}{1}{line width=0.3mm,blue}{$\Bb^{(1,3)}_1$};
				\cubex{0}{1}{1}{1}{line width=0.3mm,blue}{$\Bb^{(1,4)}_1$};
				\cubex{2}{0}{0}{2}{line width=0.5mm,black}{$\Bb_2$};
				\draw[thick,<->,dashed, color=blue] (0,-0.2,0) -- (2,-0.2,0) node[midway, below]{$l$};
				\draw[thick,<->,dashed, color=blue] (2,-0.2,0) -- (4,-0.2,0) node[midway, below]{$l$};
			\end{tikzpicture}
		}
	\end{figure}
	The proof for this case starts by subdividing the box $\Bb_1$ into eight smaller boxes, each with side length $\frac{l}{2}$. The four boxes that do not share a boundary with the box $\Bb_2$ are named $\Bb_1^{(1,1)},\Bb_1^{(1,2)},\Bb_1^{(1,3)}$ and $\Bb_1^{(1,4)}$. The boxes that share a face with box $\Bb_2$ are given as  $\bar{\Bb}_1^{(1)}=\bar{\Bb}_1^{(1,1)}\cup\bar{\Bb}_1^{(1,2)}\cup\bar{\Bb}_1^{(1,3)}\cup\bar{\Bb}_1^{(1,4)}$ It can be shown as in previous cases that for each box $\Bb_1^{(1,1)},\Bb_1^{(1,2)},\Bb_1^{(1,3)}$ and $\Bb_1^{(1,4)}$ we can circumscribe a sphere of radius $\dfrac{l\sqrt{3}}{4}$. Now the boxes $B_{1}^{(1,i)}$ are separated from the source points by a distance of $\dfrac{3l}{4}$. In the multipole expansions lemma for the boxes $\Bb_1^{(1,j)}$ ($1\leq j \leq 4$),  we have $c = \sqrt{3}$. The potential due to charges inside $\Bb_1^{(1,j)}$ be $\phi^{(1,j)}$ and by using~\cref{tm:multipole,tm:prank}, we have potential $\tilde{\phi}^{(1,j)}$ through a rank-$\tau$ approximation , $\tilde{A}^{(1,j)}_1$, such that,
 \begin{displaymath}
     \abs{\phi^{(1,j)}_i-\tilde{\phi}^{(1,j)}_i}\leq \dfrac{4Q}{3l\bkt{\sqrt{3}-1}}\bkt{\sqrt{\dfrac{1}{3}}}^{p+1}.
 \end{displaymath}
	If we choose $p+1 \in \ceil{\dfrac{\log\bkt{\dfrac{4Q}{3l\epsilon_1\bkt{\sqrt{3}-1}}}}{\log{\bkt{\sqrt{3}}}}}$ ensures 
		$\abs{\phi^{(1,j)}_i-\bkt{\tilde{A}_1^{(1,j)}q}_i} < \epsilon_1.$
    The matrix $\tilde{A}^{(1,j)}_1$ has a rank bounded above by  $\ceil{\dfrac{\log\bkt{\dfrac{4Q}{3l\epsilon_1\bkt{\sqrt{3}-1}}}}{\log{\bkt{\sqrt{3}}}}}^2$. The potential at box $\Bb_2$ due to charges in $\Bb_1^{(1)}$ is $\phi^{(1)} = \dsum_{j=0}^{4}\phi^{(1,j)}$ which is approximated using rank-$\tau$ approximation through $\tilde{\phi}^{(1)}= \dsum_{j=0}^{4}\tilde{\phi}^{(1,j)}$ then, 
    $\abs{\phi^{(1)}_1 - \tilde{\phi}^{(1)}}\leq\dsum_{j=0}^{4}\abs{\phi^{(1,j)}_1 - \tilde{\phi}^{(1,j)}} < 4\epsilon_1.$ We construct the matrix $\tilde{A}^{(1)}$ using the rank-$\tau$ approximations $\tilde{A}^{(1,1)},\cdots,\tilde{A}^{(1,4)}$ corresponding to the boxes $\Bb_1^{(1,1)},\cdots,\Bb_1^{(1,4)}$ whose rank is bounded above by $\tau_1 = 4\ceil{\dfrac{\log\bkt{\dfrac{4Q}{3l\epsilon_1\bkt{\sqrt{3}-1}}}}{\log{\bkt{\sqrt{3}}}}}^2$ such that,
    $\abs{\phi^{(1)}_i-\bkt{\tilde{A}^{(1)}q}_i} < 4\epsilon_1.$
    The box $\bar{\Bb}^{(1)} = \Bb^{(2)}\cup\bar{\Bb}^{(2)}$, where $\Bb^{(2)}= \bigcup_{j=1}^{4^2}\Bb^{(2,j)}$ boxes do not share a boundary and $\bar{\Bb}^{(2)}= \bigcup_{j=1}^{4^2}\bar{\Bb}^{(2,j)}$ shares face with box $\Bb_2$. The boxes $\Bb^{(2,j)}$ and $\bar{\Bb}^{(2,j)}$ have side length $\dfrac{l}{4}$. Similarly, for a level $k$, we have  $\bar{\Bb}^{(k-1)} = \Bb^{(k)}\cup\bar{\Bb}^{(k)}$ and $\Bb^{(k)}= \bigcup_{j=1}^{4^k}\Bb^{(2,j)}$ do not share a boundary whereas $\bar{\Bb}^{(k)}= \bigcup_{j=1}^{4^k}\bar{\Bb}^{(2,j)}$ shares boundary with box $\Bb_2$. Applying~\cref{tm:multipole,tm:prank}, we have potential at $\Bb_2$ due to charges inside $\Bb^{(k)}$ as $\phi^{(k)} = \dsum_{j=0}^{4^k}\phi^{(k,j)}$. The potential $\phi^{(k)}$ is approximated using a rank-$\tau$ approximation (i.e., $\tilde{A}^{(k,j)}$) by $\tilde{\phi}^{(k)} = \dsum_{j=0}^{4^k}\tilde{\phi}^{(k,j)}$ such that,
    \begin{displaymath}
        \abs{\phi^{(k,j)} - \tilde{\phi}^{(k,j)}} \leq \dfrac{2.2^jQ}{3l\bkt{\sqrt{3}-1}}\bkt{\dfrac{1}{\sqrt{3}}}^{p+1}.
    \end{displaymath}
If we choose $p+1 = \ceil{\dfrac{\log\bkt{\dfrac{2.2^jQ}{3l\epsilon_j\bkt{\sqrt{3}-1}}}}{\log{\bkt{\sqrt{3}}}}}$ makes $\abs{\phi^{(k,j)}-\bkt{\tilde{A}^{(k,j)}q}_i} < \epsilon_j.$
    Hence at each level, we have, 
    \begin{equation}
        \begin{split}
		\abs{\phi_i^{(k)} - \tilde{\phi}_i^{(k)}} &\leq \dsum_{j=1}^{4^k}\abs{\phi^{(k,j)}_i - \tilde{\phi}^{(k,j)}_i}\\
		 & < 4^k\epsilon_k.
	\end{split}
    \end{equation}
    At level $k$, we can construct a matrix $\tilde{A}^{(k)}$ using the rank-$\tau$ approximations using $\tilde{A}^{(k,1)},\cdots,\tilde{A}^{(k,4^k)}$ whose rank is bounded above by 
    \begin{displaymath}
        \tau^{(k)} = 4^k\ceil{\dfrac{\log\bkt{\dfrac{2.2^kQ}{3l\epsilon_j\bkt{\sqrt{3}-1}}}}{\log{\bkt{\sqrt{3}}}}}^2,
    \end{displaymath}
such that $\abs{\phi^{(k)}-\bkt{\tilde{A}^{(k)}q^{(k)}}_i} < \epsilon_j.$ Hierarchically subdividing the box $\kappa=\log_8\bkt{N}$ times, results in $\bar{\Bb}^{(\kappa)}$ without any charges. The potential $\phi=\dsum_{k=1}^{\kappa}\phi^{(k)}$ is approximated through rank-$\tau$ approximations $\tilde{\phi}=\dsum_{k=1}^{\kappa}\tilde{\phi}^{(k)}$ with $\epsilon_1 = \epsilon_k/2^k$ then,
    \begin{equation}
        \begin{split}
		\abs{\phi_i - \tilde{\phi}_i} &\leq \dsum_{k=1}^{\kappa}\abs{\phi^{(k)}_i - \tilde{\phi}^{(k)}_i}\\
		 & < \dfrac{4}{7}(8^{\kappa}-1)\epsilon_1 < 8^{\kappa}\epsilon_1
	\end{split}
    \end{equation}
    Thus to get error bound in $\epsilon$, we choose $\epsilon_1 = \dfrac{\epsilon}{N}$, with which we can construct a matrix $\tilde{A}$ through the approximants $\tilde{A}^{(1)},\cdots,\tilde{A}^{(\kappa)}$ whose rank is bounded above by
    \begin{displaymath}
        \tau = N^{2/3}\ceil{\dfrac{\log\bkt{\dfrac{2QN}{3l\epsilon\bkt{\sqrt{3}-1}}}}{\log{\bkt{\sqrt{3}}}}}^2,
    \end{displaymath}
then, $\abs{\phi_i - \bkt{\tilde{A}q}_i} < \epsilon$. Hence with $\kappa=\log_8\bkt{N}$, the rank of $\tilde{A}$ is given by
\begin{displaymath}
    \boxed{\tau \in \mathcal{O}\bkt{N^{2/3}\log^2\bkt{\dfrac{2QN}{l\epsilon}}},}\quad\text{with,}\quad \abs{\phi_i - \bkt{\tilde{A}q}_i} < \epsilon.
\end{displaymath}
\end{proof}

\begin{figure}[!htbp]
      \begin{center} 
    \includegraphics[width=0.6\linewidth]{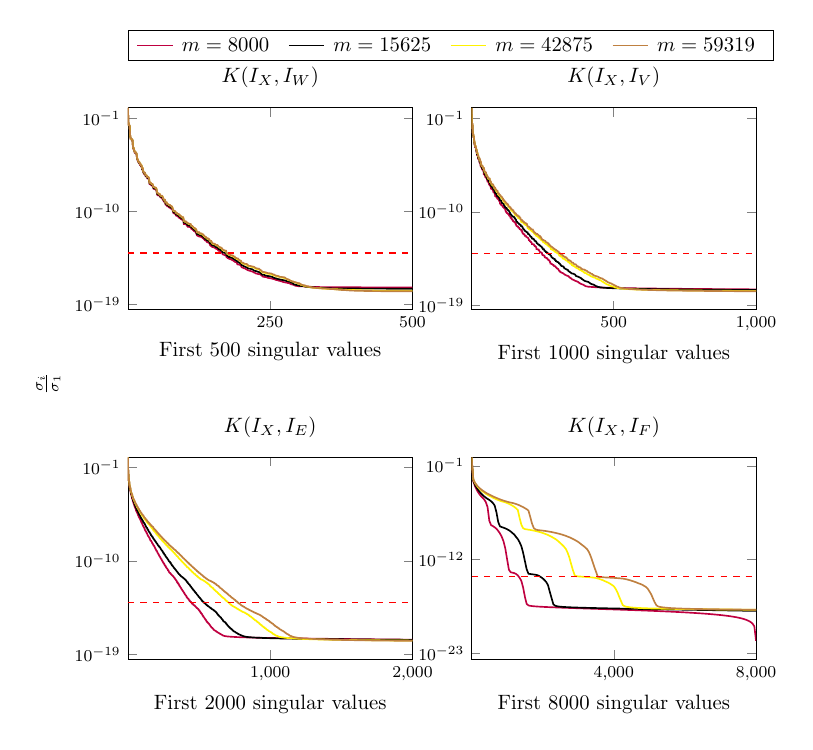}
        \end{center}
        \caption{Plot of singular values $\sigma_{i}$ normalized with the first singular value versus index $i$, with the matrix entries defined as in Equation~\eqref{eq:Laplace3D}.
        } \label{singularValues}
\end{figure}

We illustrate the decay of singular values of the interaction matrix for all the four cases of interactions (discussed in~\cref{tm:ranks3D}) for the Laplacian kernel in $3$D in Figure~\ref{singularValues}. It is to be observed that the decay of the singular values is fastest for the well-separated interaction $K(I_{X}, I_{W})$, slowest for the face-sharing interaction $K(I_{X}, I_{F})$ and that of the vertex sharing interaction $K(I_{X}, I_{V})$ is almost as fast as that of the well-separated interaction. Although we assumed a uniform distribution of particles in the domain,~\Cref{tm:ranks3D} also holds for a quasi-uniform distribution of particles.\footnote{Consider a hypercube $\Bb\subset \mathbb{R}^3$ contains $N$ particles. The particles inside $\Bb$ are said to be quasi-uniform distributed if exactly one particle is located inside each smallest hypercube resulting from the hierarchical subdivision of the hypercube $\Bb$ using an $\log_8\bkt{N}$ level octree.} 
 
\section{HODLR3D}\label{sec:HODLR3D}
\begin{figure}[!htbp]
    \centering
    \subfloat[Level 1]{
        \includegraphics[scale=0.3]{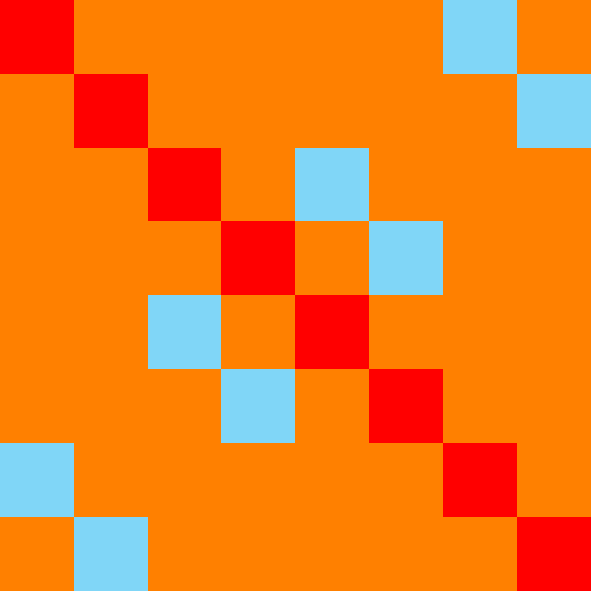}
        }\quad%
    \subfloat[Level 2]{
        \includegraphics[scale=0.3]{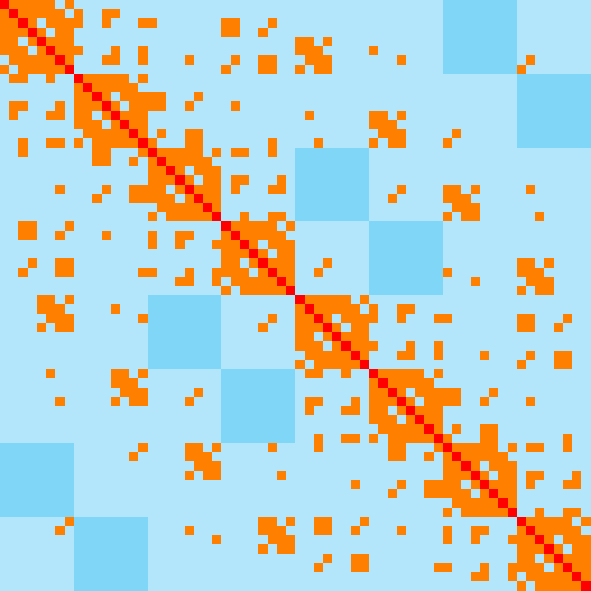}
        }\quad%
    \subfloat{
        \begin{tikzpicture}
            [
            box/.style={rectangle,draw=black, minimum size=0.25cm},scale=0.2
            ]
            \node[box,fill=red,,font=\tiny,label=right:Full-rank matrix (self interaction),  anchor=west] at (0,6){};
            \node[box,fill=orange,font=\tiny,label=right:Full-rank matrix (edge \& face sharing interaction),  anchor=west] at (0,4){};
            \node[box,fill=cyan!50,,font=\tiny,label=right:Low-rank matrix (level $1$),  anchor=west] at (0,2){};
            \node[box,fill=cyan!30,font=\tiny,label=right:Low-rank matrix (level $2$),  anchor=west] at (0,0){};
        \end{tikzpicture}
    }
    \caption{HODLR3D matrix low-rank structure at levels 1 and 2}
    \label{fig:HODLR3D_matrix}
\end{figure}

The objective of HODLR3D is to develop an almost linear complexity algorithm for Hierarchical matrices by leveraging the fact that the ranks of vertex-sharing interactions and well-separated interactions grow slowly with $N$. We have seen in the previous section  that the ranks of edge and face-sharing interactions grow as $\mathcal{O}(\sqrt[3]{N}\log^{2}(N))$ and $\mathcal{O}(\sqrt[3]{N^{2}}\log^{2}(N))$ respectively. And the rank of vertex sharing interaction grows as $\mathcal{O}(\log^{3}{N})$. Further, the rank of well-separated interaction does not grow with $N$. A similar rank growth for the interaction matrices arising out of a wide range of kernels in $3$D has been proved in~\cite{khan2022numerical}. Therefore, by choosing to compress only the well-separated and vertex-sharing interactions, we construct our HODLR3D, a Hierarchical matrix representation for kernel matrices from $3$D that yields almost linear complexity matrix algorithms. Before we describe the HODLR3D matrix-vector product algorithm, we define some notations in~\cref{table:Notations} that will be used in the rest of the section. 
\begin{definition}{Admissibility condition for HODLR3D:}\label{def:admh3d}
    Two clusters $\Cl$ and $\mathcal{C}^{(l)}_{j}$ where $i\neq j$ are admissible clusters, iff either they do not share a boundary, or they share a boundary which can be at the most a vertex.
    \end{definition}
In~\cref{fig:HODLR3D_matrix}, we illustrate the low-rank structure of the HODLR3D matrix at levels 1 and 2. 

\begin{table}[!htbp]
    \centering
    \caption{List of notations followed in the rest of the section}
    \label{table:Notations}
    \begin{tabularx}{\textwidth}{|l|X|}
        \hline
        $\Cl$ & \textit{cluster} of particles lying in cube $i$ at level $l$ of the octree \\ \hline
        $\mathcal{V}_{i}^{\ell}$ & $\{\mathcal{C}_{j}^{(l)}: \text{cubes $i$ and $j$ share only a vertex}\}$ \\ \hline
        $\mathcal{E}_{i}^{\ell}$ & $\{\mathcal{C}_{j}^{(l)}: \text{cubes $i$ and $j$ share only an edge}\}$ \\ \hline
        $\mathcal{F}_{i}^{\ell}$ & $\{\mathcal{C}_{j}^{(l)}: \text{cubes $i$ and $j$ share a face}\}$ \\ \hline
        $child(\Cl)$ & $\{\mathcal{C}_{j}^{(l+1)}: \text{cube $j$ is a child of cube $i$}\}$ \\ \hline
        $parent(\Cl)$ & $\mathcal{C}_{j}^{(l-1)}$, where cube $j$ is the parent of cube $i$ \\ \hline
        $siblings(\Cl)$ & $child(parent(\Cl))$\textbackslash $\Cl$ \\ \hline
        $\hat{C}_{\Cl}$ & Clan set of cluster $\Cl$ that is defined as\\
        &{\footnotesize $ \hat{C}_{\Cl} = \{siblings(\Cl)\bigcup \{ child(D):D\in\{\mathcal{E}_{parent(\Cl)}\bigcup \mathcal{F}_{parent(\Cl)}\}\}\}$}\\ \hline
        $\mathcal{I}_\Cl$ & Interaction list of cluster $\Cl$ that is defined as\\
        & {\footnotesize$\mathcal{I}_\Cl = \hat{C}_\Cl\cap \bkt{\mathcal{V}_\Cl  \cup \mathcal{W}_\Cl}$} \\ 
        \hline
    \end{tabularx}
 \end{table}
\Cref{tab:dense3d} provides us with the total number of dense matrix blocks for different hierarchical matrices based on an $L$-level octree. The total number of dense matrix blocks in $\mathcal{H}$-matrix with $\eta=\sqrt3$ is higher than HODLR3D. Similarly,~\Cref{tab:lrfar,tab:lrvert,tab:lr3def} shows that the total number of matrix blocks approximated as low-rank matrices is higher for $\mathcal{H}$-matrix with $\eta=\sqrt3$ compared to HODLR3D. This means that using HODLR3D, we can construct almost linear scaling algorithms by processing fewer individual matrix blocks, which is an appealing feature for extending the HODLR3D-based algorithms to parallel machines.

\begin{table*}[!htbp]
\caption{Total number of dense full-rank matrix blocks for different hierarchical matrices in $3$D based on an octree $\mathcal{T}^L$}
\label{tab:dense3d}
    \begin{tabular}{@{}cc@{}}
    \toprule
    Hierarchical matrix                       & Total number of dense blocks \\ \midrule
    HODLR (edge, vertex, face)                & $8^L$                         \\
    HODLR (edge, vertex)                      & $9 \cdot 8^{L} - 2\cdot 4^{L+1}$                             \\
    HOLDR3D                                   & $\frac{1}{3}\bkt{67 \cdot 8^{L} - 120 \cdot 4^{L} + 56 \cdot 2^{L}}$                           \\
    $\mathcal{H}$-matrix with $\eta=\sqrt{3}$ & $ \frac{1}{7}\bkt{223 \cdot 8^{L} - 126 \cdot 4^{L+1} + 49 \cdot 2^{L+3} - 104}$                            \\ \bottomrule
    \end{tabular}
\end{table*}

\begin{table*}[!htbp]
    \caption{Total number of low-rank matrix blocks due to well-separated interactions for different hierarchical matrices in $3$D based on an octree $\mathcal{T}^L$} 
    \label{tab:lrfar}
    \begin{tabular}{@{}cc@{}}
    \toprule
    Hierarchical matrix                       & Total number of low-rank blocks \\ \midrule
    HODLR (edge, vertex, face) &
    -- \\
    HODLR (edge, vertex) & $\frac{1}{21}\bkt{18\cdot 8^{L+2} - 42\cdot 4^{L+3} + 1536}$
    \\
    HOLDR3D & $\frac{1}{21}\bkt{444\cdot 8^{L+1} - 63\cdot 4^{L+4} + 735\cdot 2^{L+5} - 10944}$
    \\
    $\mathcal{H}$-matrix with $\eta=\sqrt{3}$ & $\frac{1}{7}\bkt{223\cdot 8^{L+1} - 630\cdot 4^{L+2}x + 1519\cdot 2^{L+4} - 6552L - 16008}$
    \\ \bottomrule
    \end{tabular}
\end{table*}

\begin{table*}[!htbp]
    \caption{Total number of low-rank matrix blocks due to vertex-sharing interactions for different hierarchical matrices in $3$D based on an octree $\mathcal{T}^L$} 
    \label{tab:lrvert}
    \begin{tabular}{@{}cc@{}}
    \toprule
    Hierarchical matrix                       & Total number of low-rank blocks \\ \midrule
    HODLR (edge, vertex, face)                & $\frac{8}{7}\bkt{8^{L}-1}$                         \\
    HODLR (edge, vertex)                      & $\frac{1}{21}\bkt{15\cdot 8^{L+1} - 14\cdot 4^{L+2} + 104}$           \\
    HOLDR3D   & $\frac{1}{21}\bkt{25\cdot 8^{L+1} - 42\cdot 4^{L+2} + 49\cdot 2^{L+4} - 312}$
    \\
    $\mathcal{H}$-matrix with $\eta=\sqrt{3}$ & --   \\ \bottomrule
    \end{tabular}
\end{table*}
             
\begin{table*}[!htbp]
    \caption{Total number of low-rank blocks for different hierarchical matrices in $3$D}
    \label{tab:lr3def}
    \begin{tabular}{@{}ccc@{}}
    \toprule
    \multirow{2}{*}{Hierarchical matrix} &
        \multicolumn{2}{c}{Total number of low-rank blocks due to} \\ \cmidrule(l){2-3} 
        &
        edge-sharing interactions &
        face-sharing interactions \\ \midrule
    HODLR (edge, vertex, face) 
        & $\frac{16}{7}\bkt{8^{L}-1}$
        & $\frac{32}{7}\bkt{8^{L}-1}$\\
    HODLR (edge, vertex)
    & $\frac{1}{21}\bkt{30\cdot 8^{L+1} - 7\cdot 4^{L+3} + 208}$
    & -- \\
    HOLDR3D 
    & -- & --\\
    $\mathcal{H}$-matrix with $\eta=\sqrt{3}$  & -- & -- \\ \bottomrule
    \end{tabular}
\end{table*}

\subsection{HODLR3D algorithm for Matrix-Vector product}\label{ssec:HODLR3DRepres}
The scale of problems in 3D is remarkably higher than that observed in 1D or 2D. For instance, if we consider $1000$ particles in each dimension, then the number of points in 3D with a tensor product grid is $10^{9}$.
Hence, storing the low-rank factors of all the low-rank sub-blocks of such large systems may not be feasible, as it will exhaust the RAM eventually for large system sizes. With limited RAM available, one must avoid explicitly storing the low-rank factors. Due to this, the low-rank compression technique that we describe in this article differs from that of the HODLR2D algorithm described in article~\cite{H2DKandappan}.

The first task is to perform the initialization, as described in Algorithm~\ref{alg:initialization}. 
We use ACA~\cite{bebendorf2003adaptive,zhao2005adaptive,tyrtyshnikov2000incomplete,bebendorf2000approximation} to compress the low-rank sub-blocks. 
For each pair of clusters $\Cl$ and $\mathcal{C}_{j}^{(l)}\in \mathcal{I}_\Cl$ at all levels of the $\mathcal{T}^{L}$ tree, with index sets $X$ and $Y$ respectively, the interaction matrix $K(X,Y)$ is compressed using ACA with tolerance $\epsilon$ as follows.
\begin{equation}\label{eq:prank}
    K(X,Y) \approx K(X,\tau_{XY}) K(\sigma_{XY}, \tau_{XY})^{-1} K(\sigma_{XY},Y)
\end{equation}
Here $\tau^{XY}$ and $\sigma^{XY}$ are the pivots of the approximation. 
\begin{equation}\label{eq:LR}
    K(\sigma_{XY}, \tau_{XY}) = L_{XY} R_{XY}
\end{equation}
We also compute $L_{XY}$ and $R_{XY}$, the LU factors of $K(\sigma_{XY}, \tau_{XY})$, as stated in Equation~\eqref{eq:LR}.
The LU factors of $K(\sigma_{XY}, \tau_{XY})$ can be obtained as a by-product of the ACA routine and need not be computed separately~\cite{bebendorf2009recompression}. It is to be noted that in the ACA routine, we only store the pivots $\tau^{XY}$, $\sigma^{XY}$ and the matrices $L_{XY}, R_{XY}$ instead of $K(X,\tau_{XY})$, $K(\sigma_{XY}, \tau_{XY})^{-1}$, and $K(\sigma_{XY},Y)$. 

The next task is to perform the matrix-vector product as described in Algorithm~\ref{alg:matvec}. Let the vector to be applied to the matrix be $\psi$. In the matrix-vector product routine, matrices $K(X,\tau_{XY})$ and $K(\sigma_{XY},Y)$ are assembled using the pivots stored in the initialization routine. And the low-rank sub-matrix-vector products are computed as $K(X,\tau_{XY})(R^{-1}_{XY}(L^{-1}_{XY}K(\sigma_{XY},Y)\psi(Y)))$.
\begin{remark}
Note that when we apply the inverse of upper/lower triangular matrices, we are performing backward/forward substitution. The computational complexity of the matrix-vector product using its hierarchical representation with $\log\bkt{N}$ levels computed using the low-rank compression technique of~\cite{H2DKandappan}, where the low-rank factors of the matrix are stored in the initialization phase, is $\mathcal{O}\bkt{pN\log\bkt{N}}$ in both memory and time. Whereas the computational complexity of the matrix-vector product stated in this article, where only the pivots and certain intermediate matrices are stored in the initialization phase, is $\mathcal{O}\bkt{p^2 N+pN\log\bkt{N}}$ in time and $\mathcal{O}\bkt{p^2 N}$ in memory.
\end{remark}

\begin{algorithm}[!htbp]
	\caption{HODLR3D Initialization}\label{alg:initialization}
	\begin{algorithmic}[1]
		\Procedure{InitializeHODLR3D}{$N_{\max}$,$\epsilon$}
		\State{} \Comment{$N_{\max}$ is the maximum number of particles at leaf level;}
		\State{} 
		\State {Form $\mathcal{T}^L$; where $L= \min \left\{l : \abs{\mathcal{C}^{(l)}_{i}} < N_{\max};\forall i \in \{0,1,2,\ldots,8^{l}-1\}\right\}$}
		\State {For each node in the octree identify the Interaction list $\mathcal{I}_\Cl$ and the sets $\mathcal{E}_\Cl$ and $\mathcal{F}_\Cl$.}
\For{\texttt{$l=1:L$}} 
				\For{\texttt{$i=0:8^l-1$}}
					\State $X\gets \text{Index set of } \mathcal{C}_{i}^{(l)}$
					\For{\texttt{j in $\mathcal{I}_{\mathcal{C}_{i}^{(l)}}$}}
						\State $Y\gets \text{Index set of } \mathcal{C}_{j}^{(l)}$
						\State Perform ACA with tolerance $\epsilon$ on the matrix $K(X,Y)$ to identify the pivots $\tau_{XY}$ and $\sigma_{XY}$ and find the LU factorization of the matrix $K(\tau_{XY},\sigma_{XY}) = L_{XY}R_{XY}$.
					\EndFor
				\EndFor
			\EndFor
		\EndProcedure
	\end{algorithmic}
\end{algorithm}
\begin{algorithm}[!htbp]
	\caption{HODLR3D matrix-vector Product $K\psi = b$}\label{alg:matvec}
	\begin{algorithmic}[1]
		\Procedure{MatVec}{$\psi$}
			\For{\texttt{i=0:$8^L-1$}} \Comment{Full-rank Mat-Vec product}
				\State $X\gets \text{Index set of } \mathcal{C}_{i}^{(L)}$
				\State Form the dense matrix $K(X,X)$
				\State $b(X) = b(X) + K(X,X) \times \psi(X)$
				\For{\texttt{j in $\mathcal{E}_{\mathcal{C}_{i}^{(L)}}$}}
					\State $Y\gets \text{Index set of } \mathcal{C}_{j}^{(L)}$
					\State Form the dense matrix $K(X,Y)$
				    \State $b(X) = b(X) + K(X,Y) \times \psi(Y)$
				\EndFor
				\For{\texttt{j in $\mathcal{F}_{\mathcal{C}_{i}^{(L)}}$}}
					\State $Y\gets \text{Index set of } \mathcal{C}_{j}^{(L)}$
					\State Form the dense matrix $K(X,Y)$
				    \State $b(X) = b(X) + K(X,Y) \times \psi(Y)$
				\EndFor
			\EndFor
			\For{\texttt{$l=1:L$}} \Comment{Low-rank Mat-Vec product}
				\For{\texttt{$i=0:8^l-1$}}
					\State $X\gets \text{Index set of } \mathcal{C}_{i}^{(l)}$
					\For{\texttt{j in $\mathcal{I}_{\mathcal{C}_{i}^{(l)}}$}}
						\State $Y\gets \text{Index set of } \mathcal{C}_{j}^{(l)}$
						\State $b(X) = b(X) + K(X,\tau_{XY}) R_{XY}^{-1}L_{XY}^{-1}  K(\sigma_{XY},Y)\psi(Y)$
					\EndFor
				\EndFor
			\EndFor
			\State \textbf{return} $b$
		\EndProcedure
	\end{algorithmic}
\end{algorithm}
\section{Numerical Results} 
\label{sec:numericalResults}
We perform three experiments to demonstrate the performance of HODLR3D. 
In the first experiment, we demonstrate the HODLR3D matrix-vector product, where we consider three kernels i) 3D Laplacian ii) kernel $1/r^{4}$ iii) real part of the 3D Helmholtz kernel with the wave number set to $1$, $\frac{\cos(r)}{r}$. We compare the performance of HODLR3D with that of i) HODLR and ii) a $\mathcal{H}$ matrix that compresses only those interactions between cubes that are well-separated. In the second experiment, we demonstrate HODLR3D accelerated iterative solver for an integral equation. In the third experiment, we demonstrate the parallel scalability of HODLR3D. 

In~\cref{table:ResultsNotations}, we state some notations that we use in this section.
\begin{table}[!htbp]
\centering
\caption{List of notations followed in this section}
\label{table:ResultsNotations}
  \begin{tabularx}{\textwidth}{|p{30mm}|X|}
    \hline
    $N$ & System size that denotes the number of particles in the computational domain.\\ \hline
    $\epsilon$ & tolerance set for ACA routine \\ \hline
    Memory & The memory needed to store the matrix in HODLR3D representation.\\ \hline
    CR & Compression ratio that denotes the ratio of the number of floating point numbers that need to be stored for the HODLR3D representation to $N^2$ \\ \hline
    Initialization Time & It includes the time taken by the Initialize routine and the time taken to form the dense matrices corresponding to the edge-sharing, face-sharing, and self-interactions. \\ \hline
    Matrix-Vector product time & 
    It is the time taken by the Matrix-vector product routine minus the time taken to form the dense matrices corresponding to the edge-sharing, face-sharing, and self-interactions.\\ \hline
    Maximum rank & The maximum rank among all the interactions that were compressed while building the HODLR3D representation.\\ \hline
    Relative error & Relative error in the solution are measured using $\|.\|_{2}$. \\ \hline
  \end{tabularx}

 \end{table}
For all the experiments, we consider the following settings:
\begin{itemize}
    \item $\textib{B}$ is considered to be the cube $[-1,1]^{3}$.
    \item $N_{max}$ is chosen to be $216$.
\end{itemize}
All the experiments were run on an Intel Xeon Gold 2.5GHz processor with 8 OpenMP threads. All the results in this section are reproducible; we direct our readers to refer to \url{https://hodlr3d.readthedocs.io/en/latest/reproducibility.html} for further information.
\subsection{HODLR3D matrix-vector product}
\label{sec:sec5_1} The first numerical experiment demonstrates the scalability of the matrix-vector product through HODLR3D, HODLR and $\mathcal{H}$-matrix with strong admissibility condition. A uniform distribution of particles is considered in the domain $\Bb \subseteq [-1,1]^3$. The matrix entries $A_{ij}$ are generated using the Green's function of the Laplacian in $3$D, i.e., $A_{ij} = \dfrac{1}{\|\vec{r}_i - \vec{r}_j\|_{2}}$, where $\vec{r_{i}} \in \Bb$ $\forall i\in\{1,2,...,N\}$. We consider the vector $\Vec{x}$ to be applied to the matrix to be a random vector. 

In~\cref{fig:oneOverR_2}, we illustrate the scaling of maximum rank, CR, memory, initialization time, matrix-vector product time, and relative error with $N$ for the three kernels. We set $\epsilon$ to $10^{-7}$ for all three algorithms. We have repeated the above numerical experiment for two more kernels viz., $\frac{1}{r^{4}}$ and $\frac{\cos(r)}{r}$. In~\cref{fig:oneOverR_4,fig:cosR_2}, we illustrate the scaling of maximum rank, CR, memory, initialization time, matrix-vector product time, and relative error with $N$ for the kernels $\dfrac{1}{r^4}$ and $\dfrac{\cos\bkt{r}}{r}$.

\begin{figure}[!htbp]
    \centering
    \subfloat[]{\begin{tikzpicture}[every mark/.append style={mark size=1pt}]
	\begin{axis}[xlabel={$N$},
	ylabel={Maximum rank},
		legend style={
                at={(0.8,0.4)},
               anchor=south,
               legend columns=1,
               cells={anchor=east},
              font=\footnotesize,
               rounded corners=1pt,
               },
		xmode = log,
	    ymode = log,
		]
		
		\addplot[
		color=blue,
		mark=square*,
		] coordinates {
(8000,4.600000e+01)
(27000,5.400000e+01)
(64000,5.700000e+01)
(125000,5.900000e+01)
(216000,6.300000e+01)
(343000,6.300000e+01)
(512000,6.100000e+01)
(729000,6.400000e+01)
(1000000,6.400000e+01)
(1331000,6.400000e+01)
(1728000,6.300000e+01)
(2197000,6.700000e+01)
(2744000,6.500000e+01)
(3375000,6.700000e+01)
		};
		\addplot[
		color=red,
		mark=*,
		] coordinates {
(8000,7.500000e+01)
(27000,9.300000e+01)
(64000,1.030000e+02)
(125000,1.050000e+02)
(216000,1.090000e+02)
(343000,1.200000e+02)
(512000,1.210000e+02)
(729000,1.240000e+02)
(1000000,1.370000e+02)
(1331000,1.350000e+02)
(1728000,1.450000e+02)
(2197000,1.500000e+02)
(2744000,1.250000e+02)
(3375000,1.390000e+02)
		};
						\addplot[
		color=black,
		mark=triangle*,
		] coordinates {
(8000,2.540000e+02)
(27000,4.660000e+02)
(64000,7.200000e+02)
(125000,9.370000e+02)
(216000,1.146000e+03)
(343000,1.574000e+03)
(512000,1.871000e+03)
(729000,2.154000e+03)
(1000000,2.720000e+03)
(1331000,3.328000e+03)
		};
		\addplot[mark=none, black, dashed][
		domain = 8000:729000,
		] {(pow(10,0.2)*pow(x,2/3)};
		\legend{H, HODLR3D, HODLR, $\sqrt[3]{N^{2}}$}
	\end{axis}
\end{tikzpicture}}\quad%
    \subfloat[]{\begin{tikzpicture}[every mark/.append style={mark size=1pt}]
	\begin{axis}[xlabel={$N$},
	ylabel={CR},
		legend style={
                at={(0.3,0.1)},
               anchor=south,
               legend columns=1,
               cells={anchor=east},
               font=\footnotesize,
               rounded corners=1pt,
               },
		xmode = log,
	    ymode = log,
		]

		\addplot[
		color=blue,
		mark=square*,
		] coordinates {
(8000,6.944330e-01)
(27000,4.587260e-01)
(64000,2.141660e-01)
(125000,1.704400e-01)
(216000,1.047520e-01)
(343000,6.899990e-02)
(512000,4.741950e-02)
(729000,3.537080e-02)
(1000000,3.303410e-02)
(1331000,2.561130e-02)
(1728000,1.997660e-02)
(2197000,1.601110e-02)
(2744000,1.297700e-02)
(3375000,1.068920e-02)
		};
		\addplot[
		color=red,
		mark=*,
		] coordinates {
(8000,6.518320e-01)
(27000,4.086830e-01)
(64000,1.914400e-01)
(125000,1.467460e-01)
(216000,9.080690e-02)
(343000,5.999840e-02)
(512000,4.134440e-02)
(729000,3.072430e-02)
(1000000,2.813240e-02)
(1331000,2.180530e-02)
(1728000,1.710820e-02)
(2197000,1.375940e-02)
(2744000,1.116830e-02)
(3375000,9.199400e-03)
		};
		
		\addplot[
		color=black,
		mark=triangle*,
		] coordinates {
(8000,4.672320e-01)
(27000,2.373200e-01)
(64000,1.386470e-01)
(125000,9.395540e-02)
(216000,6.753870e-02)
(343000,5.323050e-02)
(512000,4.102820e-02)
(729000,3.269360e-02)
(1000000,2.859190e-02)
(1331000,2.579220e-02)
		};
		\legend{H, HODLR3D, HODLR}
	\end{axis}
\end{tikzpicture}}\quad%
    \subfloat[]{\begin{tikzpicture}[every mark/.append style={mark size=1pt}]
	\begin{axis}[xlabel={$N$},
	ylabel={Memory (GB)},
		legend style={
                at={(0.8,0.04)},
               anchor=south,
               legend columns=1,
               cells={anchor=east},
              font=\footnotesize,
               rounded corners=1pt,
               },
		xmode = log,
	    ymode = log,
		]
		
		\addplot[
		color=blue,
		mark=square*,
		] coordinates {
(8000, 64* 5.555460e-03)
(27000, 64* 4.180140e-02)
(64000, 64* 1.096530e-01)
(125000, 64* 3.328910e-01)
(216000, 64* 6.109150e-01)
(343000, 64* 1.014720e+00)
(512000, 64* 1.553840e+00)
(729000, 64* 2.349690e+00)
(1000000, 64* 4.129270e+00)
(1331000, 64* 5.671500e+00)
(1728000, 64* 7.456230e+00)
(2197000, 64* 9.660290e+00)
(2744000, 64* 1.221380e+01)
(3375000, 64* 1.521960e+01)
		};
		\addplot[
		color=red,
		mark=*,
		] coordinates {
(8000, 64* 5.214650e-03)
(27000, 64* 3.724120e-02)
(64000, 64* 9.801730e-02)
(125000, 64* 2.866140e-01)
(216000, 64* 5.295860e-01)
(343000, 64* 8.823440e-01)
(512000, 64* 1.354770e+00)
(729000, 64* 2.041020e+00)
(1000000, 64* 3.516550e+00)
(1331000, 64* 4.828670e+00)
(1728000, 64* 6.385620e+00)
(2197000, 64* 8.301730e+00)
(2744000, 64* 1.051150e+01)
(3375000, 64* 1.309840e+01)
		};
		\addplot[mark=none, red, dashed][
		domain = 8000:7077888,
		] {(pow(10,-5.5)*x*log10(x)};
		\addplot[
		color=black,
		mark=triangle*,
		] coordinates {
(8000, 64* 3.737850e-03)
(27000, 64* 2.162580e-02)
(64000, 64* 7.098740e-02)
(125000, 64* 1.835070e-01)
(216000, 64* 3.938860e-01)
(343000, 64* 7.828140e-01)
(512000, 64* 1.344410e+00)
(729000, 64* 2.171840e+00)
(1000000, 64* 3.573980e+00)
(1331000, 64* 5.711560e+00)
		};
		\addplot[mark=none, black, dashed][
		domain = 8000:729000,
		] {(pow(10,-7.3)*pow(x,5/3)*log10(x)};
		
		\legend{H, HODLR3D, $N\log(N)$, HODLR, $N^{5/3}\log(N)$}
	\end{axis}
\end{tikzpicture}}\quad%
    \subfloat[]{\begin{tikzpicture}[every mark/.append style={mark size=1pt}]
	\begin{axis}[xlabel={$N$},
	ylabel={Initialization Time (s)},
		legend style={
                at={(0.8,0.04)},
               anchor=south,
               legend columns=1,
               cells={anchor=east},
              font=\footnotesize,
               rounded corners=1pt,
               },
		xmode = log,
	    ymode = log,
		]
	
		\addplot[
		color=blue,
		mark=square*,
		] coordinates {
(8000,6.714460e-01)
(27000,5.915870e+00)
(64000,1.420900e+01)
(125000,5.478910e+01)
(216000,7.704310e+01)
(343000,1.195490e+02)
(512000,1.686490e+02)
(729000,2.440110e+02)
(1000000,6.317040e+02)
(1331000,7.869570e+02)
(1728000,9.749840e+02)
(2197000,1.237550e+03)
(2744000,1.462450e+03)
(3375000,1.731990e+03)
		};
		\addplot[
		color=red,
		mark=*,
		] coordinates {
(8000,3.021340e-01)
(27000,2.883990e+00)
(64000,6.809100e+00)
(125000,2.569230e+01)
(216000,4.590210e+01)
(343000,7.582490e+01)
(512000,1.167270e+02)
(729000,1.779520e+02)
(1000000,3.667210e+02)
(1331000,5.154560e+02)
(1728000,6.930260e+02)
(2197000,9.086790e+02)
(2744000,1.137680e+03)
(3375000,1.435810e+03)
		};
		\addplot[mark=none, red, dashed][
		domain = 8000:7077888,
		] {(pow(10,-3)*x*log10(x)};
				\addplot[
		color=black,
		mark=triangle*,
		] coordinates {
(8000,6.126850e-01)
(27000,7.405470e+00)
(64000,3.957100e+01)
(125000,1.315910e+02)
(216000,3.357520e+02)
(343000,9.797260e+02)
(512000,2.118530e+03)
(729000,3.937140e+03)
(1000000,8.511670e+03)
(1331000,1.678160e+04)
		};
		\addplot[mark=none, black, dashed][
		domain = 8000:1331000,
		] {(pow(10,-8.5)*pow(x,7/3)*log10(x)};

		\legend{H, HODLR3D, $N\log(N)$, HODLR, $N^{7/3}\log(N)$}
	\end{axis}
\end{tikzpicture}}\quad%
    \subfloat[]{\begin{tikzpicture}[every mark/.append style={mark size=1pt}]
	\begin{axis}[xlabel={$N$},
	ylabel={Matrix-Vector product time (s)},
		legend style={
                at={(0.8,0.04)},
               anchor=south,
               legend columns=1,
               cells={anchor=east},
               font=\footnotesize,
               rounded corners=1pt,
               },
		xmode = log,
	    ymode = log,
		]
		
		\addplot[
		color=blue,
		mark=square*,
		] coordinates {
(8000,1.632320e-02)
(27000,6.933940e-02)
(64000,1.595910e-01)
(125000,4.980190e-01)
(216000,6.356110e-01)
(343000,7.836670e-01)
(512000,9.884950e-01)
(729000,1.461920e+00)
(1000000,4.896630e+00)
(1331000,5.929290e+00)
(1728000,6.931390e+00)
(2197000,8.272450e+00)
(2744000,9.439730e+00)
(3375000,1.188010e+01)
		};
		\addplot[
		color=red,
		mark=*,
		] coordinates {

(8000,3.277770e-03)
(27000,3.389720e-02)
(64000,6.175130e-02)
(125000,3.093530e-01)
(216000,4.622850e-01)
(343000,6.186310e-01)
(512000,8.460310e-01)
(729000,1.185320e+00)
(1000000,3.518580e+00)
(1331000,4.353170e+00)
(1728000,5.341960e+00)
(2197000,6.564770e+00)
(2744000,7.775710e+00)
(3375000,9.231310e+00)
		};
\addplot[mark=none, red, dashed][
		domain = 8000:7077888,
		] {(pow(10,-5.5)*x*log10(x)};		
				\addplot[
		color=black,
		mark=triangle*,
		] coordinates {
(8000,2.417450e-03)
(27000,1.591920e-02)
(64000,4.741680e-02)
(125000,1.319790e-01)
(216000,2.656580e-01)
(343000,5.174560e-01)
(512000,9.071160e-01)
(729000,1.428040e+00)
(1000000,2.604100e+00)
(1331000,7.037790e+00)
		};
\addplot[mark=none, black, dashed][
		domain = 8000:729000,
		] {(pow(10,-7.5)*pow(x,5/3)*log10(x)};
		
		\legend{H, HODLR3D, $N\log(N)$, HODLR, $N^{5/3}\log(N)$}
\end{axis}
\end{tikzpicture}}\quad%
    \subfloat[]{\begin{tikzpicture}[every mark/.append style={mark size=1pt}]
	\begin{axis}[xlabel={$N$},
	ylabel={Relative error},
		legend style={
                at={(0.8,0.04)},
               anchor=south,
               legend columns=1,
               cells={anchor=east},
              font=\footnotesize,
               rounded corners=1pt,
               },
		xmode = log,
	    ymode = log,
	    ymin = 1e-16,
	    ymax = 1e-0,
		]

		\addplot[
		color=blue,
		mark=square*,
		] coordinates {
(8000,4.190530e-08)
(27000,8.260650e-09)
(64000,1.650610e-08)
(125000,1.020340e-08)
(216000,7.377760e-09)
(343000,5.629230e-09)
(512000,5.752400e-09)
(729000,5.012510e-09)
(1000000,4.911500e-09)
(1331000,5.951520e-09)
(1728000,6.160010e-09)
(2197000,4.853100e-09)
(2744000,5.274050e-09)
(3375000,5.585000e-09)
		};
		\addplot[
		color=red,
		mark=*,
		] coordinates {
(8000,7.227830e-08)
(27000,1.658720e-08)
(64000,2.605630e-08)
(125000,2.512710e-08)
(216000,3.665400e-08)
(343000,3.330930e-08)
(512000,4.694020e-08)
(729000,5.451460e-08)
(1000000,4.561370e-08)
(1331000,6.466630e-08)
(1728000,1.304910e-07)
(2197000,9.467980e-08)
(2744000,2.163560e-07)
(3375000,1.280640e-07)
		};
		
		\addplot[
		color=black,
		mark=triangle*,
		] coordinates {
(8000,4.531300e-07)
(27000,4.149330e-06)
(64000,6.874090e-06)
(125000,4.021570e-05)
(216000,3.210050e-04)
(343000,2.048430e-03)
(512000,3.470940e-03)
(729000,4.395700e-03)
(1000000,1.063030e-03)
(1331000,4.853190e-03)
		};		
		\legend{H, HODLR3D, HODLR}
	\end{axis}
\end{tikzpicture}}%
    \caption{Various benchmarks of HODLR3D matrix-vector product in comparison with those of HODLR and $\mathcal{H}$ matrix-vector products for the kernel $\frac{1}{r}$ with $\epsilon$ of the three algorithms set to $10^{-7}$} 
    \label{fig:oneOverR_2}
\end{figure}
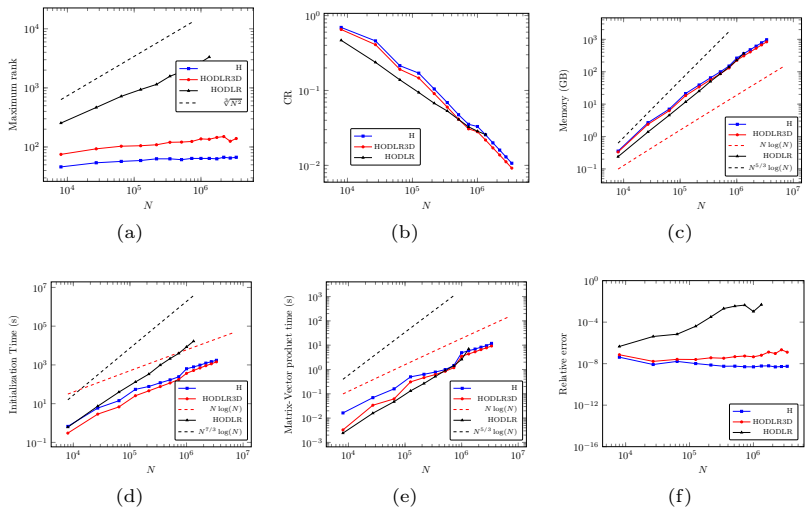
\begin{figure}[!htbp]
    \centering
    \subfloat[]{\begin{tikzpicture}[every mark/.append style={mark size=1pt}]
	\begin{axis}[xlabel={$N$},
	ylabel={Maximum rank},
		legend style={
                at={(0.8,0.4)},
               anchor=south,
               legend columns=1,
               cells={anchor=east},
              font=\footnotesize,
               rounded corners=1pt,
               },
		xmode = log,
	    ymode = log,
		]
		
		\addplot[
		color=blue,
		mark=square*,
		] coordinates {
(8000,6.500000e+01)
(27000,9.800000e+01)
(64000,9.000000e+01)
(125000,1.140000e+02)
(216000,1.230000e+02)
(343000,1.290000e+02)
(512000,1.220000e+02)
(729000,1.340000e+02)
(1000000,1.290000e+02)
		};
		\addplot[
		color=red,
		mark=*,
		] coordinates {
(8000,1.350000e+02)
(27000,1.740000e+02)
(64000,1.870000e+02)
(125000,2.210000e+02)
(216000,2.470000e+02)
(343000,2.480000e+02)
(512000,2.600000e+02)
(729000,2.900000e+02)
(1000000,2.760000e+02)
(1331000,2.520000e+02)
(1728000,3.050000e+02)
		};
\addplot[
		color=black,
		mark=triangle*,
		] coordinates {
(8000,2.500000e+02)
(27000,5.130000e+02)
(64000,9.250000e+02)
(125000,1.400000e+03)
(216000,1.876000e+03)
(343000,2.499000e+03)
(512000,3.433000e+03)
(729000,4.230000e+03)
		};
		\addplot[mark=none, black, dashed][
		domain = 8000:512000,
		] {(pow(10,0.2)*pow(x,2/3)};
		\legend{H, HODLR3D, HODLR, $\sqrt[3]{N^{2}}$}	\end{axis}
\end{tikzpicture}}\quad%
    \subfloat[]{\begin{tikzpicture}[every mark/.append style={mark size=1pt}]
	\begin{axis}[xlabel={$N$},
	ylabel={CR},
		legend style={
                at={(0.3,0.1)},
               anchor=south,
               legend columns=1,
               cells={anchor=east},
               font=\footnotesize,
               rounded corners=1pt,
               },
		xmode = log,
	    ymode = log,
		]

		\addplot[
		color=blue,
		mark=square*,
		] coordinates {
(8000,9.616210e-01)
(27000,6.454190e-01)
(64000,3.228080e-01)
(125000,2.419280e-01)
(216000,1.570840e-01)
(343000,1.071640e-01)
(512000,7.489180e-02)
(729000,5.520650e-02)
(1000000,4.948590e-02)
		};
		\addplot[
		color=red,
		mark=*,
		] coordinates {
(8000,9.030500e-01)
(27000,5.686530e-01)
(64000,2.871140e-01)
(125000,2.073970e-01)
(216000,1.352020e-01)
(343000,9.217040e-02)
(512000,6.475690e-02)
(729000,4.765240e-02)
(1000000,4.184530e-02)
(1331000,3.291040e-02)
(1728000,2.638050e-02)
		};
\addplot[
		color=black,
		mark=triangle*,
		] coordinates {
(8000,4.761110e-01)
(27000,2.702150e-01)
(64000,1.701760e-01)
(125000,1.263950e-01)
(216000,9.544790e-02)
(343000,7.772110e-02)
(512000,6.211540e-02)
(729000,4.972400e-02)
		};
		\legend{H, HODLR3D, HODLR}	\end{axis}
\end{tikzpicture}}\quad%
    \subfloat[]{\begin{tikzpicture}[every mark/.append style={mark size=1pt}]
	\begin{axis}[xlabel={$N$},
	ylabel={Memory (GB)},
		legend style={
                at={(0.8,0.04)},
               anchor=south,
               legend columns=1,
               cells={anchor=east},
              font=\footnotesize,
               rounded corners=1pt,
               },
		xmode = log,
	    ymode = log,
		]
		
		\addplot[
		color=blue,
		mark=square*,
		] coordinates {
(8000, 64* 7.692970e-03)
(27000, 64* 5.881380e-02)
(64000, 64* 1.652780e-01)
(125000, 64* 4.725150e-01)
(216000, 64* 9.161140e-01)
(343000, 64* 1.575970e+00)
(512000, 64* 2.454060e+00)
(729000, 64* 3.667380e+00)
(1000000, 64* 6.185730e+00)
		};
		\addplot[
		color=red,
		mark=*,
		] coordinates {
(8000, 64* 7.224400e-03)
(27000, 64* 5.181850e-02)
(64000, 64* 1.470030e-01)
(125000, 64* 4.050720e-01)
(216000, 64* 7.885000e-01)
(343000, 64* 1.355470e+00)
(512000, 64* 2.121950e+00)
(729000, 64* 3.165550e+00)
(1000000, 64* 5.230670e+00)
(1331000, 64* 7.287840e+00)
(1728000, 64* 9.846460e+00)
		};
		\addplot[mark=none, red, dashed][
		domain = 8000:7077888,
		] {(pow(10,-5.4)*x*log10(x)};
\addplot[
		color=black,
		mark=triangle*,
		] coordinates {
(8000, 64* 3.808890e-03)
(27000, 64* 2.462330e-02)
(64000, 64* 8.712990e-02)
(125000, 64* 2.468660e-01)
(216000, 64* 5.566520e-01)
(343000, 64* 1.142980e+00)
(512000, 64* 2.035400e+00)
(729000, 64* 3.303170e+00)
		};
\addplot[mark=none, black, dashed][
		domain = 8000:512000,
		] {(pow(10,-7.3)*pow(x,5/3)*log10(x)};
		
		\legend{H, HODLR3D, $N\log(N)$, HODLR, $N^{5/3}\log(N)$}	
		\end{axis}
\end{tikzpicture}}\quad%
    \subfloat[]{\begin{tikzpicture}[every mark/.append style={mark size=1pt}]
	\begin{axis}[xlabel={$N$},
	ylabel={Initialization Time (s)},
		legend style={
                at={(0.9,0.04)},
               anchor=south,
               legend columns=1,
               cells={anchor=east},
              font=\footnotesize,
               rounded corners=1pt,
               },
		xmode = log,
	    ymode = log,
		]
		
		\addplot[
		color=blue,
		mark=square*,
		] coordinates {
(8000,4.434360e-01)
(27000,4.318410e+00)
(64000,1.172240e+01)
(125000,4.039690e+01)
(216000,7.602250e+01)
(343000,1.334890e+02)
(512000,2.138280e+02)
(729000,3.284080e+02)
(1000000,6.121000e+02)
		};
		\addplot[
		color=red,
		mark=*,
		] coordinates {
(8000,2.010450e+00)
(27000,4.860330e+00)
(64000,1.558000e+01)
(125000,5.564210e+01)
(216000,1.009340e+02)
(343000,1.959380e+02)
(512000,2.722920e+02)
(729000,4.867980e+02)
(1000000,8.512030e+02)
(1331000,1.157520e+03)
(1728000,2.333450e+03)
		};
		\addplot[mark=none, red, dashed][
		domain = 8000:6859000,
		] {(pow(10,-3.3)*x*log10(x)};
\addplot[
		color=black,
		mark=triangle*,
		] coordinates {
(8000,8.228200e-01)
(27000,1.243880e+01)
(64000,6.506080e+01)
(125000,2.732090e+02)
(216000,8.247620e+02)
(343000,2.401190e+03)
(512000,6.734540e+03)
(729000,1.502860e+04)
		};
\addplot[mark=none, black, dashed][
		domain = 8000:512000,
		] {(pow(10,-7.5)*pow(x,7/3)*log10(x)};
		
		\legend{H, HODLR3D, $N\log(N)$, HODLR, $N^{7/3}\log(N)$}
	\end{axis}
\end{tikzpicture}}\quad%
    \subfloat[]{\begin{tikzpicture}[every mark/.append style={mark size=1pt}]
	\begin{axis}[xlabel={$N$},
	ylabel={Matrix-Vector product time (s)},
		legend style={
                at={(0.8,0.04)},
               anchor=south,
               legend columns=1,
               cells={anchor=east},
              font=\footnotesize,
               rounded corners=1pt,
               },
		xmode = log,
	    ymode = log,
		]
		
		\addplot[
		color=blue,
		mark=square*,
		] coordinates {
(8000,4.689940e-03)
(27000,4.722670e-02)
(64000,9.427520e-02)
(125000,4.270280e-01)
(216000,6.491340e-01)
(343000,9.995530e-01)
(512000,1.404840e+00)
(729000,1.888220e+00)
(1000000,5.073860e+00)
		};
		\addplot[
		color=red,
		mark=*,
		] coordinates {
(8000,5.992550e-03)
(27000,8.779600e-02)
(64000,1.238790e-01)
(125000,5.726860e-01)
(216000,1.087810e+00)
(343000,1.188240e+00)
(512000,2.314390e+00)
(729000,2.383300e+00)
(1000000,7.364930e+00)
(1331000,9.463700e+00)
(1728000,5.335460e+01)
		};
		\addplot[mark=none, red, dashed][
		domain = 8000:4096000,
		] {(pow(10,-5.5)*x*log10(x)};
\addplot[
		color=black,
		mark=triangle*,
		] coordinates {
(8000,2.334330e-03)
(27000,2.408720e-02)
(64000,5.733630e-02)
(125000,1.774650e-01)
(216000,4.360260e-01)
(343000,8.199220e-01)
(512000,1.483020e+00)
(729000,3.632280e+00)
		};
\addplot[mark=none, black, dashed][
		domain = 8000:512000,
		] {(pow(10,-7.2)*pow(x,5/3)*log10(x)};
		
		\legend{H, HODLR3D, $N\log(N)$, HODLR, $N^{5/3}\log(N)$}
		\end{axis}
\end{tikzpicture}}\quad%
    \subfloat[]{\begin{tikzpicture}[every mark/.append style={mark size=1pt}]
	\begin{axis}[xlabel={$N$},
	ylabel={Relative error},
		legend style={
                at={(0.8,0.4)},
               anchor=south,
               legend columns=1,
               cells={anchor=east},
              font=\footnotesize,
               rounded corners=1pt,
               },
		xmode = log,
	    ymode = log,
		]
		
		\addplot[
		color=blue,
		mark=square*,
		] coordinates {
(8000,1.199520e-09)
(27000,1.981540e-08)
(64000,1.308370e-09)
(125000,6.399840e-09)
(216000,1.833530e-08)
(343000,1.188450e-08)
(512000,1.372290e-09)
(729000,3.910950e-09)
(1000000,6.425530e-09)
		};
		\addplot[
		color=red,
		mark=*,
		] coordinates {
(8000,1.147650e-08)
(27000,2.427980e-08)
(64000,3.523860e-08)
(125000,1.886480e-07)
(216000,3.926220e-07)
(343000,4.130920e-07)
(512000,7.718360e-07)
(729000,8.874070e-07)
(1000000,1.656790e-06)
(1331000,1.254110e-06)
(1728000,2.015970e-06)
		};
\addplot[
		color=black,
		mark=triangle*,
		] coordinates {
(8000,1.080530e-04)
(27000,1.757390e-04)
(64000,3.301860e-04)
(125000,4.310950e-04)
(216000,8.307520e-04)
(343000,1.077500e-03)
(512000,1.544970e-03)
(729000,3.248470e-03)
		};
		\legend{H, HODLR3D, HODLR}	
		\end{axis}
\end{tikzpicture}}%
    \caption{Various benchmarks of HODLR3D matrix-vector product in comparison with those of HODLR and $\mathcal{H}$ matrix-vector products for the kernel $\frac{1}{r^4}$ with $\epsilon$ of the three algorithms set to $10^{-7}$} 
    \label{fig:oneOverR_4}
\end{figure}
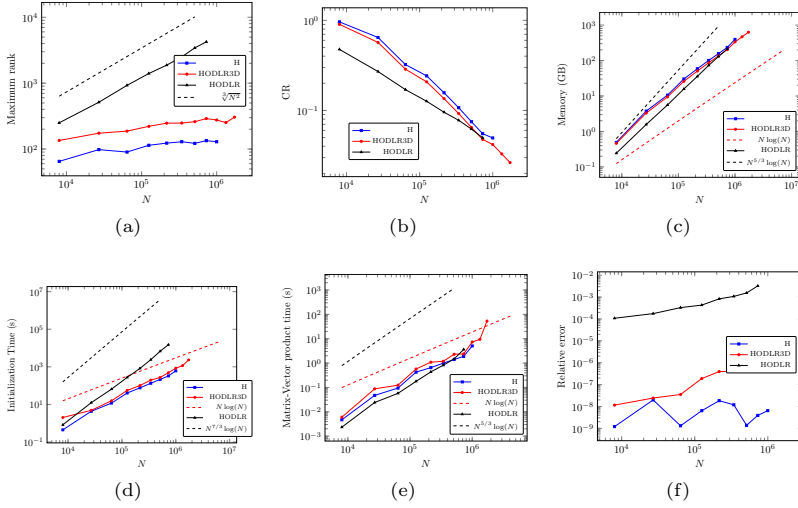
\begin{figure}[!htbp]
    \centering
    \subfloat[]{\begin{tikzpicture}[every mark/.append style={mark size=1pt}]
	\begin{axis}[xlabel={$N$},
	ylabel={Maximum rank},
		legend style={
                at={(0.8,0.4)},
               anchor=south,
               legend columns=1,
               cells={anchor=east},
              font=\footnotesize,
               rounded corners=1pt,
               },
		xmode = log,
	    ymode = log,
		]
		
		\addplot[
		color=blue,
		mark=square*,
		] coordinates {
(8000,4.700000e+01)
(27000,5.900000e+01)
(64000,5.900000e+01)
(125000,6.300000e+01)
(216000,6.200000e+01)
(343000,6.600000e+01)
(512000,6.300000e+01)
(729000,6.800000e+01)
(1000000,6.300000e+01)
(1331000,6.600000e+01)
(1728000,6.200000e+01)
(2197000,6.700000e+01)
(2744000,6.700000e+01)
(3375000,7.000000e+01)
(4096000,6.700000e+01)
		};
		\addplot[
		color=red,
		mark=*,
		] coordinates {
(8000,8.700000e+01)
(27000,9.600000e+01)
(64000,1.110000e+02)
(125000,1.250000e+02)
(216000,1.200000e+02)
(343000,1.330000e+02)
(512000,1.380000e+02)
(729000,1.440000e+02)
(1000000,1.450000e+02)
(1331000,1.430000e+02)
(1728000,1.490000e+02)
(2197000,1.560000e+02)
(2744000,1.400000e+02)
(3375000,1.550000e+02)
(4096000,1.520000e+02)
		};
\addplot[
		color=black,
		mark=triangle*,
		] coordinates {
(8000,2.450000e+02)
(27000,4.930000e+02)
(64000,7.790000e+02)
(125000,9.810000e+02)
(216000,1.397000e+03)
(343000,1.644000e+03)
(512000,1.970000e+03)
(729000,2.512000e+03)
(1000000,2.829000e+03)
(1331000,3.449000e+03)
		};
		\addplot[mark=none, black, dashed][
		domain = 8000:729000,
		] {(pow(10,0)*pow(x,2/3)};
		\legend{H, HODLR3D, HODLR, $\sqrt[3]{N^{2}}$}
		\end{axis}
\end{tikzpicture}}\quad%
    \subfloat[]{\begin{tikzpicture}[every mark/.append style={mark size=1pt}]
	\begin{axis}[	xlabel={$N$},
	ylabel={CR},
		legend style={
                at={(0.3,0.1)},
               anchor=south,
               legend columns=1,
               cells={anchor=east},
               font=\footnotesize,
               rounded corners=1pt,
               },
		xmode = log,
	    ymode = log,
		]
		
		\addplot[
		color=blue,
		mark=square*,
		] coordinates {
(8000,7.114030e-01)
(27000,4.723270e-01)
(64000,2.199700e-01)
(125000,1.753810e-01)
(216000,1.075590e-01)
(343000,7.067230e-02)
(512000,4.870770e-02)
(729000,3.620860e-02)
(1000000,3.361260e-02)
(1331000,2.602040e-02)
(1728000,2.029770e-02)
(2197000,1.620950e-02)
(2744000,1.320570e-02)
(3375000,1.087280e-02)
(4096000,8.982900e-03)
		};
		\addplot[
		color=red,
		mark=*,
		] coordinates {
(8000,6.700290e-01)
(27000,4.204250e-01)
(64000,1.969170e-01)
(125000,1.506700e-01)
(216000,9.321190e-02)
(343000,6.137730e-02)
(512000,4.253440e-02)
(729000,3.136940e-02)
(1000000,2.858430e-02)
(1331000,2.215510e-02)
(1728000,1.738790e-02)
(2197000,1.392930e-02)
(2744000,1.136070e-02)
(3375000,9.367930e-03)
(4096000,7.762380e-03)
		};
\addplot[
		color=black,
		mark=triangle*,
		] coordinates {
(8000,4.723330e-01)
(27000,2.555450e-01)
(64000,1.412660e-01)
(125000,1.019190e-01)
(216000,7.379730e-02)
(343000,5.430090e-02)
(512000,4.349610e-02)
(729000,3.567540e-02)
(1000000,3.037780e-02)
(1331000,2.662550e-02)
		};
		\legend{H, HODLR3D, HODLR}
		\end{axis}
\end{tikzpicture}}\quad%
    \subfloat[]{\begin{tikzpicture}[every mark/.append style={mark size=1pt}]
	\begin{axis}[xlabel={$N$},
	ylabel={Memory (GB)},
		legend style={
                at={(0.8,0.04)},
               anchor=south,
               legend columns=1,
               cells={anchor=east},
              font=\footnotesize,
               rounded corners=1pt,
               },
		xmode = log,
	    ymode = log,
		]
		
		\addplot[
		color=blue,
		mark=square*,
		] coordinates {
(8000, 64* 5.691220e-03)
(27000, 64* 4.304080e-02)
(64000, 64* 1.126250e-01)
(125000, 64* 3.425410e-01)
(216000, 64* 6.272830e-01)
(343000, 64* 1.039320e+00)
(512000, 64* 1.596050e+00)
(729000, 64* 2.405340e+00)
(1000000, 64* 4.201580e+00)
(1331000, 64* 5.762090e+00)
(1728000, 64* 7.576070e+00)
(2197000, 64* 9.780010e+00)
(2744000, 64* 1.242910e+01)
(3375000, 64* 1.548100e+01)
(4096000, 64* 1.883850e+01)
};
		\addplot[
		color=red,
		mark=*,
		] coordinates {
(8000, 64* 5.360230e-03)
(27000, 64* 3.831130e-02)
(64000, 64* 1.008210e-01)
(125000, 64* 2.942770e-01)
(216000, 64* 5.436120e-01)
(343000, 64* 9.026220e-01)
(512000, 64* 1.393770e+00)
(729000, 64* 2.083880e+00)
(1000000, 64* 3.573040e+00)
(1331000, 64* 4.906140e+00)
(1728000, 64* 6.489990e+00)
(2197000, 64* 8.404230e+00)
(2744000, 64* 1.069260e+01)
(3375000, 64* 1.333830e+01)
(4096000, 64* 1.627890e+01)
		};
		\addplot[mark=none, red, dashed][
		domain = 8000:7077888,
		] {(pow(10,-5.5)*x*log10(x)};
\addplot[
		color=black,
		mark=triangle*,
		] coordinates {
(8000, 64* 3.778660e-03)
(27000, 64* 2.328650e-02)
(64000, 64* 7.232840e-02)
(125000, 64* 1.990600e-01)
(216000, 64* 4.303860e-01)
(343000, 64* 7.985550e-01)
(512000, 64* 1.425280e+00)
(729000, 64* 2.369920e+00)
(1000000, 64* 3.797220e+00)
(1331000, 64* 5.896080e+00)
		};
\addplot[mark=none, black, dashed][
		domain = 8000:729000,
		] {(pow(10,-7.3)*pow(x,5/3)*log10(x)};
		\legend{H, HODLR3D, $N\log(N)$, HODLR, $N^{5/3}\log(N)$}	
		\end{axis}
\end{tikzpicture}}\quad%
    \subfloat[]{\begin{tikzpicture}[every mark/.append style={mark size=1pt}]
	\begin{axis}[	xlabel={$N$},
	ylabel={Initialization Time (s)},
		legend style={
                at={(0.8,0.04)},
               anchor=south,
               legend columns=1,
               cells={anchor=east},
              font=\footnotesize,
               rounded corners=1pt,
               },
		xmode = log,
	    ymode = log,
		]
		
		\addplot[
		color=blue,
		mark=square*,
		] coordinates {

(8000,6.793450e-01)
(27000,7.283170e+00)
(64000,1.609360e+01)
(125000,6.384510e+01)
(216000,1.083300e+02)
(343000,1.693780e+02)
(512000,2.435560e+02)
(729000,3.483120e+02)
(1000000,8.265870e+02)
(1331000,1.120190e+03)
(1728000,1.512240e+03)
(2197000,1.924360e+03)
(2744000,2.438660e+03)
(3375000,2.959680e+03)
(4096000,3.604970e+03)
		};
		\addplot[
		color=red,
		mark=*,
		] coordinates {
(8000,1.088490e+00)
(27000,7.936960e+00)
(64000,1.878240e+01)
(125000,6.618000e+01)
(216000,1.097770e+02)
(343000,2.058400e+02)
(512000,3.000250e+02)
(729000,5.664160e+02)
(1000000,1.051450e+03)
(1331000,1.155440e+03)
(1728000,1.889940e+03)
(2197000,2.127270e+03)
(2744000,2.771280e+03)
(3375000,3.780140e+03)
(4096000,3.442880e+03)
		};
		\addplot[mark=none, red, dashed][
		domain = 8000:7077888,
		] {(pow(10,-3.7)*x*log10(x)};
\addplot[
		color=black,
		mark=triangle*,
		] coordinates {
(8000,1.615620e+00)
(27000,1.601090e+01)
(64000,6.774760e+01)
(125000,1.898470e+02)
(216000,5.304100e+02)
(343000,1.196450e+03)
(512000,2.559040e+03)
(729000,5.769670e+03)
(1000000,1.006440e+04)
(1331000,1.957770e+04)
		};
\addplot[mark=none, black, dashed][
		domain = 8000:729000,
		] {(pow(10,-8.6)*pow(x,7/3)*log10(x)};
		\legend{H, HODLR3D, $N\log(N)$, HODLR, $N^{7/3}\log(N)$}
		
		\end{axis}
\end{tikzpicture}}\quad%
    \subfloat[]{\begin{tikzpicture}[every mark/.append style={mark size=1pt}]
	\begin{axis}[xlabel={$N$},
	ylabel={Matrix-Vector product time (s)},
		legend style={
                at={(0.8,0.04)},
               anchor=south,
               legend columns=1,
               cells={anchor=east},
              font=\footnotesize,
               rounded corners=1pt,
               },
		xmode = log,
	    ymode = log,
		]
		
		\addplot[
		color=blue,
		mark=square*,
		] coordinates {
(8000,4.757980e-03)
(27000,5.749370e-02)
(64000,9.992730e-02)
(125000,5.747720e-01)
(216000,7.534080e-01)
(343000,1.138220e+00)
(512000,1.422670e+00)
(729000,2.034450e+00)
(1000000,6.724870e+00)
(1331000,8.120280e+00)
(1728000,9.704670e+00)
(2197000,1.230760e+01)
(2744000,1.458090e+01)
(3375000,1.684800e+01)
(4096000,1.908340e+01)
		};
		\addplot[
		color=red,
		mark=*,
		] coordinates {
(8000,8.773200e-03)
(27000,4.552450e-02)
(64000,8.662700e-02)
(125000,4.833160e-01)
(216000,7.015930e-01)
(343000,9.399370e-01)
(512000,1.610010e+00)
(729000,1.954570e+00)
(1000000,5.430500e+00)
(1331000,7.298090e+00)
(1728000,8.000190e+00)
(2197000,7.761070e+00)
(2744000,1.118030e+01)
(3375000,1.246610e+01)
(4096000,1.681350e+01)
		};
		\addplot[mark=none, red, dashed][
		domain = 8000:7077888,
		] {(pow(10,-5.5)*x*log10(x)};
\addplot[
		color=black,
		mark=triangle*,
		] coordinates {
(8000,2.669190e-03)
(27000,6.885180e-02)
(64000,1.351910e-01)
(125000,1.371670e-01)
(216000,3.330310e-01)
(343000,5.345130e-01)
(512000,1.024120e+00)
(729000,2.140480e+00)
(1000000,3.720680e+00)
(1331000,6.202880e+00)
		};
\addplot[mark=none, black, dashed][
		domain = 8000:729000,
		] {(pow(10,-7.5)*pow(x,5/3)*log10(x)};
		\legend{H, HODLR3D, $N\log(N)$, HODLR, $N^{5/3}\log(N)$}
		\end{axis}
\end{tikzpicture}}\quad%
    \subfloat[]{\begin{tikzpicture}[every mark/.append style={mark size=1pt}]
	\begin{axis}[	xlabel={$N$},
	ylabel={Relative error},
		legend style={
                at={(0.8,0.04)},
               anchor=south,
               legend columns=1,
               cells={anchor=east},
              font=\footnotesize,
               rounded corners=1pt,
               },
		xmode = log,
	    ymode = log,
	    ymin = 1e-16,
	    ymax = 1e-0,
		]
		
		\addplot[
		color=blue,
		mark=square*,
		] coordinates {
(8000,1.698700e-08)
(27000,1.391860e-08)
(64000,1.866130e-08)
(125000,1.451070e-08)
(216000,1.080680e-08)
(343000,9.144300e-09)
(512000,7.522630e-09)
(729000,9.717860e-09)
(1000000,7.166050e-09)
(1331000,7.767820e-09)
(1728000,6.009470e-09)
(2197000,6.370800e-09)
(2744000,4.507280e-09)
(3375000,7.829810e-09)
(4096000,5.640310e-09)
		};
		\addplot[
		color=red,
		mark=*,
		] coordinates {
(8000,3.760970e-08)
(27000,2.691660e-08)
(64000,3.486200e-08)
(125000,4.229370e-08)
(216000,1.864620e-08)
(343000,3.674880e-08)
(512000,5.535820e-08)
(729000,6.707400e-08)
(1000000,9.333490e-08)
(1331000,1.242860e-07)
(1728000,1.820600e-07)
(2197000,7.458330e-08)
(2744000,1.033520e-07)
(3375000,2.446950e-07)
(4096000,2.235300e-07)
		};
\addplot[
		color=black,
		mark=triangle*,
		] coordinates {
(8000,5.227730e-07)
(27000,3.304780e-06)
(64000,1.352050e-05)
(125000,1.383050e-05)
(216000,2.876520e-05)
(343000,1.665050e-04)
(512000,2.535310e-03)
(729000,1.864210e-03)
(1000000,2.047450e-03)
(1331000,1.618230e-03)
		};
		\legend{H, HODLR3D, HODLR}
		\end{axis}
\end{tikzpicture}}%
    \caption{Various benchmarks of HODLR3D matrix-vector product in comparison with those of HODLR and $\mathcal{H}$ matrix-vector products for the kernel $\frac{\cos\bkt{r}}{r}$ with $\epsilon$ of the three algorithms set to $10^{-7}$} 
    \label{fig:cosR_2}
\end{figure}
\subsection{HODLR3D accelerated iterative solver for integral equations in 3D}
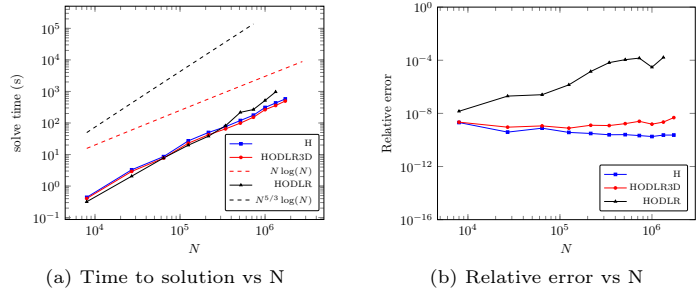
\begin{figure}[!htbp]
    \centering
    \subfloat[Time to solution vs N]{\begin{tikzpicture}[every mark/.append style={mark size=1pt}]
	\begin{axis}[xlabel={$N$},
	ylabel={solve time (s)},
		legend style={
                at={(0.8,0.04)},
               anchor=south,
               legend columns=1,
               cells={anchor=east},
              font=\footnotesize,
               rounded corners=1pt,
               },
		xmode = log,
	    ymode = log,
		]
		
		\addplot[
		color=blue,
		mark=square*,
		] coordinates {
(8000,0.439511)
(27000,3.280980)
(64000,8.583250)
(125000,27.109500)
(216000,50.037600)
(343000,76.351000)
(512000,119.691000)
(729000,176.542000)
(1000000,310.767000)
(1331000,431.025000)
(1728000,577.755000)
		};
		\addplot[
		color=red,
		mark=*,
		] coordinates {
(8000,4.045770e-01)
(27000,2.911540e+00)
(64000,7.698440e+00)
(125000,2.284190e+01)
(216000,4.230780e+01)
(343000,6.505390e+01)
(512000,9.912750e+01)
(729000,1.521170e+02)
(1000000,2.641030e+02)
(1331000,3.589520e+02)
(1728000,4.963620e+02)
		};
		\addplot[mark=none, red, dashed][
		domain = 8000:2744000,
		] {(pow(10,-3.3)*x*log10(x)};
\addplot[
		color=black,
		mark=triangle*,
		] coordinates {
(8000,3.188200e-01)
(27000,2.056960e+00)
(64000,7.764010e+00)
(125000,1.977550e+01)
(216000,3.813890e+01)
(343000,8.196370e+01)
(512000,2.174850e+02)
(729000,2.673730e+02)
(1000000,5.188950e+02)
(1331000,9.683570e+02)
		};
\addplot[mark=none, black, dashed][
		domain = 8000:729000,
		] {(pow(10,-5.4)*pow(x,5/3)*log10(x)};
		\legend{H, HODLR3D, $N\log(N)$, HODLR, $N^{5/3}\log(N)$}
		\end{axis}
\end{tikzpicture}}\qquad%
    \subfloat[Relative error vs N]{\begin{tikzpicture}[every mark/.append style={mark size=1pt}]
	\begin{axis}[	xlabel={$N$},
	ylabel={Relative error},
		legend style={
                at={(0.8,0.04)},
               anchor=south,
               legend columns=1,
               cells={anchor=east},
              font=\footnotesize,
               rounded corners=1pt,
               },
		xmode = log,
	    ymode = log,
	    ymin = 1e-16,
	    ymax = 1e-0,
		]
		
		\addplot[
		color=blue,
		mark=square*,
		] coordinates {
(8000,2.052220e-09)
(27000,3.892780e-10)
(64000,7.673900e-10)
(125000,3.696760e-10)
(216000,3.074250e-10)
(343000,2.458860e-10)
(512000,2.512590e-10)
(729000,2.136640e-10)
(1000000,1.810420e-10)
(1331000,2.326170e-10)
(1728000,2.336750e-10)
		};
		\addplot[
		color=red,
		mark=*,
		] coordinates {
(8000,2.190290e-09)
(27000,9.167730e-10)
(64000,1.141590e-09)
(125000,7.707800e-10)
(216000,1.275440e-09)
(343000,1.211290e-09)
(512000,1.689690e-09)
(729000,2.528760e-09)
(1000000,1.545610e-09)
(1331000,2.215180e-09)
(1728000,4.801720e-09)
		};
\addplot[
		color=black,
		mark=triangle*,
		] coordinates {
(8000,1.445170e-08)
(27000,2.001610e-07)
(64000,2.550850e-07)
(125000,1.405670e-06)
(216000,1.425070e-05)
(343000,6.648470e-05)
(512000,1.106860e-04)
(729000,1.445630e-04)
(1000000,2.998490e-05)
(1331000,1.638370e-04)
		};
		\legend{H, HODLR3D, HODLR}
		\end{axis}
\end{tikzpicture}}%
    \caption{Results for the HODLR3D accelerated iterative solver for the integral equation of equation~\ref{eq:IE} with $\epsilon$ of the three algorithms set to $10^{-7}$} 
    \label{fig:solvePlots2}
\end{figure}

We solve the Fredholm integral equation of the second kind, equation~\eqref{eq:IE}, with $\bm{K}(x,y)=\frac{1}{\|x-y\|_{2}}$, where $x,y\in \mathbb{R}^3$.
\begin{equation}\label{eq:IE}
    \sigma(x)+\int_{\textib{B}} \bm{K}(x,y)\sigma(y)dy = f(x), \qquad \Bb\subset \mathbb{R}^3
\end{equation}
We use a piecewise constant collocation method with collocation points on a uniform tensor grid  on $3$D and an appropriate quadrature to integrate the singularity. The Fredholm integral equation is thus discretised to obtain a linear system of the form
\begin{equation}\label{eq:ILDiscrete}
    A\vec{\sigma}=\vec{f}
\end{equation}
We solve for $\vec{\sigma}$ using GMRES~\cite{barrett1994templates,saad1986gmres}, an iterative technique. Each iteration of GMRES involves computing a matrix-vector product. We employ HODLR3D to perform this computation to accelerate the solver.
In order to find the error in the solution, we consider a random vector $\vec{\sigma}$ and compute the right-hand side, $\vec{f}$, exactly up to roundoff. We then use $\vec{f}$ as the right-hand side of equation~\eqref{eq:ILDiscrete}, to solve the equation. Let the computed $\vec{\sigma}$ be $\vec{\sigma_{c}}$. We compute the relative forward error in 2-norm sense $\frac{\|\vec{\sigma_{c}}-\vec{\sigma}\|_{2}}{\|\vec{\sigma}\|_{2}}$. We terminate the GMRES routine when the relative residual, $\frac{\|A\vec{\sigma_{c}}-\vec{f}\|_{2}}{\|\vec{f}\|_{2}}$, is less than $10^{-10}$. In~\cref{fig:solvePlots2}, we illustrate the scaling of solve time and relative error with $N$, wherein $\epsilon$ for ACA of the three algorithms is set to $10^{-7}$.

\subsection{Inferences}
It is to be observed from~\cref{fig:oneOverR_2,fig:oneOverR_4,fig:cosR_2} that the maximum rank for HODLR scales roughly with $\sqrt[3]{N^{2}}$. The initialization time scales roughly as $\mathcal{O}(N^{7/3}\log(N))$ and the memory, matrix-vector product time, and solve time scale roughly as $\mathcal{O}(N^{5/3}\log(N))$. In contrast, the maximum ranks for both HODLR3D and $\mathcal{H}$ matrices are almost constant and do not scale with any power of $N$. The complexities of memory, initialization time, the matrix-vector product time and the solve time for both HODLR3D and $\mathcal{H}$-matrix are almost linear. It is also to be observed that the memory and timing performance of HODLR3D is almost the same as that of $\mathcal{H}$-matrix. Based on the results shown in~\Cref{fig:oneOverR_2,fig:oneOverR_4,fig:cosR_2,fig:solvePlots2}, HODLR3D is an attractive alternative to $\mathcal{H}$-matrix with $\eta=\sqrt3$ for large $N$-body problems.

\section{Parallel HODLR3D}
\label{sec:parH3D}
   This section describes the parallel HODLR3D matrix-vector product algorithm (\cref{alg:matvec}) on a distributed memory system.
   We refer the readers~\cite {izadi2012hierarchical,li2020distributed}
   for the literature on the parallelization of the matrix-vector product of hierarchical matrices. In the parallel HODLR3D matrix-vector product, we try to maintain its parallel efficiency almost constant irrespective of the system size. We achieve this by considering each node in the hierarchical tree $\mathcal{T}^L$ as a data-independent computational unit and the processors perform an almost equal amount of computations at each level of the hierarchical tree. The computational unit that we refer to here is in context to the sub-matrix vector product concerned with a particular node in the hierarchical tree in~\cref{alg:matvec}. We also assume that the computations involved with respect to a node at a particular level in the hierarchical tree are the same with other nodes in that level.  
   
   We consider the number of parallel processors $n_p$ as $2^k,$ with $k>0$. For a level $l$ in $\mathcal{T}^L$, there are $8^l$ computational units. If $8^l\geq 2^k$, then each processor has $\ceil{\dfrac{8^l}{n_p}}$ computational units.
   If $8^l < 2^k$ then, each computational unit is shared by $\ceil{\dfrac{n_p}{8^l}}$ processors. 
   For the computational units at a level $l$, where $8^l < 2^k$, the computations are shared by a group of processors. Hence, we need an explicit work-sharing routine. This is performed as follows: For the dense matrix-vector product as in line numbers 5,9 and 14 of~\cref{alg:matvec}, we divide the columns of the matrix across the processors and perform the matrix-vector product. If $r$ is the rank of the low-rank representation in line $22$ of~\cref{alg:matvec}, then we perform the low-rank matrix-vector product through the following steps.
   \begin{itemize}
        \item Row-wise parallel matrix-vector product of $z_1=K(\sigma_{XY},Y)z$, i.e., the rows of the matrix $K(\sigma_{XY},Y)$ are equally shared among the processors. This involves a computational cost of $\mathcal{O}\bkt{\frac{rN}{n_p}}$.
        \item Communicate among other MPI processes and perform $z_2 = R^{-1}_{XY}\bkt{L^{-1}_{XY}z_1}$ which involves communication cost of $\mathcal{O}\bkt{n_pr}$, and a computational cost of $\mathcal{O}\bkt{r^2}$
        \item Perform  $K(X,\tau_{XY})z_2$, which involves a computational cost of $\mathcal{O}\bkt{\frac{rN}{n_p}}$
    \end{itemize}

    In summary, based on the level in the hierarchical tree, the computational unit is either processed by a single processor or a group of processors. On the whole, the HODLR3D-based matrix-vector product is equivalent to row-wise sharing of the matrix among the processors and the communication step just involves the MPI Gather operation among all the processors. 
    
    To illustrate the scalability of parallel HODLR3D, we look at the same numerical experiment in~\Cref{sec:sec5_1} (with Laplacian kernel i.e., $\frac{1}{\magn{x-y}_2}\quad x,y \in \mathbb{R}^3$ and $\epsilon=10^{-6}$). The parallel HODLR3D algorithm is developed in C++ with OpenMPI 4.14. Each MPI process runs on a single core of  Intel Xeon Gold 6248, a 2.5 GHz processor. We vary the system size $N$ and the number of MPI process $n_p$ and tabulate the maximum and average time taken to perform HODLR3D matrix-vector product in~\cref{tab:par_mat_vec}. From~\cref{tab:par_mat_vec}, we infer that for $n_p\leq 8$, we have good parallel efficiency as at each level the computational units are independent; Whereas for $n_p > 8$, the parallel efficiency drops as $n_p$ increases since we have more levels where the computational units are shared by a group of processors. For more details on the process of initialising HODLR3D structure in the distributed memory system refer to~\Cref {sec:par_init}. 

\begin{table}[!htbp]
    \centering
    \caption{Parallel HODLR3D Matrix-Vector product}
\label{tab:par_mat_vec}
    \begin{tabular}{@{}ccccccc@{}}
    \toprule
                                        &                                     & \multicolumn{5}{c}{\textbf{Time taken for}} \\ \cmidrule(l){3-7} 
                                 &                               & \multicolumn{3}{c}{\textbf{MatVec per process}} & \multicolumn{2}{c}{\textbf{MPI Collective Communication}} \\ \cmidrule(l){3-7} 
    \multirow{-3}{*}{\textbf{N}} & \multirow{-3}{*}{\textbf{np}} & \textbf{Avg}  & \textbf{Max} & \textbf{Speedup} & \textbf{Avg}                & \textbf{Max}                \\ \midrule
                                        & \textbf{2}  & 0.290    & 0.290    & 1  & 0.004  & 0.008   \\
                                        & \textbf{4}                          & 0.131    & 0.140    & 2  & 0.352  & 0.717   \\
                                        & \textbf{8}                          & 0.072    & 0.078    & 4  & 0.798  & 0.964   \\
                                        & \textbf{16}                         & 0.037    & 0.046    & 6  & 0.021  & 0.030   \\
                                        & \textbf{32} & 0.022    & 0.034    & 8  & 0.085  & 0.144   \\
    \multirow{-6}{*}{\textbf{125000}}   & \textbf{64} & 0.009    & 0.013    & 23 & 0.008  & 0.009   \\ \midrule
                                        & \textbf{2}  & 4.418    & 4.434    & 1  & 0.035  & 0.065   \\
                                        & \textbf{4}                          & 2.092    & 2.192    & 2  & 3.824  & 8.041   \\
                                        & \textbf{8}                          & 1.083    & 1.163    & 4  & 1.832  & 3.160   \\
                                        & \textbf{16}                         & 0.532    & 0.608    & 7  & 0.459  & 0.482   \\
                                        & \textbf{32} & 0.302    & 0.409    & 11 & 0.097  & 0.154   \\
    \multirow{-6}{*}{\textbf{1000000}}  & \textbf{64} & 0.142    & 0.201    & 22 & 0.032  & 0.045   \\ \midrule
                                        & \textbf{2}  & 28.560   & 28.879   & 1  & 0.758  & 1.499   \\
                                        & \textbf{4}                          & 13.993   & 14.903   & 2  & 13.419 & 28.002  \\
                                        & \textbf{8}                          & 6.940    & 7.847    & 4  & 4.231  & 13.093  \\
                                        & \textbf{16}                         & 3.138    & 3.725    & 8  & 0.287  & 0.311   \\
                                        & \textbf{32} & 1.848    & 2.405    & 12 & 0.106  & 0.265   \\
    \multirow{-6}{*}{\textbf{3375000}}  & \textbf{64} & 0.921    & 1.361    & 21 & 0.094  & 0.131   \\ \midrule
                                        & \textbf{2}  & 62.914   & 63.257   & 1  & 1.113  & 2.184   \\
                                        & \textbf{4}                          & 30.203   & 31.106   & 2  & 37.588 & 74.930  \\
                                        & \textbf{8}                          & 15.851   & 18.092   & 3  & 16.925 & 43.202  \\
                                        & \textbf{16}                         & 7.150    & 8.136    & 8  & 3.897  & 3.975   \\
                                        & \textbf{32} & 4.001    & 4.993    & 13 & 0.214  & 0.345   \\
    \multirow{-6}{*}{\textbf{8000000}}  & \textbf{64} & 2.482    & 4.349    & 15 & 0.401  & 0.495   \\ \midrule
                                        & \textbf{2}  & 151.724  & 152.148  & 1  & 1.249  & 2.414   \\
                                        & \textbf{4}                          & 72.447   & 73.242   & 2  & 62.934 & 128.980 \\
                                        & \textbf{8}                          & 36.441   & 40.670   & 4  & 29.887 & 72.478  \\
                                        & \textbf{16}                         & 17.627   & 19.924   & 8  & 7.720  & 8.017   \\
                                        & \textbf{32} & 9.715    & 11.630   & 13 & 0.269  & 0.378   \\
    \multirow{-6}{*}{\textbf{15625000}} & \textbf{64} & 5.193    & 7.145    & 21 & 0.517  & 0.735   \\ \bottomrule
    \end{tabular}
    \end{table}
\section{Conclusions}
In this article, we introduced HODLR3D, a new class of hierarchical matrices for problems arising in three dimensions. We proved upper bounds for the ranks of different off-diagonal blocks of the matrix when the underlying kernel function is the 3D Laplacian kernel. The main highlight of the theorem is that the ranks of the vertex-sharing interactions do not scale with any power of $N$. Based on this observation, the HODLR3D matrix is introduced, whose construction and matrix-vector product algorithm scale almost linearly.

We further compared the performances of HODLR, HODLR3D and a $\mathcal{H}$ matrix in computing matrix-vector product and solving an integral equation. It can be observed from the numerical results that HODLR3D performs better than HODLR, and HODLR3D can be considered competitive to the $\mathcal{H}$ matrices with strong admissibility condition.

Like HODLR2D~\cite{H2DKandappan} and HODLR3D, it is possible to extend the idea of compressing vertex-sharing blocks to higher dimensions. For example, consider $n$-dimensional hypercube which can be partitioned using $2^{n}$-tree. The neighbors for a hypercube are those that share a $1$-cube (vertex), $2$-cube (edge), $3$-cube (face), ..., $(n-2)$-cube, $(n-1)$-cube with it. A similar observation, as in Section~\ref{sec:rankGrowth}, holds true in $n$-dimensions too, i.e., the rank of vertex sharing interactions scale as $\mathcal{O}(\log(N)\log^{n}(\log(N)))$ and that of the well-separated interactions do not scale with $N$. And the ranks of higher order-cube sharing interactions, say $n'$-cube sharing interactions ($1\leq n'\leq n-1$), scale with $\mathcal{O}(N^{n'/n}\log^{n}(N))$. For the proof, we refer the readers to~\cite{khan2022numerical}. An almost linear complexity HODLR$n$D algorithm can thus be constructed by compressing only the well-separated and vertex-sharing interactions. The algorithm for computing matrix-vector product is similar to the one described in Section~\ref{sec:HODLR3D}.

In addition to the fast matrix-vector products, one can construct a fast direct solver similar to the Inverse Fast Multipole Method (IFMM)~\cite{ambikasaran2014inverse} by leveraging the off-diagonal low-rank structure of HODLR3D. Since the rank of off-diagonal blocks does not scale with the system size,  the direct solver based on HODLR3D will be a promising alternative to other direct solvers available~\cite{SA_FDS_2013,gujjula2023algebraic,ambikasaran2014inverse}.
\section{Declarations}
\subsection*{Availability of data and materials}
The code used to generate the results are available at the following Repository.
\begin{itemize}
    \item Code Repository: \url{https://github.com/SAFRAN-LAB/HODLR3D}
    \item Documentation on Reproducibility of results:\url{https://hodlr3d.readthedocs.io/en/latest/reproducibility.html} 
\end{itemize}
\subsection*{Funding}
Vaishnavi Gujjula acknowledges the support of Women Leading IITM (India) 2022 in Mathematics (SB22230053MAIITM008880). Sivaram Ambikasaran acknowledges the support of Young Scientist Research Award from Board of Research in Nuclear Sciences, Department of Atomic Energy, India (No.34/20/03/2017-BRNS/34278) and MATRICS grant from the Science and Engineering Research Board, India (Sanction number: MTR/2019/001241).
\subsection*{Acknowledgments}
The authors acknowledge HPCE, IIT Madras for providing access to the AQUA cluster. The authors would like to thank Ritesh Khan for his valuable comments on the draft of this article.

\begin{appendices}
\section{Numerical Experiment on HODLR3D matrix-vector product}
In this section, we repeat the numerical experiment in~\Cref{ssec:HODLR3DRepres} for the different hierarchical structures considered viz., HODLR, HODLR3D and $\mathcal{H}$ matrix such that in matrix-vector product the forward relative error is of the same order. The intention of this numerical experiment is to understand the performance and scalability of different hierarchical structures for the same matrix-vector product forward relative error. We make sure that the relative forward error of the three algorithms that we compare, HODLR3D, HODLR, and $\mathcal{H}$ matrix, are nearly equal so that the rest of the benchmarks can be compared and an inference can be made. To achieve this, we use different values of $\epsilon$ in the ACA routine of the three hierarchical structures, in the range of $10^{-6}-10^{-10}$. We perform this incrementally and record various benchmarks of the hierarchical structures, such that they have a forward relative error of the same order. The kernels that we use to perform the numerical experiment are 
\begin{itemize}
    \item Green's function for Laplace equation in 3D  which is $\dfrac1{r}$
    \item $\dfrac1{r^4}$
    \item Real part of the Green's function for Helmholtz equation in 3D, which is $\dfrac{\cos\bkt{r}}{r}$
\end{itemize}
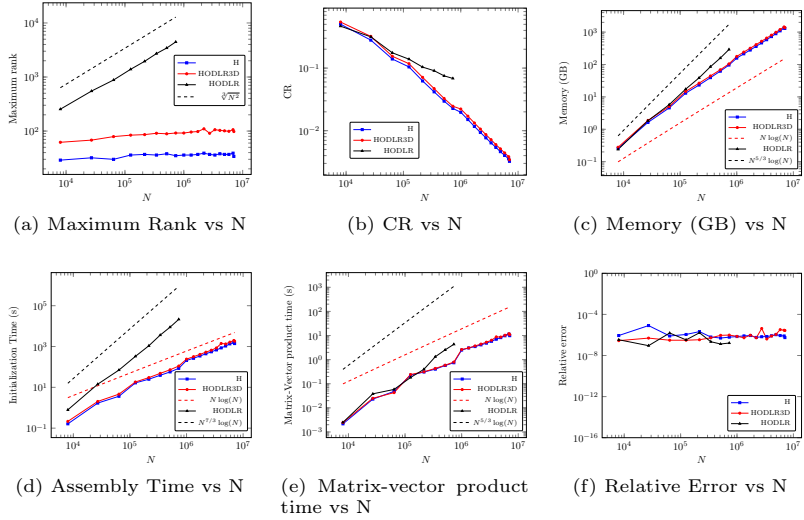
\begin{figure}[!htbp]
    \centering
    \subfloat[Maximum Rank vs N]{\begin{tikzpicture}[every mark/.append style={mark size=1pt}]
	\begin{axis}[xlabel={$N$},
	ylabel={Maximum rank},
		legend style={
                at={(0.8,0.4)},
               anchor=south,
               legend columns=1,
               cells={anchor=east},
              font=\footnotesize,
               rounded corners=1pt,
               },
		xmode = log,
	    ymode = log,
		]
		
		\addplot[
		color=blue,
		mark=square*,
		] coordinates {
(8000,29.000000)
(27000,32.000000)
(64000,30.000000)
(125000,36.000000)
(216000,37.000000)
(343000,36.000000)
(512000,38.000000)
(729000,35.000000)
(1000000,36.000000)
(1331000,36.000000)
(1728000,37.000000)
(2197000,39.000000)
(2744000,37.000000)
(3375000,36.000000)
(4096000,38.000000)
(4913000,37.000000)
(5832000,37.000000)
(6859000,39.000000)
(7077888,34.000000)
		};
		\addplot[
		color=red,
		mark=*,
		] coordinates {
(8000,62.000000)
(27000,68.000000)
(64000,79.000000)
(125000,84.000000)
(216000,86.000000)
(343000,91.000000)
(512000,89.000000)
(729000,92.000000)
(1000000,92.000000)
(1331000,96.000000)
(1728000,98.000000)
(2197000,110.000000)
(2744000,91.000000)
(3375000,106.000000)
(4096000,103.000000)
(4913000,101.000000)
(5832000,99.000000)
(6859000,106.000000)
(7077888,98.000000)
		};
						\addplot[
		color=black,
		mark=triangle*,
		] coordinates {
(8000,254.000000)
(27000,549.000000)
(64000,880.000000)
(125000,1390.000000)
(216000,1951.000000)
(343000,2713.000000)
(512000,3438.000000)
(729000,4459.000000)
		};
		\addplot[mark=none, black, dashed][
		domain = 8000:729000,
		] {(pow(10,0.2)*pow(x,2/3)};
		\legend{H, HODLR3D, HODLR, $\sqrt[3]{N^{2}}$}
	\end{axis}
\end{tikzpicture}}\quad%
    \subfloat[CR vs N]{\begin{tikzpicture}[every mark/.append style={mark size=1pt}]
	\begin{axis}[xlabel={$N$},
	ylabel={CR},
		legend style={
                at={(0.3,0.1)},
               anchor=south,
               legend columns=1,
               cells={anchor=east},
               font=\footnotesize,
               rounded corners=1pt,
               },
		xmode = log,
	    ymode = log,
		]

		\addplot[
		color=blue,
		mark=square*,
		] coordinates {
(8000,0.502707)
(27000,0.280057)
(64000,0.139726)
(125000,0.104397)
(216000,0.062265)
(343000,0.041743)
(512000,0.029662)
(729000,0.022654)
(1000000,0.019699)
(1331000,0.015056)
(1728000,0.011699)
(2197000,0.009397)
(2744000,0.007687)
(3375000,0.006385)
(4096000,0.005408)
(4913000,0.004623)
(5832000,0.004021)
(6859000,0.003538)
(7077888,0.003243)
		};
		\addplot[
		color=red,
		mark=*,
		] coordinates {
(8000,0.544637)
(27000,0.321098)
(64000,0.154500)
(125000,0.117509)
(216000,0.071118)
(343000,0.047071)
(512000,0.032757)
(729000,0.024578)
(1000000,0.022203)
(1331000,0.017068)
(1728000,0.013292)
(2197000,0.010621)
(2744000,0.008660)
(3375000,0.007187)
(4096000,0.005995)
(4913000,0.005126)
(5832000,0.004415)
(6859000,0.003841)
(7077888,0.003444)
		};
		
		\addplot[
		color=black,
		mark=triangle*,
		] coordinates {
(8000,0.467232)
(27000,0.315538)
(64000,0.175540)
(125000,0.138782)
(216000,0.103787)
(343000,0.090862)
(512000,0.075243)
(729000,0.068362)
		};
		\legend{H, HODLR3D, HODLR}
	\end{axis}
\end{tikzpicture}}\quad%
    \subfloat[Memory (GB) vs N]{\begin{tikzpicture}[every mark/.append style={mark size=1pt}]
	\begin{axis}[xlabel={$N$},
	ylabel={Memory (GB)},
		legend style={
                at={(0.8,0.04)},
               anchor=south,
               legend columns=1,
               cells={anchor=east},
              font=\footnotesize,
               rounded corners=1pt,
               },
		xmode = log,
	    ymode = log,
		]
		
		\addplot[
		color=blue,
		mark=square*,
		] coordinates {
(8000, 64* 0.004022)
(27000, 64* 0.025520)
(64000, 64* 0.071540)
(125000, 64* 0.203900)
(216000, 64* 0.363127)
(343000, 64* 0.613879)
(512000, 64* 0.971961)
(729000, 64* 1.504880)
(1000000, 64* 2.462360)
(1331000, 64* 3.334160)
(1728000, 64* 4.366680)
(2197000, 64* 5.669890)
(2744000, 64* 7.235060)
(3375000, 64* 9.091370)
(4096000, 64* 11.340700)
(4913000, 64* 13.947600)
(5832000, 64* 17.095400)
(6859000, 64* 20.806000)
(7077888, 64* 20.304900)
		};
		\addplot[
		color=red,
		mark=*,
		] coordinates {
(8000, 64* 0.004357)
(27000, 64* 0.029260)
(64000, 64* 0.079104)
(125000, 64* 0.229510)
(216000, 64* 0.414762)
(343000, 64* 0.692229)
(512000, 64* 1.073370)
(729000, 64* 1.632720)
(1000000, 64* 2.775380)
(1331000, 64* 3.779570)
(1728000, 64* 4.961260)
(2197000, 64* 6.408010)
(2744000, 64* 8.150790)
(3375000, 64* 10.233400)
(4096000, 64* 12.572900)
(4913000, 64* 15.465900)
(5832000, 64* 18.769600)
(6859000, 64* 22.585300)
(7077888, 64* 21.563500)
		};
		\addplot[mark=none, red, dashed][
		domain = 8000:7077888,
		] {(pow(10,-5.5)*x*log10(x)};
		\addplot[
		color=black,
		mark=triangle*,
		] coordinates {
(8000, 64* 0.003738)
(27000, 64* 0.028753)
(64000, 64* 0.089876)
(125000, 64* 0.271060)
(216000, 64* 0.605285)
(343000, 64* 1.336230)
(512000, 64* 2.465570)
(729000, 64* 4.541270)
		};
		\addplot[mark=none, black, dashed][
		domain = 8000:729000,
		] {(pow(10,-7.3)*pow(x,5/3)*log10(x)};
		
		\legend{H, HODLR3D, $N\log(N)$, HODLR, $N^{5/3}\log(N)$}
	\end{axis}
\end{tikzpicture}}\quad%
    \subfloat[Assembly Time vs N]{\begin{tikzpicture}[every mark/.append style={mark size=1pt}]
	\begin{axis}[xlabel={$N$},
	ylabel={Initialization Time (s)},
		legend style={
                at={(0.8,0.04)},
               anchor=south,
               legend columns=1,
               cells={anchor=east},
              font=\footnotesize,
               rounded corners=1pt,
               },
		xmode = log,
	    ymode = log,
		]
	
		\addplot[
		color=blue,
		mark=square*,
		] coordinates {
(8000,0.162569)
(27000,1.694250)
(64000,3.651780)
(125000,16.551400)
(216000,24.990800)
(343000,38.992400)
(512000,56.731200)
(729000,85.300200)
(1000000,208.652000)
(1331000,265.847000)
(1728000,343.271000)
(2197000,439.119000)
(2744000,550.053000)
(3375000,679.747000)
(4096000,855.731000)
(4913000,1058.700000)
(5832000,1343.950000)
(6859000,1567.290000)
(7077888,1393.570000)
		};
		\addplot[
		color=red,
		mark=*,
		] coordinates {
(8000,0.215109)
(27000,2.036520)
(64000,4.640080)
(125000,18.308200)
(216000,30.006300)
(343000,48.329400)
(512000,72.403200)
(729000,110.581000)
(1000000,241.125000)
(1331000,318.527000)
(1728000,410.406000)
(2197000,529.091000)
(2744000,660.748000)
(3375000,850.421000)
(4096000,1410.420000)
(4913000,1336.740000)
(5832000,1638.500000)
(6859000,1977.900000)
(7077888,1815.840000)
		};
		\addplot[mark=none, red, dashed][
		domain = 8000:7077888,
		] {(pow(10,-4)*x*log10(x)};
				\addplot[
		color=black,
		mark=triangle*,
		] coordinates {
(8000,0.791514)
(27000,14.297000)
(64000,70.286400)
(125000,329.610000)
(216000,1089.170000)
(343000,3634.420000)
(512000,8869.180000)
(729000,21490.700000)
		};
		\addplot[mark=none, black, dashed][
		domain = 8000:729000,
		] {(pow(10,-8.5)*pow(x,7/3)*log10(x)};

		\legend{H, HODLR3D, $N\log(N)$, HODLR, $N^{7/3}\log(N)$}
	\end{axis}
\end{tikzpicture}}\quad%
    \subfloat[Matrix-vector product time vs N]{\begin{tikzpicture}[every mark/.append style={mark size=1pt}]
	\begin{axis}[xlabel={$N$},
	ylabel={Matrix-Vector product time (s)},
		legend style={
                at={(0.8,0.04)},
               anchor=south,
               legend columns=1,
               cells={anchor=east},
               font=\footnotesize,
               rounded corners=1pt,
               },
		xmode = log,
	    ymode = log,
		]
		
		\addplot[
		color=blue,
		mark=square*,
		] coordinates {
(8000,0.002173)
(27000,0.023090)
(64000,0.048179)
(125000,0.226921)
(216000,0.306543)
(343000,0.403285)
(512000,0.572150)
(729000,0.742048)
(1000000,2.497260)
(1331000,3.058610)
(1728000,3.477820)
(2197000,4.094780)
(2744000,4.834670)
(3375000,5.783720)
(4096000,7.049000)
(4913000,8.172140)
(5832000,9.744460)
(6859000,11.481300)
(7077888,9.850940)
		};
		\addplot[
		color=red,
		mark=*,
		] coordinates {

(8000,0.002473)
(27000,0.024943)
(64000,0.042809)
(125000,0.241087)
(216000,0.313825)
(343000,0.432521)
(512000,0.589152)
(729000,0.793872)
(1000000,2.624690)
(1331000,3.047470)
(1728000,3.704850)
(2197000,4.415450)
(2744000,5.281030)
(3375000,6.385070)
(4096000,8.797410)
(4913000,8.673140)
(5832000,10.282400)
(6859000,11.893100)
(7077888,10.903200)
		};
\addplot[mark=none, red, dashed][
		domain = 8000:7077888,
		] {(pow(10,-5.5)*x*log10(x)};		
				\addplot[
		color=black,
		mark=triangle*,
		] coordinates {
(8000,0.002499)
(27000,0.038193)
(64000,0.058186)
(125000,0.181204)
(216000,0.398232)
(343000,1.322740)
(512000,2.547690)
(729000,4.328180)
		};
\addplot[mark=none, black, dashed][
		domain = 8000:729000,
		] {(pow(10,-7.5)*pow(x,5/3)*log10(x)};
		
		\legend{H, HODLR3D, $N\log(N)$, HODLR, $N^{5/3}\log(N)$}
\end{axis}
\end{tikzpicture}}\quad%
    \subfloat[Relative Error vs N]{\begin{tikzpicture}[every mark/.append style={mark size=1pt}]
	\begin{axis}[xlabel={$N$},
	ylabel={Relative error},
		legend style={
                at={(0.8,0.04)},
               anchor=south,
               legend columns=1,
               cells={anchor=east},
              font=\footnotesize,
               rounded corners=1pt,
               },
		xmode = log,
	    ymode = log,
	    ymin = 1e-16,
	    ymax = 1e-0,
		]

		\addplot[
		color=blue,
		mark=square*,
		] coordinates {
(8000,8.781170e-07)
(27000,8.181320e-06)
(64000,7.946040e-07)
(125000,1.113950e-06)
(216000,2.113240e-06)
(343000,6.333840e-07)
(512000,5.006640e-07)
(729000,6.115090e-07)
(1000000,6.981040e-07)
(1331000,8.089660e-07)
(1728000,8.762710e-07)
(2197000,5.748790e-07)
(2744000,6.699870e-07)
(3375000,7.094760e-07)
(4096000,8.286130e-07)
(4913000,1.061750e-06)
(5832000,8.717860e-07)
(6859000,7.981120e-07)
(7077888,5.521490e-07)
		};
		\addplot[
		color=red,
		mark=*,
		] coordinates {
(8000,2.791760e-07)
(27000,4.940570e-07)
(64000,3.074780e-07)
(125000,3.070430e-07)
(216000,3.419200e-07)
(343000,5.240070e-07)
(512000,9.188420e-07)
(729000,9.728950e-07)
(1000000,7.314820e-07)
(1331000,5.318200e-07)
(1728000,9.487780e-07)
(2197000,5.292880e-07)
(2744000,4.261560e-06)
(3375000,4.043310e-07)
(4096000,7.805670e-07)
(4913000,1.150470e-06)
(5832000,3.285690e-06)
(6859000,2.886940e-06)
(7077888,2.725030e-06)
		};
		
		\addplot[
		color=black,
		mark=triangle*,
		] coordinates {
(8000,3.412190e-07)
(27000,9.071260e-08)
(64000,1.543320e-06)
(125000,3.322150e-07)
(216000,1.630430e-06)
(343000,2.205170e-07)
(512000,1.318990e-07)
(729000,1.639850e-07)
		};		
		\legend{H, HODLR3D, HODLR}
	\end{axis}
\end{tikzpicture}}%
    \caption{Various benchmarks of HODLR3D matrix-vector product in comparison with those of HODLR and $\mathcal{H}$ matrix-vector products for the kernel $\frac{1}{r}$ with the relative forward errors of the three algorithms to be of the same order}
    \label{fig:oneOverR}
\end{figure}
\begin{figure}[!htbp]
    \centering
    \subfloat[Maximum rank vs N]{\begin{tikzpicture}[every mark/.append style={mark size=1pt}]
	\begin{axis}[xlabel={$N$},
	ylabel={Maximum rank},
		legend style={
                at={(0.8,0.4)},
               anchor=south,
               legend columns=1,
               cells={anchor=east},
              font=\footnotesize,
               rounded corners=1pt,
               },
		xmode = log,
	    ymode = log,
		]
		
		\addplot[
		color=blue,
		mark=square*,
		] coordinates {
(8000,43.000000)
(27000,50.000000)
(64000,49.000000)
(125000,62.000000)
(216000,56.000000)
(343000,64.000000)
(512000,62.000000)
(729000,65.000000)
(1000000,72.000000)
(1331000,65.000000)
(1728000,66.000000)
		};
		\addplot[
		color=red,
		mark=*,
		] coordinates {
(8000,107.000000)
(27000,135.000000)
(64000,138.000000)
(125000,148.000000)
(216000,173.000000)
(343000,183.000000)
(512000,206.000000)
(729000,191.000000)
(1000000,220.000000)
(1331000,222.000000)
(1728000,208.000000)
		};
\addplot[
		color=black,
		mark=triangle*,
		] coordinates {
(8000,381.000000)
(27000,906.000000)
(64000,1482.000000)
(125000,2317.000000)
(216000,3536.000000)
(343000,4634.000000)
(512000,6000.000000)
		};
		\addplot[mark=none, black, dashed][
		domain = 8000:512000,
		] {(pow(10,0.2)*pow(x,2/3)};
		\legend{H, HODLR3D, HODLR, $\sqrt[3]{N^{2}}$}	\end{axis}
\end{tikzpicture}}\quad%
    \subfloat[CR vs N]{\begin{tikzpicture}[every mark/.append style={mark size=1pt}]
	\begin{axis}[xlabel={$N$},
	ylabel={CR},
		legend style={
                at={(0.3,0.1)},
               anchor=south,
               legend columns=1,
               cells={anchor=east},
               font=\footnotesize,
               rounded corners=1pt,
               },
		xmode = log,
	    ymode = log,
		]

		\addplot[
		color=blue,
		mark=square*,
		] coordinates {
(8000,0.632226)
(27000,0.384119)
(64000,0.190866)
(125000,0.145229)
(216000,0.088868)
(343000,0.059626)
(512000,0.029894)
(729000,0.023235)
(1000000,0.019757)
(1331000,0.015216)
(1728000,0.011771)
(2197000,0.009486)
(2744000,0.007843)
(3375000,0.006558)
(4096000,0.005491)
		};
		\addplot[
		color=red,
		mark=*,
		] coordinates {
(8000,0.719382)
(27000,0.448290)
(64000,0.220833)
(125000,0.163754)
(216000,0.102847)
(343000,0.069620)
(512000,0.048761)
(729000,0.035969)
(1000000,0.032270)
(1331000,0.024819)
(1728000,0.019779)
(4096000,0.009739)
		};
\addplot[
		color=black,
		mark=triangle*,
		] coordinates {
(8000,0.765451)
(27000,0.466623)
(64000,0.301861)
(125000,0.248619)
(216000,0.195355)
(343000,0.151355)
(512000,0.123864)
		};
		\legend{H, HODLR3D, HODLR}	\end{axis}
\end{tikzpicture}}\quad%
    \subfloat[Memory (GB) vs N]{\begin{tikzpicture}[every mark/.append style={mark size=1pt}]
	\begin{axis}[xlabel={$N$},
	ylabel={Memory (GB)},
		legend style={
                at={(0.8,0.04)},
               anchor=south,
               legend columns=1,
               cells={anchor=east},
              font=\footnotesize,
               rounded corners=1pt,
               },
		xmode = log,
	    ymode = log,
		]
		
		\addplot[
		color=blue,
		mark=square*,
		] coordinates {
(8000, 64* 0.005058)
(27000, 64* 0.035003)
(64000, 64* 0.097724)
(125000, 64* 0.283650)
(216000, 64* 0.518278)
(343000, 64* 0.876862)
(512000, 64* 0.979559)
(729000, 64* 1.543490)
(1000000, 64* 2.469610)
(1331000, 64* 3.369570)
(1728000, 64* 4.393640)
(2197000, 64* 5.723660)
(2744000, 64* 7.381680)
(3375000, 64* 9.338010)
(4096000, 64* 11.514600)
		};
		\addplot[
		color=red,
		mark=*,
		] coordinates {
(8000, 64* 0.005755)
(27000, 64* 0.040850)
(64000, 64* 0.113066)
(125000, 64* 0.319833)
(216000, 64* 0.599801)
(343000, 64* 1.023840)
(512000, 64* 1.597820)
(729000, 64* 2.389400)
(1000000, 64* 4.033750)
(1331000, 64* 5.495990)
(1728000, 64* 7.382610)
(4096000, 64* 20.424200)
		};
		\addplot[mark=none, red, dashed][
		domain = 8000:7077888,
		] {(pow(10,-5.4)*x*log10(x)};
\addplot[
		color=black,
		mark=triangle*,
		] coordinates {
(8000,0.006124)
(27000,0.042521)
(64000,0.154553)
(125000,0.485584)
(216000,1.139310)
(343000,2.225850)
(512000,4.058760)
		};
\addplot[mark=none, black, dashed][
		domain = 8000:512000,
		] {(pow(10,-7.3)*pow(x,5/3)*log10(x)};
		
		\legend{H, HODLR3D, $N\log(N)$, HODLR, $N^{5/3}\log(N)$}	
		\end{axis}
\end{tikzpicture}}\quad%
    \subfloat[Assembly time vs N]{\begin{tikzpicture}[every mark/.append style={mark size=1pt}]
	\begin{axis}[xlabel={$N$},
	ylabel={Initialization Time (s)},
		legend style={
                at={(0.9,0.04)},
               anchor=south,
               legend columns=1,
               cells={anchor=east},
              font=\footnotesize,
               rounded corners=1pt,
               },
		xmode = log,
	    ymode = log,
		]
		
		\addplot[
		color=blue,
		mark=square*,
		] coordinates {
(8000,0.226666)
(27000,2.350580)
(64000,5.458470)
(125000,22.698100)
(216000,36.655000)
(343000,59.784600)
(512000,54.966600)
(729000,85.866900)
(1000000,202.607000)
(1331000,262.623000)
(1728000,326.372000)
(2197000,411.335000)
(2744000,517.289000)
(3375000,656.802000)
(4096000,835.521000)
		};
		\addplot[
		color=red,
		mark=*,
		] coordinates {
(8000,0.673134)
(27000,5.831420)
(64000,14.133800)
(125000,51.310500)
(216000,88.180800)
(343000,143.210000)
(512000,208.924000)
(729000,323.843000)
(1000000,632.140000)
(1331000,830.274000)
(1728000,1113.140000)
(4096000,1679.210000)
		};
		\addplot[mark=none, red, dashed][
		domain = 8000:6859000,
		] {(pow(10,-3.3)*x*log10(x)};
\addplot[
		color=black,
		mark=triangle*,
		] coordinates {
(8000,2.043330)
(27000,34.796200)
(64000,203.021000)
(125000,1003.040000)
(216000,3780.380000)
(343000,9918.900000)
(512000,25528.700000)
		};
\addplot[mark=none, black, dashed][
		domain = 8000:512000,
		] {(pow(10,-7.5)*pow(x,7/3)*log10(x)};
		
		\legend{H, HODLR3D, $N\log(N)$, HODLR, $N^{7/3}\log(N)$}
	\end{axis}
\end{tikzpicture}}\quad%
    \subfloat[Matrix-vector product time vs N]{\begin{tikzpicture}[every mark/.append style={mark size=1pt}]
	\begin{axis}[xlabel={$N$},
	ylabel={Matrix-Vector product time (s)},
		legend style={
                at={(0.8,0.04)},
               anchor=south,
               legend columns=1,
               cells={anchor=east},
              font=\footnotesize,
               rounded corners=1pt,
               },
		xmode = log,
	    ymode = log,
		]
		
		\addplot[
		color=blue,
		mark=square*,
		] coordinates {
(8000,0.002957)
(27000,0.030034)
(64000,0.053410)
(125000,0.319039)
(216000,0.411752)
(343000,0.558148)
(512000,0.533417)
(729000,0.749918)
(1000000,2.443930)
(1331000,2.926650)
(1728000,3.400300)
(2197000,4.038470)
(2744000,4.807550)
(3375000,5.718990)
(4096000,6.685590)
		};
		\addplot[
		color=red,
		mark=*,
		] coordinates {
(8000,0.022766)
(27000,0.061886)
(64000,0.120075)
(125000,0.545283)
(216000,0.653549)
(343000,0.928369)
(512000,1.208730)
(729000,1.638360)
(1000000,4.497310)
(1331000,5.595460)
(1728000,6.293340)
(4096000,14.260800)
		};
		\addplot[mark=none, red, dashed][
		domain = 8000:4096000,
		] {(pow(10,-5.5)*x*log10(x)};
\addplot[
		color=black,
		mark=triangle*,
		] coordinates {
(8000,0.004544)
(27000,0.031384)
(64000,0.104389)
(125000,0.362827)
(216000,0.893305)
(343000,1.566290)
(512000,3.222250)
		};
\addplot[mark=none, black, dashed][
		domain = 8000:512000,
		] {(pow(10,-7.2)*pow(x,5/3)*log10(x)};
		
		\legend{H, HODLR3D, $N\log(N)$, HODLR, $N^{5/3}\log(N)$}
		\end{axis}
\end{tikzpicture}}\quad%
    \subfloat[Relative error vs N]{\begin{tikzpicture}[every mark/.append style={mark size=1pt}]
	\begin{axis}[xlabel={$N$},
	ylabel={Relative error},
		legend style={
                at={(0.8,0.04)},
               anchor=south,
               legend columns=1,
               cells={anchor=east},
              font=\footnotesize,
               rounded corners=1pt,
               },
		xmode = log,
	    ymode = log,
	    ymin = 1e-16,
	    ymax = 1e-0,
		]
		
		\addplot[
		color=blue,
		mark=square*,
		] coordinates {
(8000,5.307450e-08)
(27000,7.596810e-07)
(64000,1.137690e-07)
(125000,7.521850e-07)
(216000,9.368660e-07)
(343000,8.278650e-07)
(512000,2.182760e-06)
(729000,8.070650e-06)
(1000000,1.259530e-05)
(1331000,2.574730e-05)
(1728000,6.146410e-05)
(2197000,8.353590e-05)
(2744000,6.280750e-05)
(3375000,5.040420e-05)
(4096000,4.532280e-06)
		};
		\addplot[
		color=red,
		mark=*,
		] coordinates {
(8000,9.439170e-08)
(27000,4.660330e-07)
(64000,7.414970e-07)
(125000,3.284480e-06)
(216000,4.389440e-06)
(343000,5.836740e-06)
(512000,8.697790e-06)
(729000,1.101590e-05)
(1000000,1.629840e-05)
(1331000,1.623180e-05)
(1728000,2.034470e-05)
(4096000,3.981880e-06)
		};
\addplot[
		color=black,
		mark=triangle*,
		] coordinates {
(8000,4.430180e-08)
(27000,1.057320e-05)
(64000,7.974380e-07)
(125000,1.653310e-07)
(216000,3.592410e-07)
(343000,8.092150e-07)
(512000,2.406340e-05)
		};
		\legend{H, HODLR3D, HODLR}	\end{axis}
\end{tikzpicture}}%
    \caption{Various benchmarks of HODLR3D matrix-vector product in comparison with those of HODLR and $\mathcal{H}$ matrix-vector products for the kernel $\frac{1}{r^4}$ with the relative forward errors of the three algorithms to be of the same order}
    \label{fig:oneOverR4}
\end{figure}
\begin{figure}[!htbp]
    \centering
    \subfloat[Maximum rank vs N]{\begin{tikzpicture}[every mark/.append style={mark size=1pt}]
	\begin{axis}[xlabel={$N$},
	ylabel={Maximum rank},
		legend style={
                at={(0.8,0.4)},
               anchor=south,
               legend columns=1,
               cells={anchor=east},
              font=\footnotesize,
               rounded corners=1pt,
               },
		xmode = log,
	    ymode = log,
		]
		
		\addplot[
		color=blue,
		mark=square*,
		] coordinates {
(8000,39.000000)
(27000,44.000000)
(64000,47.000000)
(125000,49.000000)
(216000,50.000000)
(343000,51.000000)
(512000,51.000000)
(729000,38.000000)
(1000000,38.000000)
(1331000,43.000000)
(1728000,37.000000)
(2197000,42.000000)
(2744000,38.000000)
(3375000,43.000000)
(4096000,40.000000)
(4913000,42.000000)
(5832000,37.000000)
(6859000,41.000000)
		};
		\addplot[
		color=red,
		mark=*,
		] coordinates {
(8000,67.000000)
(27000,76.000000)
(64000,87.000000)
(125000,92.000000)
(216000,97.000000)
(343000,104.000000)
(512000,98.000000)
(729000,104.000000)
(1000000,105.000000)
(1331000,100.000000)
(1728000,115.000000)
(2197000,115.000000)
(2744000,110.000000)
(3375000,113.000000)
(4096000,105.000000)
(4913000,122.000000)
(5832000,128.000000)
(6859000,125.000000)
		};
\addplot[
		color=black,
		mark=triangle*,
		] coordinates {
(8000,245.000000)
(27000,493.000000)
(64000,910.000000)
(125000,1343.000000)
(216000,1961.000000)
(343000,2578.000000)
(512000,3374.000000)
(729000,4294.000000)
		};
		\addplot[mark=none, black, dashed][
		domain = 8000:729000,
		] {(pow(10,0)*pow(x,2/3)};
		\legend{H, HODLR3D, HODLR, $\sqrt[3]{N^{2}}$}
		\end{axis}
\end{tikzpicture}}\quad%
    \subfloat[CR vs N]{\begin{tikzpicture}[every mark/.append style={mark size=1pt}]
	\begin{axis}[	xlabel={$N$},
	ylabel={CR},
		legend style={
                at={(0.3,0.1)},
               anchor=south,
               legend columns=1,
               cells={anchor=east},
               font=\footnotesize,
               rounded corners=1pt,
               },
		xmode = log,
	    ymode = log,
		]
		
		\addplot[
		color=blue,
		mark=square*,
		] coordinates {
(8000,0.617191)
(27000,0.378431)
(64000,0.180640)
(125000,0.140613)
(216000,0.084602)
(343000,0.055738)
(512000,0.038920)
(729000,0.023204)
(1000000,0.020078)
(1331000,0.015422)
(1728000,0.012008)
(2197000,0.009590)
(2744000,0.007856)
(3375000,0.006544)
(4096000,0.005493)
(4913000,0.004722)
(5832000,0.004102)
(6859000,0.003639)

		};
		\addplot[
		color=red,
		mark=*,
		] coordinates {
(8000,0.574767)
(27000,0.337191)
(64000,0.160442)
(125000,0.121011)
(216000,0.073221)
(343000,0.048365)
(512000,0.033773)
(729000,0.025194)
(1000000,0.022643)
(1331000,0.017435)
(1728000,0.013534)
(2197000,0.010858)
(2744000,0.008857)
(3375000,0.007318)
(4096000,0.006128)
(4913000,0.005227)
(5832000,0.004517)
(6859000,0.003920)
		};
\addplot[
		color=black,
		mark=triangle*,
		] coordinates {
(8000,0.472333)
(27000,0.255545)
(64000,0.182242)
(125000,0.129209)
(216000,0.108863)
(343000,0.084138)
(512000,0.071956)
(729000,0.061018)
		};
		\legend{H, HODLR3D, HODLR}
		\end{axis}
\end{tikzpicture}}\quad%
    \subfloat[Memory (GB) vs N]{\begin{tikzpicture}[every mark/.append style={mark size=1pt}]
	\begin{axis}[xlabel={$N$},
	ylabel={Memory (GB)},
		legend style={
                at={(0.8,0.04)},
               anchor=south,
               legend columns=1,
               cells={anchor=east},
              font=\footnotesize,
               rounded corners=1pt,
               },
		xmode = log,
	    ymode = log,
		]
		
		\addplot[
		color=blue,
		mark=square*,
		] coordinates {
(8000, 64* 0.004938)
(27000, 64* 0.034485)
(64000, 64* 0.092488)
(125000, 64* 0.274634)
(216000, 64* 0.493400)
(343000, 64* 0.819692)
(512000, 64* 1.275320)
(729000, 64* 1.541420)
(1000000, 64* 2.509760)
(1331000, 64* 3.415170)
(1728000, 64* 4.481950)
(2197000, 64* 5.785830)
(2744000, 64* 7.393720)
(3375000, 64* 9.318050)
(4096000, 64* 11.520200)
(4913000, 64* 14.247300)
(5832000, 64* 17.440200)
(6859000, 64* 21.401100)
		};
		\addplot[
		color=red,
		mark=*,
		] coordinates {
(8000, 64* 0.004598)
(27000, 64* 0.030727)
(64000, 64* 0.082146)
(125000, 64* 0.236350)
(216000, 64* 0.427027)
(343000, 64* 0.711266)
(512000, 64* 1.106670)
(729000, 64* 1.673670)
(1000000, 64* 2.830410)
(1331000, 64* 3.860840)
(1728000, 64* 5.051590)
(2197000, 64* 6.551020)
(2744000, 64* 8.336580)
(3375000, 64* 10.420000)
(4096000, 64* 12.850700)
(4913000, 64* 15.771400)
(5832000, 64* 19.204000)
(6859000, 64* 23.049800)
		};
		\addplot[mark=none, red, dashed][
		domain = 8000:7077888,
		] {(pow(10,-5.5)*x*log10(x)};
\addplot[
		color=black,
		mark=triangle*,
		] coordinates {
(8000, 64* 0.003779)
(27000, 64* 0.023287)
(64000, 64* 0.093308)
(125000, 64* 0.252361)
(216000, 64* 0.634886)
(343000, 64* 1.237340)
(512000, 64* 2.357860)
(729000, 64* 4.053440)
		};
\addplot[mark=none, black, dashed][
		domain = 8000:729000,
		] {(pow(10,-7.3)*pow(x,5/3)*log10(x)};
		\legend{H, HODLR3D, $N\log(N)$, HODLR, $N^{5/3}\log(N)$}	
		\end{axis}
\end{tikzpicture}}\quad%
    \subfloat[Assembly time vs N]{\begin{tikzpicture}[every mark/.append style={mark size=1pt}]
	\begin{axis}[	xlabel={$N$},
	ylabel={Initialization Time (s)},
		legend style={
                at={(0.8,0.04)},
               anchor=south,
               legend columns=1,
               cells={anchor=east},
              font=\footnotesize,
               rounded corners=1pt,
               },
		xmode = log,
	    ymode = log,
		]
		
		\addplot[
		color=blue,
		mark=square*,
		] coordinates {

(8000,0.415710)
(27000,3.477850)
(64000,8.277660)
(125000,29.496100)
(216000,49.453900)
(343000,79.519100)
(512000,120.327000)
(729000,165.943000)
(1000000,390.868000)
(1331000,511.824000)
(1728000,620.029000)
(2197000,776.704000)
(2744000,986.599000)
(3375000,1217.880000)
(4096000,1481.300000)
(4913000,1834.890000)
(5832000,2448.040000)
(6859000,2485.020000)
		};
		\addplot[
		color=red,
		mark=*,
		] coordinates {
(8000,0.485883)
(27000,3.733520)
(64000,9.097790)
(125000,29.897500)
(216000,54.459700)
(343000,90.565200)
(512000,132.844000)
(729000,194.015000)
(1000000,391.616000)
(1331000,535.940000)
(1728000,738.490000)
(2197000,832.440000)
(2744000,1122.220000)
(3375000,1398.250000)
(4096000,1734.860000)
(4913000,2204.200000)
(5832000,2696.010000)
(6859000,3122.400000)
		};
		\addplot[mark=none, red, dashed][
		domain = 8000:7077888,
		] {(pow(10,-3.7)*x*log10(x)};
\addplot[
		color=black,
		mark=triangle*,
		] coordinates {
(8000,0.997900)
(27000,11.279300)
(64000,80.588200)
(125000,331.653000)
(216000,1173.040000)
(343000,3108.130000)
(512000,7919.490000)
(729000,18345.600000)
		};
\addplot[mark=none, black, dashed][
		domain = 8000:729000,
		] {(pow(10,-8.6)*pow(x,7/3)*log10(x)};
		\legend{H, HODLR3D, $N\log(N)$, HODLR, $N^{7/3}\log(N)$}
		
		\end{axis}
\end{tikzpicture}}\quad%
    \subfloat[Matrix-vector product time vs N]{\begin{tikzpicture}[every mark/.append style={mark size=1pt}]
	\begin{axis}[xlabel={$N$},
	ylabel={Matrix-Vector product time (s)},
		legend style={
                at={(0.8,0.04)},
               anchor=south,
               legend columns=1,
               cells={anchor=east},
              font=\footnotesize,
               rounded corners=1pt,
               },
		xmode = log,
	    ymode = log,
		]
		
		\addplot[
		color=blue,
		mark=square*,
		] coordinates {
(8000,0.002999)
(27000,0.030109)
(64000,0.045156)
(125000,0.267557)
(216000,0.372911)
(343000,0.465301)
(512000,0.677534)
(729000,0.997064)
(1000000,2.974710)
(1331000,4.362850)
(1728000,4.052130)
(2197000,4.732000)
(2744000,5.841680)
(3375000,7.356550)
(4096000,6.929720)
(4913000,8.587820)
(5832000,12.358000)
(6859000,11.149200)
		};
		\addplot[
		color=red,
		mark=*,
		] coordinates {
(8000,0.002788)
(27000,0.026326)
(64000,0.046711)
(125000,0.244085)
(216000,0.414865)
(343000,0.522523)
(512000,0.652037)
(729000,0.870077)
(1000000,2.607710)
(1331000,3.121390)
(1728000,3.839270)
(2197000,4.662530)
(2744000,5.792950)
(3375000,6.832690)
(4096000,8.439400)
(4913000,9.605090)
(5832000,12.032300)
(6859000,13.236200)
		};
		\addplot[mark=none, red, dashed][
		domain = 8000:7077888,
		] {(pow(10,-5.5)*x*log10(x)};
\addplot[
		color=black,
		mark=triangle*,
		] coordinates {
(8000,0.002446)
(27000,0.018124)
(64000,0.061574)
(125000,0.183393)
(216000,0.469178)
(343000,0.788424)
(512000,1.657480)
(729000,4.099990)
		};
\addplot[mark=none, black, dashed][
		domain = 8000:729000,
		] {(pow(10,-7.5)*pow(x,5/3)*log10(x)};
		\legend{H, HODLR3D, $N\log(N)$, HODLR, $N^{5/3}\log(N)$}
		\end{axis}
\end{tikzpicture}}\quad%
    \subfloat[Relative error vs N]{\begin{tikzpicture}[every mark/.append style={mark size=1pt}]
	\begin{axis}[	xlabel={$N$},
	ylabel={Relative error},
		legend style={
                at={(0.8,0.04)},
               anchor=south,
               legend columns=1,
               cells={anchor=east},
              font=\footnotesize,
               rounded corners=1pt,
               },
		xmode = log,
	    ymode = log,
	    ymin = 1e-16,
	    ymax = 1e-0,
		]
		
		\addplot[
		color=blue,
		mark=square*,
		] coordinates {
(8000,2.557150e-07)
(27000,1.202710e-07)
(64000,1.370560e-07)
(125000,1.230890e-07)
(216000,5.331510e-07)
(343000,1.123290e-07)
(512000,5.729630e-08)
(729000,1.751730e-06)
(1000000,2.027060e-06)
(1331000,2.202360e-06)
(1728000,2.307030e-06)
(2197000,1.690290e-06)
(2744000,1.829160e-06)
(3375000,1.463420e-06)
(4096000,1.213900e-06)
(4913000,1.479610e-06)
(5832000,2.596010e-06)
(6859000,1.748060e-06)
		};
		\addplot[
		color=red,
		mark=*,
		] coordinates {
(8000,2.657270e-07)
(27000,4.404200e-07)
(64000,3.362950e-07)
(125000,4.484270e-07)
(216000,7.017800e-07)
(343000,9.128360e-07)
(512000,4.942570e-07)
(729000,1.300020e-06)
(1000000,1.400770e-06)
(1331000,1.041650e-06)
(1728000,1.125510e-06)
(2197000,1.085220e-06)
(2744000,2.611650e-06)
(3375000,2.357380e-06)
(4096000,1.979970e-06)
(4913000,1.679470e-06)
(5832000,8.381420e-07)
(6859000,1.209460e-06)
		};
\addplot[
		color=black,
		mark=triangle*,
		] coordinates {
(8000,5.823700e-07)
(27000,2.719750e-06)
(64000,7.445260e-08)
(125000,2.823590e-06)
(216000,6.860880e-08)
(343000,1.586930e-06)
(512000,5.050740e-07)
(729000,4.830780e-07)
		};
		\legend{H, HODLR3D, HODLR}
		\end{axis}
\end{tikzpicture}}%
    \caption{Various benchmarks of HODLR3D matrix-vector product in comparison with those of HODLR and $\mathcal{H}$ matrix-vector products for the kernel $\frac{cos\bkt{r}}{r}$ with the relative forward errors of the three algorithms to be of the same order} 
    \label{fig:cosR}
\end{figure}
From~\cref{fig:oneOverR,fig:oneOverR4,fig:cosR}, we observe that by maintaining the relative forward error to be nearly equal, the computational complexity for the matrix-vector product using HODLR3D and $\mathcal{H}$-matrix representation still roughly scales $\mathcal{O}\bkt{N\log\bkt{N}}$, which is not the case with HODLR in 3D.
\section{HODLR3D Initialization in Distributed memory systems}
\label{sec:par_init}
As discussed in~\Cref{sec:parH3D}, we consider the nodes in a particular level of the hierarchical tree as data-independent computational units. For level $l$, where the number of nodes in a level $\bkt{8^l}$ is greater than $n_p$ MPI processes, each MPI process has $\ceil{\dfrac{8^l}{n_p}}$ computational units. For level $l$, where the number of nodes in that level is lesser than $n_p$ MPI processes, each node in level $l$ is shared by $\ceil{\dfrac{n_p}{8^l}}$ MPI processes. The low-rank compression involved with the shared node is performed separately by each MPI process that shares  that node. This is performed to eliminate the communication involved and reduce idle time.~\cref{tab:par_h3d_init} shows the time taken by parallel HODLR3D to initialize the data structure.

\begin{table}[!htbp]
\centering
\caption{Parallel HODLR3D Initialization}
\label{tab:par_h3d_init}
\begin{tabular}{@{}ccccccc@{}}
\toprule
\multirow{3}{*}{\textbf{N}}        & \multirow{3}{*}{$n_p$} & \multicolumn{5}{c}{Time taken for}                                                  \\ \cmidrule(l){3-7} 
                                   &                        & \multicolumn{3}{c}{Form Low Rank Basis} & \multicolumn{2}{c}{Generate Matrix Entry} \\ \cmidrule(l){3-7} 
                                   &                        & Avg          & Max         & Speedup    & Avg                 & Max                 \\ \midrule
\multirow{6}{*}{\textbf{125000}}   & 2                      & 31.04        & 31.07       & 1          & 4.687               & 4.690               \\
                                   & 4                      & 13.62        & 15.09       & 2          & 1.980               & 2.323               \\
                                   & 8                      & 7.46         & 7.72        & 4          & 1.212               & 1.980               \\
                                   & 16                     & 3.80         & 4.43        & 7          & 0.575               & 0.762               \\
                                   & 32                     & 2.10         & 2.86        & 11         & 0.281               & 0.375               \\
                                   & 64                     & 1.18         & 1.70        & 18         & 0.134               & 0.197               \\ \midrule
\multirow{6}{*}{\textbf{1000000}}  & 2                      & 589.65       & 589.94      & 1          & 59.857              & 59.870              \\
                                   & 4                      & 250.02       & 272.61      & 2          & 26.148              & 29.859              \\
                                   & 8                      & 128.09       & 132.75      & 4          & 14.896              & 16.560              \\
                                   & 16                     & 61.22        & 70.36       & 8          & 7.241               & 8.813               \\
                                   & 32                     & 34.31        & 42.90       & 14         & 3.659               & 4.584               \\
                                   & 64                     & 18.79        & 24.92       & 24         & 1.776               & 2.306               \\ \midrule
\multirow{6}{*}{\textbf{3375000}}  & 2                      & 2456.92      & 2482.36     & 1          & 321.421             & 321.837             \\
                                   & 4                      & 989.55       & 1040.38     & 2          & 143.905             & 157.784             \\
                                   & 8                      & 492.24       & 522.79      & 5          & 80.070              & 82.361              \\
                                   & 16                     & 236.63       & 268.47      & 9          & 37.660              & 44.408              \\
                                   & 32                     & 144.93       & 191.63      & 13         & 21.309              & 26.590              \\
                                   & 64                     & 78.82        & 117.63      & 21         & 10.554              & 14.459              \\ \midrule
\multirow{6}{*}{\textbf{8000000}}  & 2                      & 8600.18      & 8620.25     & 1          & 772.062             & 772.723             \\
                                   & 4                      & 3880.94      & 4017.35     & 2          & 346.048             & 383.079             \\
                                   & 8                      & 1948.46      & 2085.96     & 4          & 193.935             & 208.576             \\
                                   & 16                     & 902.95       & 1038.41     & 8          & 87.325              & 100.494             \\
                                   & 32                     & 501.63       & 611.37      & 14         & 49.599              & 58.468              \\
                                   & 64                     & 257.57       & 366.81      & 24         & 26.609              & 39.912              \\ \midrule
\multirow{6}{*}{\textbf{15625000}} & 2                      & 17794.70     & 17821.10    & 1          & 1815.010            & 1815.670            \\
                                   & 4                      & 8058.94      & 8467.18     & 2          & 826.699             & 891.079             \\
                                   & 8                      & 4102.86      & 4348.15     & 4          & 452.395             & 475.912             \\
                                   & 16                     & 1870.53      & 2137.44     & 8          & 212.278             & 242.730             \\
                                   & 32                     & 1027.04      & 1284.76     & 14         & 120.396             & 143.222             \\
                                   & 64                     & 555.55       & 801.46      & 22         & 60.863              & 84.852              \\ \bottomrule
\end{tabular}
\end{table}

\end{appendices}

\bibliography{references}


\begin{thebibliography}{36}
\ifx \bisbn   \undefined \def \bisbn  #1{ISBN #1}\fi
\ifx \binits  \undefined \def \binits#1{#1}\fi
\ifx \bauthor  \undefined \def \bauthor#1{#1}\fi
\ifx \batitle  \undefined \def \batitle#1{#1}\fi
\ifx \bjtitle  \undefined \def \bjtitle#1{#1}\fi
\ifx \bvolume  \undefined \def \bvolume#1{\textbf{#1}}\fi
\ifx \byear  \undefined \def \byear#1{#1}\fi
\ifx \bissue  \undefined \def \bissue#1{#1}\fi
\ifx \bfpage  \undefined \def \bfpage#1{#1}\fi
\ifx \blpage  \undefined \def \blpage #1{#1}\fi
\ifx \burl  \undefined \def \burl#1{\textsf{#1}}\fi
\ifx \doiurl  \undefined \def \doiurl#1{\url{https://doi.org/#1}}\fi
\ifx \betal  \undefined \def \betal{\textit{et al.}}\fi
\ifx \binstitute  \undefined \def \binstitute#1{#1}\fi
\ifx \binstitutionaled  \undefined \def \binstitutionaled#1{#1}\fi
\ifx \bctitle  \undefined \def \bctitle#1{#1}\fi
\ifx \beditor  \undefined \def \beditor#1{#1}\fi
\ifx \bpublisher  \undefined \def \bpublisher#1{#1}\fi
\ifx \bbtitle  \undefined \def \bbtitle#1{#1}\fi
\ifx \bedition  \undefined \def \bedition#1{#1}\fi
\ifx \bseriesno  \undefined \def \bseriesno#1{#1}\fi
\ifx \blocation  \undefined \def \blocation#1{#1}\fi
\ifx \bsertitle  \undefined \def \bsertitle#1{#1}\fi
\ifx \bsnm \undefined \def \bsnm#1{#1}\fi
\ifx \bsuffix \undefined \def \bsuffix#1{#1}\fi
\ifx \bparticle \undefined \def \bparticle#1{#1}\fi
\ifx \barticle \undefined \def \barticle#1{#1}\fi
\bibcommenthead
\ifx \bconfdate \undefined \def \bconfdate #1{#1}\fi
\ifx \botherref \undefined \def \botherref #1{#1}\fi
\ifx \url \undefined \def \url#1{\textsf{#1}}\fi
\ifx \bchapter \undefined \def \bchapter#1{#1}\fi
\ifx \bbook \undefined \def \bbook#1{#1}\fi
\ifx \bcomment \undefined \def \bcomment#1{#1}\fi
\ifx \oauthor \undefined \def \oauthor#1{#1}\fi
\ifx \citeauthoryear \undefined \def \citeauthoryear#1{#1}\fi
\ifx \endbibitem  \undefined \def \endbibitem {}\fi
\ifx \bconflocation  \undefined \def \bconflocation#1{#1}\fi
\ifx \arxivurl  \undefined \def \arxivurl#1{\textsf{#1}}\fi
\csname PreBibitemsHook\endcsname

\bibitem[\protect\citeauthoryear{Gray and Moore}{2000}]{gray2000n}
\begin{botherref}
\oauthor{\bsnm{Gray}, \binits{A.}},
\oauthor{\bsnm{Moore}, \binits{A.}}:
N-body'problems in statistical learning.
Advances in neural information processing systems
\textbf{13}
(2000)
\end{botherref}
\endbibitem

\bibitem[\protect\citeauthoryear{Litvinenko
  et~al.}{2019}]{litvinenko2019likelihood}
\begin{barticle}
\bauthor{\bsnm{Litvinenko}, \binits{A.}},
\bauthor{\bsnm{Sun}, \binits{Y.}},
\bauthor{\bsnm{Genton}, \binits{M.G.}},
\bauthor{\bsnm{Keyes}, \binits{D.E.}}:
\batitle{Likelihood approximation with hierarchical matrices for large spatial
  datasets}.
\bjtitle{Computational Statistics \& Data Analysis}
\bvolume{137},
\bfpage{115}--\blpage{132}
(\byear{2019})
\end{barticle}
\endbibitem

\bibitem[\protect\citeauthoryear{Coulier and
  Darve}{2016}]{coulier2016efficient}
\begin{barticle}
\bauthor{\bsnm{Coulier}, \binits{P.}},
\bauthor{\bsnm{Darve}, \binits{E.}}:
\batitle{Efficient mesh deformation based on radial basis function
  interpolation by means of the inverse fast multipole method}.
\bjtitle{Computer Methods in Applied Mechanics and Engineering}
\bvolume{308},
\bfpage{286}--\blpage{309}
(\byear{2016})
\end{barticle}
\endbibitem

\bibitem[\protect\citeauthoryear{Gumerov and
  Duraiswami}{2007}]{gumerov2007fast}
\begin{barticle}
\bauthor{\bsnm{Gumerov}, \binits{N.A.}},
\bauthor{\bsnm{Duraiswami}, \binits{R.}}:
\batitle{Fast radial basis function interpolation via preconditioned krylov
  iteration}.
\bjtitle{SIAM Journal on Scientific Computing}
\bvolume{29}(\bissue{5}),
\bfpage{1876}--\blpage{1899}
(\byear{2007})
\end{barticle}
\endbibitem

\bibitem[\protect\citeauthoryear{Hackbusch}{1999}]{hackbusch1999sparse}
\begin{barticle}
\bauthor{\bsnm{Hackbusch}, \binits{W.}}:
\batitle{A sparse matrix arithmetic based on $\mathcal{H}$-matrices. part i:
  Introduction to $\mathcal{H}$-matrices}.
\bjtitle{Computing}
\bvolume{62}(\bissue{2}),
\bfpage{89}--\blpage{108}
(\byear{1999})
\end{barticle}
\endbibitem

\bibitem[\protect\citeauthoryear{Grasedyck and
  Hackbusch}{2003}]{grasedyck2003construction}
\begin{barticle}
\bauthor{\bsnm{Grasedyck}, \binits{L.}},
\bauthor{\bsnm{Hackbusch}, \binits{W.}}:
\batitle{Construction and arithmetics of h-matrices}.
\bjtitle{Computing}
\bvolume{70}(\bissue{4}),
\bfpage{295}--\blpage{334}
(\byear{2003})
\end{barticle}
\endbibitem

\bibitem[\protect\citeauthoryear{Kandappan et~al.}{2022}]{H2DKandappan}
\begin{botherref}
\oauthor{\bsnm{Kandappan}, \binits{V.A.}},
\oauthor{\bsnm{Gujjula}, \binits{V.}},
\oauthor{\bsnm{Ambikasaran}, \binits{S.}}:
Hodlr2d: A new class of hierarchical matrices.
arXiv preprint arXiv:2204.05536
(2022)
\end{botherref}
\endbibitem

\bibitem[\protect\citeauthoryear{Barnes and Hut}{1986}]{barnes1986hierarchical}
\begin{barticle}
\bauthor{\bsnm{Barnes}, \binits{J.}},
\bauthor{\bsnm{Hut}, \binits{P.}}:
\batitle{A hierarchical o (n log n) force-calculation algorithm}.
\bjtitle{nature}
\bvolume{324}(\bissue{6096}),
\bfpage{446}--\blpage{449}
(\byear{1986})
\end{barticle}
\endbibitem

\bibitem[\protect\citeauthoryear{Greengard and
  Rokhlin}{1987}]{greengard1987fast}
\begin{barticle}
\bauthor{\bsnm{Greengard}, \binits{L.}},
\bauthor{\bsnm{Rokhlin}, \binits{V.}}:
\batitle{A fast algorithm for particle simulations}.
\bjtitle{Journal of computational physics}
\bvolume{73}(\bissue{2}),
\bfpage{325}--\blpage{348}
(\byear{1987})
\end{barticle}
\endbibitem

\bibitem[\protect\citeauthoryear{Greengard}{1988}]{greengard1988rapid}
\begin{bbook}
\bauthor{\bsnm{Greengard}, \binits{L.}}:
\bbtitle{The Rapid Evaluation of Potential Fields in Particle Systems}.
\bpublisher{MIT press}, \blocation{???}
(\byear{1988})
\end{bbook}
\endbibitem

\bibitem[\protect\citeauthoryear{Greengard and
  Rokhlin}{1997}]{greengard1997new}
\begin{barticle}
\bauthor{\bsnm{Greengard}, \binits{L.}},
\bauthor{\bsnm{Rokhlin}, \binits{V.}}:
\batitle{A new version of the fast multipole method for the laplace equation in
  three dimensions}.
\bjtitle{Acta numerica}
\bvolume{6},
\bfpage{229}--\blpage{269}
(\byear{1997})
\end{barticle}
\endbibitem

\bibitem[\protect\citeauthoryear{Ambikasaran}{2013}]{sa_thesis}
\begin{botherref}
\oauthor{\bsnm{Ambikasaran}, \binits{S.}}:
Fast algorithms for dense numerical linear algebra and applications.
PhD thesis,
Stanford University
(2013)
\end{botherref}
\endbibitem

\bibitem[\protect\citeauthoryear{Ambikasaran and Darve}{2013}]{SA_FDS_2013}
\begin{barticle}
\bauthor{\bsnm{Ambikasaran}, \binits{S.}},
\bauthor{\bsnm{Darve}, \binits{E.}}:
\batitle{An $\mathcal{O} (n \log n)$ - fast direct solver for partial
  hierarchically semi-separable matrices}.
\bjtitle{Journal of Scientific Computing}
\bvolume{57}(\bissue{3}),
\bfpage{477}--\blpage{501}
(\byear{2013})
\doiurl{10.1007/s10915-013-9714-z}
\end{barticle}
\endbibitem

\bibitem[\protect\citeauthoryear{Chandrasekaran
  et~al.}{2005}]{chandrasekaran2005some}
\begin{barticle}
\bauthor{\bsnm{Chandrasekaran}, \binits{S.}},
\bauthor{\bsnm{Dewilde}, \binits{P.}},
\bauthor{\bsnm{Gu}, \binits{M.}},
\bauthor{\bsnm{Pals}, \binits{T.}},
\bauthor{\bsnm{Sun}, \binits{X.}},
\bauthor{\bsnm{Veen}, \binits{A.-J.}},
\bauthor{\bsnm{White}, \binits{D.}}:
\batitle{Some fast algorithms for sequentially semiseparable representations}.
\bjtitle{SIAM Journal on Matrix Analysis and Applications}
\bvolume{27}(\bissue{2}),
\bfpage{341}--\blpage{364}
(\byear{2005})
\end{barticle}
\endbibitem

\bibitem[\protect\citeauthoryear{Vandebril
  et~al.}{2005a}]{vandebril2005bibliography}
\begin{barticle}
\bauthor{\bsnm{Vandebril}, \binits{R.}},
\bauthor{\bsnm{Barel}, \binits{M.V.}},
\bauthor{\bsnm{Golub}, \binits{G.}},
\bauthor{\bsnm{Mastronardi}, \binits{N.}}:
\batitle{A bibliography on semiseparable matrices}.
\bjtitle{Calcolo}
\bvolume{42}(\bissue{3}),
\bfpage{249}--\blpage{270}
(\byear{2005})
\end{barticle}
\endbibitem

\bibitem[\protect\citeauthoryear{Vandebril et~al.}{2005b}]{vandebril2005note}
\begin{barticle}
\bauthor{\bsnm{Vandebril}, \binits{R.}},
\bauthor{\bsnm{Van~Barel}, \binits{M.}},
\bauthor{\bsnm{Mastronardi}, \binits{N.}}:
\batitle{A note on the representation and definition of semiseparable
  matrices}.
\bjtitle{Numerical Linear Algebra with Applications}
\bvolume{12}(\bissue{8}),
\bfpage{839}--\blpage{858}
(\byear{2005})
\end{barticle}
\endbibitem

\bibitem[\protect\citeauthoryear{B{\"o}rm et~al.}{2003a}]{borm2003hierarchical}
\begin{barticle}
\bauthor{\bsnm{B{\"o}rm}, \binits{S.}},
\bauthor{\bsnm{Grasedyck}, \binits{L.}},
\bauthor{\bsnm{Hackbusch}, \binits{W.}}:
\batitle{Hierarchical matrices}.
\bjtitle{Lecture notes}
\bvolume{21},
\bfpage{2003}
(\byear{2003})
\end{barticle}
\endbibitem

\bibitem[\protect\citeauthoryear{B{\"o}rm et~al.}{2003b}]{borm2003introduction}
\begin{barticle}
\bauthor{\bsnm{B{\"o}rm}, \binits{S.}},
\bauthor{\bsnm{Grasedyck}, \binits{L.}},
\bauthor{\bsnm{Hackbusch}, \binits{W.}}:
\batitle{Introduction to hierarchical matrices with applications}.
\bjtitle{Engineering analysis with boundary elements}
\bvolume{27}(\bissue{5}),
\bfpage{405}--\blpage{422}
(\byear{2003})
\end{barticle}
\endbibitem

\bibitem[\protect\citeauthoryear{Hackbusch}{2015}]{hackbusch2015hierarchical}
\begin{bbook}
\bauthor{\bsnm{Hackbusch}, \binits{W.}}:
\bbtitle{Hierarchical Matrices: Algorithms and Analysis}
vol. \bseriesno{49}.
\bpublisher{Springer}, \blocation{???}
(\byear{2015})
\end{bbook}
\endbibitem

\bibitem[\protect\citeauthoryear{Yokota et~al.}{2015}]{yokota2015fast}
\begin{bchapter}
\bauthor{\bsnm{Yokota}, \binits{R.}},
\bauthor{\bsnm{Ibeid}, \binits{H.}},
\bauthor{\bsnm{Keyes}, \binits{D.}}:
\bctitle{Fast multipole method as a matrix-free hierarchical low-rank
  approximation}.
In: \bbtitle{International Workshop on Eigenvalue Problems: Algorithms,
  Software and Applications in Petascale Computing},
pp. \bfpage{267}--\blpage{286}
(\byear{2015}).
\bcomment{Springer}
\end{bchapter}
\endbibitem

\bibitem[\protect\citeauthoryear{Amestoy et~al.}{2015}]{amestoy2015improving}
\begin{barticle}
\bauthor{\bsnm{Amestoy}, \binits{P.}},
\bauthor{\bsnm{Ashcraft}, \binits{C.}},
\bauthor{\bsnm{Boiteau}, \binits{O.}},
\bauthor{\bsnm{Buttari}, \binits{A.}},
\bauthor{\bsnm{l'Excellent}, \binits{J.-Y.}},
\bauthor{\bsnm{Weisbecker}, \binits{C.}}:
\batitle{Improving multifrontal methods by means of block low-rank
  representations}.
\bjtitle{SIAM Journal on Scientific Computing}
\bvolume{37}(\bissue{3}),
\bfpage{1451}--\blpage{1474}
(\byear{2015})
\end{barticle}
\endbibitem

\bibitem[\protect\citeauthoryear{Amestoy et~al.}{2017}]{amestoy2017complexity}
\begin{barticle}
\bauthor{\bsnm{Amestoy}, \binits{P.}},
\bauthor{\bsnm{Buttari}, \binits{A.}},
\bauthor{\bsnm{l'Excellent}, \binits{J.-Y.}},
\bauthor{\bsnm{Mary}, \binits{T.}}:
\batitle{On the complexity of the block low-rank multifrontal factorization}.
\bjtitle{SIAM Journal on Scientific Computing}
\bvolume{39}(\bissue{4}),
\bfpage{1710}--\blpage{1740}
(\byear{2017})
\end{barticle}
\endbibitem

\bibitem[\protect\citeauthoryear{Khan et~al.}{2022}]{khan2022numerical}
\begin{botherref}
\oauthor{\bsnm{Khan}, \binits{R.}},
\oauthor{\bsnm{Kandappan}, \binits{V.}},
\oauthor{\bsnm{Ambikasaran}, \binits{S.}}:
Numerical rank of singular kernel functions.
arXiv preprint arXiv:2209.05819
(2022)
\end{botherref}
\endbibitem

\bibitem[\protect\citeauthoryear{Hackbusch
  et~al.}{2004}]{hackbusch2004hierarchical}
\begin{barticle}
\bauthor{\bsnm{Hackbusch}, \binits{W.}},
\bauthor{\bsnm{Khoromskij}, \binits{B.N.}},
\bauthor{\bsnm{Kriemann}, \binits{R.}}:
\batitle{Hierarchical matrices based on a weak admissibility criterion}.
\bjtitle{Computing}
\bvolume{73}(\bissue{3}),
\bfpage{207}--\blpage{243}
(\byear{2004})
\end{barticle}
\endbibitem

\bibitem[\protect\citeauthoryear{Beatson and
  Greengard}{1997}]{beatson1997short}
\begin{barticle}
\bauthor{\bsnm{Beatson}, \binits{R.}},
\bauthor{\bsnm{Greengard}, \binits{L.}}:
\batitle{A short course on fast multipole methods}.
\bjtitle{Wavelets, multilevel methods and elliptic PDEs}
\bvolume{1},
\bfpage{1}--\blpage{37}
(\byear{1997})
\end{barticle}
\endbibitem

\bibitem[\protect\citeauthoryear{Bebendorf and
  Rjasanow}{2003}]{bebendorf2003adaptive}
\begin{barticle}
\bauthor{\bsnm{Bebendorf}, \binits{M.}},
\bauthor{\bsnm{Rjasanow}, \binits{S.}}:
\batitle{Adaptive low-rank approximation of collocation matrices}.
\bjtitle{Computing}
\bvolume{70}(\bissue{1}),
\bfpage{1}--\blpage{24}
(\byear{2003})
\end{barticle}
\endbibitem

\bibitem[\protect\citeauthoryear{Zhao et~al.}{2005}]{zhao2005adaptive}
\begin{barticle}
\bauthor{\bsnm{Zhao}, \binits{K.}},
\bauthor{\bsnm{Vouvakis}, \binits{M.N.}},
\bauthor{\bsnm{Lee}, \binits{J.-F.}}:
\batitle{The adaptive cross approximation algorithm for accelerated method of
  moments computations of emc problems}.
\bjtitle{IEEE transactions on electromagnetic compatibility}
\bvolume{47}(\bissue{4}),
\bfpage{763}--\blpage{773}
(\byear{2005})
\end{barticle}
\endbibitem

\bibitem[\protect\citeauthoryear{Tyrtyshnikov}{2000}]{tyrtyshnikov2000incomplete}
\begin{barticle}
\bauthor{\bsnm{Tyrtyshnikov}, \binits{E.}}:
\batitle{Incomplete cross approximation in the mosaic-skeleton method}.
\bjtitle{Computing}
\bvolume{64}(\bissue{4}),
\bfpage{367}--\blpage{380}
(\byear{2000})
\end{barticle}
\endbibitem

\bibitem[\protect\citeauthoryear{Bebendorf}{2000}]{bebendorf2000approximation}
\begin{barticle}
\bauthor{\bsnm{Bebendorf}, \binits{M.}}:
\batitle{Approximation of boundary element matrices}.
\bjtitle{Numerische Mathematik}
\bvolume{86}(\bissue{4}),
\bfpage{565}--\blpage{589}
(\byear{2000})
\end{barticle}
\endbibitem

\bibitem[\protect\citeauthoryear{Bebendorf and
  Kunis}{2009}]{bebendorf2009recompression}
\begin{botherref}
\oauthor{\bsnm{Bebendorf}, \binits{M.}},
\oauthor{\bsnm{Kunis}, \binits{S.}}:
Recompression techniques for adaptive cross approximation.
The Journal of Integral Equations and Applications,
331--357
(2009)
\end{botherref}
\endbibitem

\bibitem[\protect\citeauthoryear{Barrett et~al.}{1994}]{barrett1994templates}
\begin{bbook}
\bauthor{\bsnm{Barrett}, \binits{R.}},
\bauthor{\bsnm{Berry}, \binits{M.}},
\bauthor{\bsnm{Chan}, \binits{T.F.}},
\bauthor{\bsnm{Demmel}, \binits{J.}},
\bauthor{\bsnm{Donato}, \binits{J.}},
\bauthor{\bsnm{Dongarra}, \binits{J.}},
\bauthor{\bsnm{Eijkhout}, \binits{V.}},
\bauthor{\bsnm{Pozo}, \binits{R.}},
\bauthor{\bsnm{Romine}, \binits{C.}},
\bauthor{\bsnm{Vorst}, \binits{H.}}:
\bbtitle{Templates for the Solution of Linear Systems: Building Blocks for
  Iterative Methods}.
\bpublisher{SIAM}, \blocation{???}
(\byear{1994})
\end{bbook}
\endbibitem

\bibitem[\protect\citeauthoryear{Saad and Schultz}{1986}]{saad1986gmres}
\begin{barticle}
\bauthor{\bsnm{Saad}, \binits{Y.}},
\bauthor{\bsnm{Schultz}, \binits{M.H.}}:
\batitle{Gmres: A generalized minimal residual algorithm for solving
  nonsymmetric linear systems}.
\bjtitle{SIAM Journal on scientific and statistical computing}
\bvolume{7}(\bissue{3}),
\bfpage{856}--\blpage{869}
(\byear{1986})
\end{barticle}
\endbibitem

\bibitem[\protect\citeauthoryear{Izadi}{2012}]{izadi2012hierarchical}
\begin{botherref}
\oauthor{\bsnm{Izadi}, \binits{M.}}:
Hierarchical matrix techniques on massively parallel computers.
Thesis
(2012)
\end{botherref}
\endbibitem

\bibitem[\protect\citeauthoryear{Li et~al.}{2020}]{li2020distributed}
\begin{botherref}
\oauthor{\bsnm{Li}, \binits{Y.}},
\oauthor{\bsnm{Poulson}, \binits{J.}},
\oauthor{\bsnm{Ying}, \binits{L.}}:
Distributed-memory $\mathcal {H} $-matrix algebra i: Data distribution and
  matrix-vector multiplication.
arXiv preprint arXiv:2008.12441
(2020)
\end{botherref}
\endbibitem

\bibitem[\protect\citeauthoryear{Ambikasaran and
  Darve}{2014}]{ambikasaran2014inverse}
\begin{botherref}
\oauthor{\bsnm{Ambikasaran}, \binits{S.}},
\oauthor{\bsnm{Darve}, \binits{E.}}:
The inverse fast multipole method.
arXiv preprint arXiv:1407.1572
(2014)
\end{botherref}
\endbibitem

\bibitem[\protect\citeauthoryear{Gujjula and
  Ambikasaran}{2023}]{gujjula2023algebraic}
\begin{botherref}
\oauthor{\bsnm{Gujjula}, \binits{V.}},
\oauthor{\bsnm{Ambikasaran}, \binits{S.}}:
Algebraic inverse fast multipole method: A fast direct solver that is better
  than hodlr based fast direct solver.
arXiv preprint arXiv:2301.12704
(2023)
\end{botherref}
\endbibitem

\end{thebibliography}

\end{document}